\documentclass{article}

\usepackage{amsmath, amssymb} %,amsthm,txfonts
\usepackage{enumerate}
\usepackage{pgf}
\usepackage{tikz}

\usepackage{subcaption}
\usepackage{float}
\floatstyle{boxed}
\restylefloat{figure}
\usepackage{placeins}

\newcommand{\bc}{\begin{C}}
\newcommand{\ec}{\end{C}}
\newcommand{\be}{\begin{equation}}
\newcommand{\ee}{\end{equation}}

\newcommand{\claim}{\begin{Cl}}
\newcommand{\eclaim}{\end{Cl}}
\newcommand{\nb}{\begin{Prop}}
\newcommand{\nbe}{\end{Prop}}
\newcommand{\bl}{\begin{LE}}
\newcommand{\el}{\end{LE}}

\newtheorem{Cl}{Claim}

\newcommand{\bd}{\begin{Def}}
\newcommand{\ed}{\end{Def}}
\newcommand{\bt}{\begin{Th}}
\newcommand{\et}{\end{Th}}

\newtheorem{Th}{Theorem}[section]
\newtheorem{LE}[Th]{Lemma}
\newtheorem{C}[Th]{Corollary}

\newtheorem{Def}[Th]{Definition}
\newcommand{\bp}{\begin{Prop}}
\newcommand{\ep}{\end{Prop}}

\newtheorem{Prop}[Th]{Proposition}

\DeclareMathSymbol{\mlq}{\mathord}{operators}{``}
\DeclareMathSymbol{\mrq}{\mathord}{operators}{`'}

\renewenvironment{itemize}
  {\begin{olditemize}
     \setlength{\parskip}{0pt}}
  {\end{olditemize}}

\begin{document}
 \title{On embedding Lambek calculus into commutative categorial grammars
}
\author{Sergey Slavnov
\\  National Research University Higher School of Economics
\\ sslavnov@yandex.ru\\} \maketitle

\begin{abstract}
We consider {\it tensor grammars}, which are an example of ``commutative'' grammars, based on the classical  (rather than intuitionistic) linear logic. They can be seen as a surface representation of abstract categorial grammars ACG in the sense that derivations of ACG translate to derivations of tensor grammars and this translation is isomorphic on the level of string languages. The basic ingredient are {\it tensor terms}, which can be seen as encoding and generalizing proof-nets. Using tensor terms makes the syntax  extremely simple and a direct geometric meaning becomes transparent.
Then we address the problem of encoding noncommutative operations in our setting. This turns out possible after enriching the system with new unary operators. The resulting system allows representing both ACG and Lambek grammars as conservative fragments, while the formalism remains, as it seems to us, rather simple and intuitive.
\end{abstract}

\section{Introduction}
\subsection{Commutative and noncommutative grammars}
The best known approach to categorial grammars is based on noncommutative variants of linear logic, most notably, on Lambek calculus ({\bf LC}) \cite{Lambek} and its variations/extensions, such as {\it non-associative Lambek calculus} \cite{Lambek_non-ass}, various {\it mixed multimodal} systems \cite{Moortgat}, {\it displacement calculus} \cite{Morrill_Displacement} etc.

{\it Abstract categorial grammars} (ACG) \cite{deGroote}, as well as their close relatives {\it $\lambda$-grammars}  \cite{Muskens} or {\it linear grammars}  \cite{PollardMichalicek}, arise from an alternative or, rather, complementary approach, based on the ordinary implicational linear logic and linear $\lambda$-calculus. They  can be called ``commutative'' in contrast to the ``noncommutative'' Lambek grammars. Commutative grammars are remarkably flexible and expressive (sometimes even too expressive for effective parsing \cite{YoshinakaKanazawa}). On the other hand they are also remarkably simple because of the much more familiar and intuitive  underlying logic.

However, as far as natural language modeling is concerned, it turns out that, in many situations, commutative grammars, such as ACG,  behave very poorly compared to noncommutative variants.  In fact, it has  even been argued that ACG are {\it descriptively inadequate} \cite{Moot_inadequacy}. Simple and striking instances of this inadequacy arise, for example, when linguistic coordination is considered.

The reason is that, as an analysis shows,
 ``commutative'' types of ACG and its relatives are too coarse to distinguish actual linguistic categories.
  Thus, if we want to model important linguistic phenomena in the commutative setting, we need somehow to enrich the formalism  with a finer structure of {\it subtypes} corresponding to the ``noncommutative'' types of {\bf LC}.

  One solution to this problem, was proposed in \cite{Worth} (but see also \cite{Kanazawa_ACG_vs_Labek}), where an explicit subtyping mechanism was added to the system. Unfortunately, this results in a rather impressive complication of the formalism (as it seems to us). Another proposed direction is, simply, to enrich a commutative system with explicit noncommutative constructions. (This suggests a comparison with the quite long known Abrusci-Ruet logic \cite{AbrusciRuet}.) {\it Hybrid type logical categorial grammars} (HTLCG) \cite{KubotaLevine} have three kinds of implication on the level of types (one commutative implication of linear logic and two noncommutative slashes of {\bf LC}) and two kinds of application on the level of terms, the usual application of $\lambda$-terms and an additional  operation of concatenation.
  Both approaches, of \cite{Worth} and of \cite{KubotaLevine}, led to, at least partially, successful developments.
  An apparent drawback though, as it seems to us, is that the attractive simplicity of ACG gets somewhat lost.

An interesting perspective comes from considering {\it first order logic} \cite{MootPiazza}, \cite{Moot_comparing}, \cite{Moot_inadequacy}. It turns out that different grammatical formalisms including Lambek grammars, ACG and HTLCG   can be faithfully  represented as fragments of {\it first order multiplicative intuitionistic linear logic} ({\bf MILL1}). This suggests another approach to combining commutative and noncommutative features, as well as provides some common ground on which different systems can be compared.

 \subsection{Content of this work}
 {\it Tensor grammars} of this work are an elaboration of the so called {\it linear logic grammars} (LLG) introduced in \cite{Slavnov_cowordisms}.

  LLG are another example of  commutative grammars, based on the {\it classical}, rather than intuitionistic, multiplicative linear logic ({\bf MLL}). They were defined in terms of certain bipartite graphs (generalizing {\bf MLL} {\it proof-nets}) with string-labeled edges.

LLG   (as well as tensor grammars of this work) can be seen as a surface representation of ACG. Derivations of ACG translate to derivations of tensor grammars and this translation is isomorphic on the  level of string languages. On the logical side, this is, simply, a reflection of the fact that implicational linear logic is a conservative fragment of classical {\bf MLL} and linear $\lambda$-terms can be represented as proof-nets. An advantage of this representation, as it seems to us, is that the syntax becomes extremely simple and a direct geometric meaning is transparent.

  In this work we introduce an in-line notation for edge-labeled graphs of LLG and reformulate the system in these new terms. Then we address the problem of encoding noncommutative operations of {\bf LC} in our setting. This turns out possible after enriching the system with new unary operators,  which results in {\it extended tensor grammars}.

\subsubsection{Tensor grammars}
{\it Tensor terms} are, basically, tuples of words (written multiplicatively, as products) with labeled endpoints. We write words in square brackets and represent the endpoints as lower and upper indices,  lower indices standing for left endpoints and upper indices for right ones. Thus, tensor terms can have the form
$$[\mbox{a}]_i^j,\quad [\mbox{a}]_i^j\cdot[\mbox{b}]_k^l\cdot[\mbox{c}]_r^s,\quad [\epsilon]_i^j $$
(where $\epsilon$ stands for the empty word) and so on.
The index notation is taken directly from usual tensor algebra.

An index in a term can be repeated at most twice, once as an upper, and once as a lower one. A repeated (i.e. {\it bound}) index means that the corresponding words are concatenated along matching endpoints. For example, we have the term equality $$[\mbox{a}]_i^j\cdot[\mbox{b}]_k^l\cdot[\mbox{c}]_j^k=[\mbox{acb}]_i^l$$
(the product is commutative).

{\it Tensor term calculus} ({\bf TTC}) equips tensor terms with types. Tensor types are {\bf MLL} formulas decorated with indices, which should match indices in corresponding terms, with the convention that upper indices in types match lower indices in terms and vice versa. A  {\it tensor typing judgement} looks, for example, as the following:
$$[\mbox{a}]_i^j\cdot[\mbox{b}]_k^l\cdot[\mbox{c}]_r^s\vdash (A_j^k\otimes B_l)\wp C^{ir}_s.$$
Typing rules, of course, are rules of {\bf MLL} decorated with terms and indices.

A {\it tensor grammar} is defined then by a set of {\it axioms}, which are tensor typing judgements as above, and the {\it sentence type} $S_i^j$ with exactly one upper and one lower index. Tensor terms of the sentence type are single words and they constitute the language of the grammar.

Of course, this  is  an oversimplified inaccurate sketch rather than  a consistent formalism. But we think that it is easy to believe that technical details can be worked out, as well as that the resulting formalism indeed provides a representation for ACG. Unfortunately it does not allow simple representation of noncommutative operations of {\bf LC}.

Before approaching noncommutativity we note that indices in our formalism  can be thought of as first order variables, which  suggests a comparison with  {\bf MILL1} and the work in \cite{MootPiazza}, \cite{Moot_comparing}. The system of  {\bf MILL1} has an extra degree of freedom because of binding operators in formulas, i.e. quantifiers. It is thanks to quantifiers  that   representation of noncommutative systems becomes possible.

This suggests that we need to extend tensor grammars with {\it index binding operators} on types.

\subsubsection{Extended tensor types}
Tensor representation makes very transparent how the ``non-commutative'' types of {\bf LC} look inside ``coarse'' implicational   types of commutative grammars.

Atomic types of {\bf LC} correspond to  tensor types with one upper and one lower index. If   $A=A^i_j$ and $B=B^k_l$ are two such types, then, by the standard definition of linear implication in {\bf MLL}, the implicational type $A\multimap B$ is the type $\overline{A}_i^j\wp B^k_l$ (where we denote linear negation as a bar). Thus, elements of the type $A\multimap B$ are tensor terms with four indices.

We can single out two important {\it subtypes} of the tensor type $A\multimap B$. The first subtype consists of terms of the form
\be\label{lambek backslash}
[u]^i_k\cdot[\epsilon]^l_j,
\ee
and the second one, of terms of the form
\be\label{lambek/}
[u]^l_j\cdot[\epsilon]^i_k.
\ee

It is easily computed that terms of  form (\ref{lambek backslash}) act  on elements of $A_i^j$ by multiplication (concatenation) on the left, and elements of form (\ref{lambek/}), by multiplication on the right.
Indeed, for a term $[v]_i^j$ of type $A^i_j$, we have
\be\label{subtypes actiom}
[v]_i^j\cdot([u]^i_k\cdot[\epsilon]^l_j)=[uv]_k^l,\quad
[v]_i^j\cdot([u]^l_j\cdot[\epsilon]^i_k)=[vu]_k^l.
\ee
It follows that the subtype defined by (\ref{lambek backslash}) corresponds to the Lambek type $B/A$, while the subtype defined by (\ref{lambek/}), to the Lambek type $A\backslash B$.

Following this insight, we  emulate noncommutative implications by means of a new type constructor $\nabla$, which {\it  binds indices in types}.

If $A$ is a tensor type with a {\it free} upper index $i$ and free lower index $j$,  we define the new type $\nabla_i^jA$, in which the indices $i$ and $j$ are no longer free.
The rule for introducing $\nabla$ is
$$\frac{t\cdot [\epsilon]_i^j\vdash\Gamma, A^{I_1iI_2}_{J_1jJ_2}}
{t\vdash\Gamma, \nabla_i^j(A^{I_1iI_2}_{J_1jJ_2})},$$
and noncommutative implications are encoded as
$$(A\backslash B)^j_l=\nabla^i_k(\overline{A}_i^j\wp B^k_l),\quad (B/A)^j_l=\nabla^i_k(B_i^j\wp \overline{A}^k_l),$$
as formulas (\ref{lambek backslash}), (\ref{lambek/}), (\ref{subtypes actiom}) suggest.

In fact, since we work in the setting of classical linear logic, we  also have to introduce the {\it dual} operator $\triangle$ of $\nabla$. Remarkably, this second binding operator serves to encode the  product of {\bf LC}.

Along these lines we develop  {\it extended tensor type calculus} ({\bf ETTC}) of tensor terms and define {\it extended tensor grammars}. Both ACG and Lambek grammars embed in {\bf ETTC}. We would expect that HTLCG and displacement grammars do embed as well, but this requires a proof.

Again, it seems to us that the formalism of {\bf ETTC} is rather simple and intuitive, most notably because of its transparent geometric meaning, and  this makes our work interesting.

\subsection{Comparing with first order logic}
It remains to be understood how does {\bf ETTC} of this work compare  with the first order formalism of \cite{MootPiazza}, \cite{Moot_comparing},  and to what degree (and in what sense) do the two formalisms overlap. Indices in tensor types certainly can be read as first order variables, but the binding operators of {\bf ETTC} do not seem to behave exactly as quantifiers. Comparing and understanding relationships of these formalisms seems to us an  important question, and we leave it for future research.

An attractive feature of {\bf ETTC} is that its typing judgements do directly represent  surface elements of language (i.e. tuples of strings), while a general {\bf MILL1} formula, sequent  or derivation   has no intrinsic linguistic interpretation whatsoever. Representations of logical grammars  actually use only  carefully selected fragments of  {\bf MILL1}  and assign interpretations only to some  specific derivations (or, rather, derivable sequents).

Typically, when deriving the translation of an {\bf LC} sequent  in {\bf MILL1}, one might have to use at intermediate steps sequents and derivations which are not translations of anything and have no linguistic interpretation at all (as it seems).
%The fragment of {\bf MILL1} that is actually used for grammars has no sequent calculus formulation.
This contrasts with {\bf ETTC} where each step in the derivation corresponds to a concrete operation with strings (which, moreover, can be conveniently visualized in the pictorial setting of labeled graphs). Thus formal language generation is decomposed into elementary steps.
From a proof-theoretic perspective, {\bf ETTC} comes with a solid {\it denotational  semantics} and this defines grammars.
{\bf MILL1} grammars, on the contrary, simply  select  a sparse class of sequents, which have already been derived,   and interpret them in a somewhat {\it ad hoc}, external way.

On the other hand, we should note that {\bf MILL1} grammars are in a way  more flexible and expressive then {\bf ETTC}, for  being able to encode as first order variables not only string endpoints, but also grammatical features, locality domains, tree depth levels etc (see \cite{MootPiazza}).

\subsubsection{Background}
We assume that the reader has some acquaintance with $\lambda$-calculus (see \cite{Barendregt}) as well as with Lambek calculus \cite{Lambek} and (multiplicative) linear logic, in particular, with basic ideas  of cut-elimination and the formalism of proof-nets \cite{Girard}, \cite{Girard2}.

\section{Tensor terms}
\subsection{Term expressions}
Throughout the paper we assume that we are given   an infinite  set $\mathit{Ind}$ of  {\it indices}. They will be used in all syntactic objects (terms, types, typing judgements) that we consider.

Now let  $T$ be  an alphabet of {\it terminal symbols} or, simply, {\it terminals}.

We will build {\it tensor term expressions} from terminal symbols, using    elements of $\mathit{Ind}$ as {\it upper} and {\it lower} indices. The set of upper, respectively lower, indices occurring in a
 term expression $t$ will be denoted as $Sup(t)$, respectively  $Sub(t)$. The set  of all indices of $t$ will be denoted  as $\mathit{Ind}(t)=Sup(t)\cup Sub(t)$.

 {\it Tensor terms} will be defined  as tensor term expressions quotiented by an appropriate equivalence.

Term expressions are defined by induction.
\begin{itemize}
  \item If $w\in T^*$ and $i,j\in Ind$,  then $t=[w]^i_j$ is a  term expression;
  \item if $w\in T^*$ then $t=[w]$ is a term expression;
  \item if $t_{(1)},t_{(2)}$ are term expressions with
  $$Sup(t_{(1)})\cap Sup(t_{(2)})=\emptyset,\quad Sub(t_{(1)})\cap Sub(t_{(2)})=\emptyset$$
  then $t=(t_{(1)}\cdot t_{(2)})$ is a term expression.
\end{itemize}
The multiplication symbol (dot)  will usually be omitted, i.e., we will write $(ts)$ for $(t\cdot s)$.

The definition  implies that an index can occur in a term expression at most twice: once as an upper one and once as a lower one.
We say that an index occurring in a term expression $t$ is  {\it free} in $t$, if it occurs in $t$ once. Otherwise we say that the index is  {\it bound}.

We denote the set of  free upper, respectively lower, indices of $t$   as $FSup(t)$, respectively $FSub(t)$.
 We denote the set of all free indices of $t$  as $F\mathit{Ind}(t)=FSup(t)\cup FSub(t)$.
We say that a  term  expression is {\it normal} if it has no bound indices.
We say that a term  expression is {\it closed} if it has no occurrence of a terminal symbol.

An {\it elementary}  term expression is  an expression of the form $[w]^j_i$ or $[w]$.
An elementary term expression $[w]^j_i$, where $j\not=i$ is called {\it regular}.
Elementary   term expressions $[w]^i_i$ and  $[w]$ are {\it singular}. Singular term expressions should be considered as pathological (in the context of this work), but we need them for consistency of definitions. We will discuss their meaning shortly.

\subsection{Terms}
We define {\it congruence} of term expressions as the smallest equivalence relation satisfying the conditions
$$((t\cdot s)\cdot k)\equiv(t\cdot (s\cdot k)),\quad
  (t\cdot s)\equiv(s\cdot t),\quad
  t\equiv s\Rightarrow (t\cdot k)\equiv(s\cdot k),$$
  $$[u]^j_i[v]^k_j\equiv[uv]^k_i,\quad [u]^i_i\equiv[u].$$
  $$[a_1a_2\ldots a_n]\equiv[a_na_1\ldots a_{n-1}],\quad\mbox{ for }a_1,\ldots, a_n\in T.$$
Congruence has a simple geometric meaning, which we will discuss shortly.

{\it Tensor terms}  (over a given alphabet of terminal symbols) are defined as equivalence classes for congruence of term expressions.

 Multiplication of terms is associative, therefore we will usually omit brackets in term expressions in the sequel.

The sets of free upper and lower indices of a term expression are easily seen to be invariant under congruence. Thus they are well-defined for terms as well. For a term $t$ we write  $FSup(t)$, $FSub(t)$ for the sets of free upper, respectively lower, indices of $t$ and $F\mathit{Ind}(t)=FSup(t)\cup Fsub(t)$.
We say that a  term is {\it closed} if it is a congruence class of a closed term expression.

A crucial role will be played by the following closed constants ({\it Kronecker deltas}), familiar from linear algebra:
$$\delta_{i}^j=[\epsilon]_i^j,\quad \delta_{i_1\ldots i_n}^{j_1\ldots j_n}=\delta_{i_1}^{j_1}\cdots \delta_{i_n}^{j_n},$$
where $\epsilon$ stands for the empty word.

We have the relations
$$\delta_i^j\cdot [w]_j^k \equiv [w]_i^k,\quad\delta_i^j\cdot [w]^i_k\equiv [w]^j_k,$$
and these imply the following property.
\bp\label{Kroneker}
Let $t$ be a tensor term expression, and
$$i_1,\ldots,i_m\in FSup(t), \quad j_1,\ldots,j_n\in FSub(t),$$
$$i_1',\ldots,i_m', j_1',\ldots,j_n'\not\in \mathit{Ind}(t).$$
Let the term expression $t'$ be obtained from $t$ by replacing the indices
$$i_1,\ldots,i_m, j_1,\ldots,j_n$$ with
$$i_1',\ldots,i_m', j_1',\ldots,j_n'$$ respectively.

Then
 $$\delta_{i_1,\ldots,i_m}^{i_1',\ldots,i_m'}\cdot\delta_{j_1',\ldots,j_n'}^{j_1,\ldots,j_n}\cdot t\equiv t'. \Box$$
\ep
Apparently, writing long sequences of indices as in the above proposition would be rather cumbersome.
In the sequel, we adopt the  convention that capital Latin letters stand for sequences of indices and small Latin letters, for individual indices.

\subsubsection{Geometric representation and normalization}
Tensor term expressions and terms have a simple intended geometric representation. ``Good'' (non-singular) term expressions can be represented as bipartite graphs whose vertices are labeled with indices and edges with words, and direction of edges is from lower indices to upper ones. Singular parts, such as $[w]^i_i$ and $[w]$, are represented as closed loops.
An elementary term expression $[w]_i^j$, $i\not=j$, corresponds to a single edge,
\begin{figure}
\begin{subfigure}{.3\textwidth}
\centering
\begin{tikzpicture}
%\node at (0,.25){$[w]^j_i$};
%\begin{scope}[shift={(0,-1.5)}]
\draw[thick,->](0,0)--(0.5,0);
\node[left] at(0,0){$i$};
\node[right] at(.5,0){$j$};
\node[above] at(.25,0){$w$};
\draw [fill] (0,0) circle [radius=0.05];
%\end{scope}
\end{tikzpicture}
\caption{$[w]^j_i$}
\end{subfigure}
\begin{subfigure}{.3\textwidth}
\centering
\begin{tikzpicture}
\draw[thick,->](0,0)--(0.5,0);
\node[left] at(0,0){$i$};
\node[right] at(.5,0){$j$};
\draw [fill] (0,0) circle [radius=0.05];
\draw[thick,->](0,.5)--(0.5,.5);
\node[left] at(0,.5){$k$};
\node[right] at(.5,0.5){$l$};
\node[above] at(.25,0){$u$};
\node[above] at(.25,.5){$v$};
\draw [fill] (0,.5) circle [radius=0.05];
\end{tikzpicture}
\caption{$[u]_i^j\cdot [v]_k^l$}
\end{subfigure}
\begin{subfigure}{.3\textwidth}
\centering
\begin{tikzpicture}
\draw[thick,-](0,0)--(0.5,0);
\node[left] at(0,0){$i$};
%\node[right] at(.5,0){$j.$};
\draw [fill] (0,0) circle [radius=0.05];
\draw[thick,->](0,.5)--(0.5,.5);
\node[left] at(0,.5){$k$};
\node[right] at(.5,0.5){$l$};
\node[above] at(.25,0){$u$};
\node[above] at(.25,.5){$v$};
\draw [fill] (0,.5) circle [radius=0.05];

%\node[left] at(0,1){$j$};
\node[right] at(.5,1){$r$};
\draw[thick,->](0,1)--(0.5,1);
\node[above] at(.25,1){$w$};
%\draw [fill] (0,1) circle [radius=0.05];
\draw[thick,-] (.5,0)to[out=0,in=90](.75,-.25)
to[out=-90,in=0](.5,-0.5)--(-.25,-0.5)
to[out=180,in=-90](-.5,-.25)--(-.5,.75)
to[out=90,in=0](0,1);

\node at(1.1,.5){$=$};

\begin{scope}[shift={(.5,0)}]
\begin{scope}[shift=({0,.25})]
\draw[thick,->](1.25,0)--(1.75,0);
\node[left] at(1.25,0){$i$};
\node[right] at(1.75,0){$r$};
\node[above] at(1.5,0){$uw$};
\draw [fill] (1.25,0) circle [radius=0.05];
\end{scope}

\begin{scope}[shift=({0,-.25})]
\draw[thick,->](1.25,1)--(1.75,1);
\node[left] at(1.25,1){$k$};
\node[right] at(1.75,1){$l$};

\node[above] at(1.5,1){$v$};
\draw [fill] (1.25,1) circle [radius=0.05];
\end{scope}
\end{scope}
\end{tikzpicture}
\caption{$[u]_i^j\cdot [v]_k^l\cdot[w]_j^r$}
\end{subfigure}
\caption{Geometric representation of terms}
\label{terms geometrically}
\end{figure}
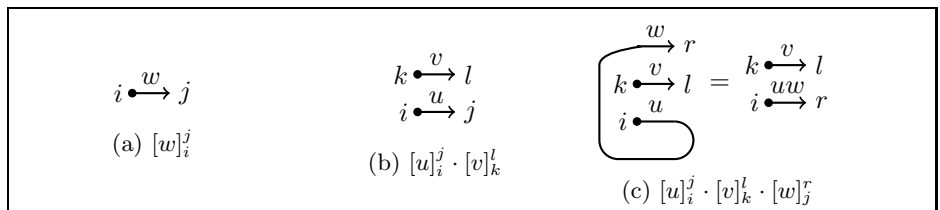
the product of two terms without common indices is represented as the disjoint union of the corresponding graphs, %the term $[u]_i^j\cdot [v]_k^l$ is represented as two edges.
and a term with bound indices corresponds to a graph obtained by gluing edges along matching vertices. This is illustrated in Figure \ref{terms geometrically}.
%Thus, the term $[u]_i^j\cdot [v]_k^l\cdot[w]_j^r$ is represented as the following.

As for singular terms, they  arise when edges are glued cyclically, for example, when the endpoints of the same edge are glued together as in $[w]_i^i$, which is the same as $[w]$.
Then a singular term  $[a_1\ldots a_n]$, where $a_1,\ldots,a_n\in T$, should be represented as a closed loop labeled with the cyclically ordered  sequence $\{a_1,\ldots,a_n\}$. The ordering is cyclic, because there is no way to say which letter is first. This is no longer a graph (at least, in the most common sense) because there are no vertices, but it has an obvious geometric meaning (it is a topological space, even a manifold).

Obviously, congruent term expressions have the same geometric representation. In fact, congruence of term expressions can be {\it defined} in terms of having the same geometric representation, and terms can be defined in terms of graphs and loops.
In general, a tensor term can be represented as  a finite set of word-labeled edges with index-labeled vertices and a finite multiset of closed loops labeled with cyclically ordered words. These geometric objects were introduced in \cite{Slavnov_cowordisms} under the name of {\it word cobordisms} or {\it cowordisms} for short, and {\it linear logic grammars} (LLG) were formulated directly in the geometric language. Tensor grammars of this work are a reformulation and an elaboration of constructions of \cite{Slavnov_cowordisms}.

The geometric representation makes  especially obvious   that
any  term  expression is congruent to a  normal one (which is unique up to associativity and commutativity of multiplication).

We say that
a term expression $t$ is regular if it is congruent to a  term expression of the form $[w_1]_{i_1}^{j_1}\cdots[w_1]_{i_n}^{j_n}$, where no index is repeated. Otherwise we say that $t$ is {\it singular}.
That is, a term expression is regular if its normal form is the product of elementary regular expressions.
We say that a tensor term is  {\it regular}, respectively {\it singular},  if it is the congruence class of a regular, respectively, singular term expression.
In a geometric language,  regular terms are those that correspond to loop-free graphs.

\section{Tensor type calculus}\label{tensor types section}
\subsection{Tensor types}
Our goal is to assign types to terms. Our types will be formulas of {\it multiplicative linear logic} ({\bf MLL}), decorated with indices intended to match free indices of terms.

We assume that we are given a set $P$ of  {\it positive atomic type symbols} or {\it positive literals}, and every  element $p\in P$ is assigned a {\it valency} $v(p)\in {\bf N}^2$.

The set $Lit$ of  {\it atomic type symbols} or {\it literals} is defined as $Lit=P\cup \overline{P}$,
where $\overline{P}=\{\overline{p}|~p\in P\}$.
Elements of $\overline{P}$,  {\it negative atomic type symbols} or {\it negative literals}, are assigned valencies by the rule: if $p\in P$ and $v(p)=(n,m)$, then  $v(\overline{p})=(m,n)$. The atomic type symbol $\overline p$ is the {\it dual} of the symbol $p$. For a negative literal $\overline p$ we denote $\overline{\overline{p}}=p$.

Types, similarly to terms, will have upper and lower indices.
Accordingly, we will denote the set of upper, respectively lower indices occurring in a type $A$ as $Sup(A)$, respectively $Sub(A)$ and the set of all indices occurring in $A$ as $\mathit{Ind}(A)=Sup(A)\cup Sub(A)$.

In a tensor type, each index occurs at most once, so there are no bound indices. However, later we will add more constructors, which allow binding indices in types as well. Therefore we will explicitly use the adjective ``free'' for type indices right from the start and use the tautological notation $FSup(A)=Sup(A)$, $FSub(A)=Sub(A)$, $Find(A)=\mathit{Ind}(A)$.

    Given a set $P$ of  positive atomic type symbols together with the valency function $v:P\to{\bf N}^2$, the set of {\it tensor types} over $P$ is defined by the following rules.
\begin{itemize}
  \item If $p\in Lit=P\cup\overline{P}$ with $v(p)=(m,n)$ and $i_1,\ldots, i_m,j_1,\ldots,j_n$ are pairwise distinct elements of $\mathit{Ind}$ then $A=p_{j_1\ldots j_n}^{i_1\ldots i_m}$ is a type;
  \item if $A$, $B$ are types  with $F\mathit{Ind}(A)\cap F\mathit{Ind}(B)=\emptyset$
  then $A\otimes B$, $A\wp B$ are types.
    \end{itemize}
    Occasionally we will use the title ``formula'' as  synonym of ``type''. (Typically,  we have a {\it subformula} $A$ of the {\it formula} $A\otimes B$, rather than  a subtype.)

A {\it tensor type symbol} (over $P$) is a tensor type (over $P$) with all  indices erased. We denote the set of tensor type symbols over $P$ as $Tp_\otimes(P)$.
{\it Type valency} extends from  atomic type symbols to all tensor type symbols in the obvious way:
\be\label{valency}
v(A\otimes B)=v(A\wp B)=v(A)+v(B),
\ee
where addition of elements of ${\bf N}^2$ is defined componentwise: $(m_1,n_1)+(m_2,n_2)=(m_1+m_2,n_1+n_2)$.

Tensor type symbols are  formulas of {\bf MLL} equipped with valencies. In order to obtain a type from a type symbol it is sufficient to specify the ordered sets of free upper and lower indices.
Accordingly, we will use the notation $A^I_J$ to denote a tensor type whose symbol is $A$, and whose ordered sets of free upper and lower indices are $I$ and $J$ respectively.

We now define {\it duality} (or {\it negation}) of tensor types.
For a finite sequence $I=(i_1,\ldots,i_n)$ of indices we denote $\overline{I}=({i_n},\ldots, {i_1})$.

 The {\it dual type } $\overline{A}$
 of a tensor type  $A$,
 is defined by induction:
\be\label{dual type}
\overline{p^I_J}=\overline{p}^{\overline{J}}_{\overline{I}}\mbox{ for }p\in Lit,\quad \overline{A\otimes B}=\overline{B}\wp\overline{A},\quad
\overline{A\wp B}=\overline{B}\otimes\overline{A}.
\ee
Note that we use the convention that duality flips tensor factors, which is typical for {\it noncommutative} variants of linear logic. This does not change the logic (the formulas $A\otimes B$ and $B\otimes A$ are provably equivalent in {\bf MLL}), but is more convenient for graphical interpretation of tensor type calculus.
We define the  {\it linear implication type}  by
\be\label{linear implication}
  A\multimap B=B\wp\overline{A}
  \ee
(again inverting the  conventional  order of subformulas).

\subsection{Typing judgements}
A {\it tensor  sequent} (over a set $P$)  is a finite set of tensor types (over $P$) whose elements have no common free indices.
Following the tradition, we write tensor sequents as sequences of types without any braces and denote them with capital Greek letters.
The sets of free upper and lower indices of $\Gamma$ will be denoted, respectively,
as
$$FSup(\Gamma)=\bigcup\limits_{A\in\Gamma}FSup(A),\quad
FSub(\Gamma)=\bigcup\limits_{A\in\Gamma}FSub(A).$$

A {\it syntactic tensor typing judgement} (over sets $P$ and $T$) is a pair $(t,\Gamma)$, denoted as $t\vdash \Gamma$, where
 $\Gamma$ is a tensor  sequent (over $P$),
  and $t$ is a term expression (over $T$), such that $FSup(t)=FSub(\Gamma)$, $FSub(t)=FSup(\Gamma)$.

We define {\it congruence} of syntactic typing judgements as the  smallest  equivalence relation on syntactic tensor typing judgements satisfying the conditions:
\begin{itemize}
  \item if $t\vdash\Gamma$, $t'\vdash\Gamma'$ are syntactic typing judgements, $i\in F\mathit{Ind}(t)$, $i'\not\in \mathit{Ind}(t)$, and $t'$, respectively $\Gamma'$, are obtained by replacing the index $i$ in $t$, respectively $\Gamma$, with $i'$, then $t\vdash\Gamma\equiv t'\vdash\Gamma'$ ({\it $\alpha$-conversion} );
  \item if $t\vdash\Gamma$, $t'\vdash\Gamma$ are syntactic typing judgements, and $t\equiv t'$, then $t\vdash\Gamma\equiv t'\vdash\Gamma$.
\end{itemize}
   A {\it  tensor typing judgement} as an equivalence class for congruence of syntactic tensor typing judgements.
      Thus, typing judgements are identified by $\alpha$-conversion and by congruence of terms.

We say that a typing judgement $t\vdash\Gamma$ is {\it regular} if the  term $t$ is regular (i.e. $t$ is represented geometrically as a loop-free graph).
\bp[On change of coordinates]\label{alpha conversion}
Typing judgements
$$t\vdash \Gamma, A^{I}_{J},\quad\delta^{{I}J'}_{I' {J}}\cdot t\vdash \Gamma, A^{I'}_{J'}$$
(where $I,J,I',J'$ are pairwise disjoint)
 are congruent.
\ep
{\bf Proof}
Let $I=i_1\ldots i_n$, $I'=i_1'\ldots i_n'$,  $J=j_1\ldots j_m$,  $J'=j_1'\ldots j_m'$, and let us write
 $t'=\delta^{{I}J'}_{I' {J}}\cdot t$ for the term in the second of the typing judgements above.

Multiplying a tensor term by a Kronecker delta amounts to renaming indices (Proposition \ref{Kroneker}). The term $t$ contains no primed indices, so $t'$ is obtained as $t$ where each primed index $i_k'$ or $j_k'$ is replaced with $i_k$ or $j_k$ respectively. The two typing judgements are congruent  by $\alpha$-conversion. $\Box$

\subsection{Typing rules}\label{tensor typing rules}
Typing judgements are derived by the following rules of {\it Tensor type calculus} ({\bf TTC}).
$$\delta^{I'\overline{J}}_{\overline{I}J'}\vdash \overline{p}^{I}_{J},p_{I'}^{J'},~p\in P~({\rm{Id}}),\quad
\frac{t\vdash \Gamma,A\quad s\vdash
\overline{A},\Theta}{ts\vdash\Gamma,\Theta} ~({\rm{Cut}}),$$
$$\frac{t\vdash
\Gamma,A,B}{t\vdash
\Gamma,A\wp B}
~
(\wp)\quad\frac{t\vdash\Gamma, A \quad s\vdash
B,\Theta}{ts\vdash\Gamma,A\otimes B,\Theta}~{ }
(\otimes).$$

  It is implicit in the rules above that all typing judgements are well defined, i.e. there are no index collisions.  For example, syntactic representatives for  the two premises of the $(\otimes)$ rule must be chosen to  have no common indices. Similarly, there should be no common indices in the contexts $\Gamma$ and $\Theta$ for the Cut rule. Note that the term expression $ts$ in the conclusion of the Cut rule will not be normal: all indices occurring in the cut-types $A$, $\overline{A}$ will be bound in $ts$.

The above presentation of the Cut rule is not always very convenient for theoretical arguments. An alternative version is obtained if we choose syntactic representatives for the premises to have no common indices.
\bp\label{alternative cut}
For any tensor type symbol $A$ the Cut rule can be written as
$$\cfrac
{ t\vdash \Gamma,A^I_J\quad s\vdash \overline{A}^{J'}_{I'},\Theta}
{\delta^{\overline{I}J'}_{I'\overline {J}}\cdot ts\vdash\Gamma,\Theta}({\rm{Cut}}),
$$
where the index sentences $I,J,I',J'$ are pairwise disjoint.
\ep
{\bf Proof} Let $s'=\delta^{\overline{I}J'}_{I'\overline {J}}\cdot s$. By Proposition  \ref{alpha conversion} on change of coordinates, the second premise in the rule above is congruent to
$s'\vdash \overline{A}^{\overline{J}}_{\overline{I}}$. On the other hand, the term in the conclusion is $s'\cdot t$. This restores the original form of the Cut rule. $\Box$

  \subsection{Graphical representation}
    Of course, the rules of {\bf TTC} are those of {\it multiplicative linear logic} ({\bf MLL}), decorated with indices and terms. One might observe that the syntax of tensor typing judgements is, basically, an in-line notation for {\bf MLL} proof-nets.

  While tensor terms have a natural interpretation as edge-labeled graphs, tensor typing judgements  suggest particular pictorial representation for such graphs. A type to the right of the turnstile can be seen as inducing an ordering on the set of indices, for example, from left to right, from   top to bottom, and this can be seen as an ordering of vertices in a picture. Thus, a typing judgement  $t\vdash F$, where $F$ is a formula,  can be depicted as the graph representing the term $t$ with vertices aligned, say, horizontally according to the ordering induced by $F$. Identifying  a sequent with the par of its formulas we get a prescription for depicting more general typing judgements of the form $t\vdash\Gamma$, where $\Gamma$ contains more than one formula. (This is not completely accurate with our current definitions though, because we defined  sequent as an unordered set of formulas.)

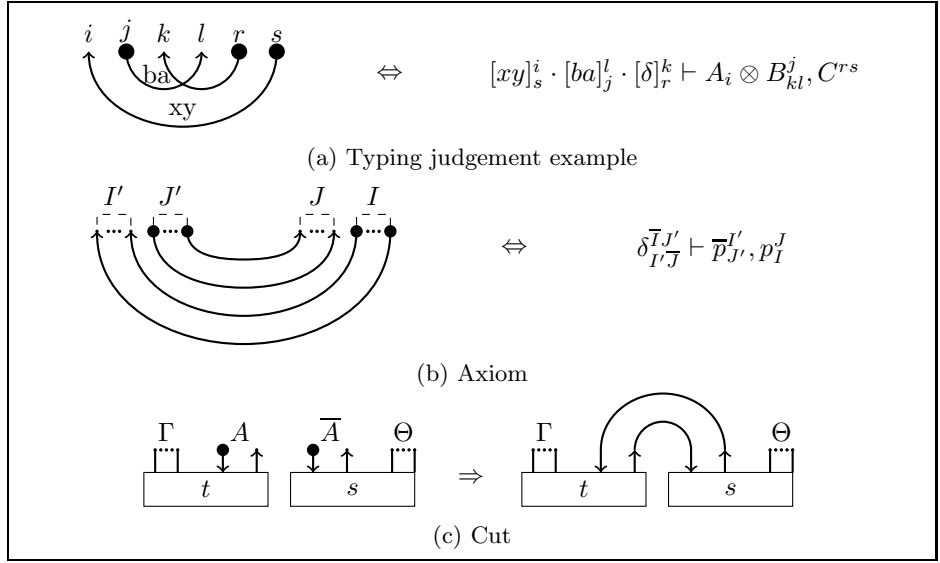
\begin{figure}
\begin{subfigure}{1\textwidth}
\centering
         \begin{tikzpicture}[xscale=2,yscale=2]
    \begin{scope}[shift={(.5,0)}]
        %\draw  (0,0) circle [radius=0.05];
        \draw[thick,->](1.25,0) to  [out=-90,in=0] (.625,-.5) to  [out=180,in=-90](0,0);
        \node at(.625,-.4){xy};
        %\node at(.75,.35){a};
        \node at(.45,-.15){ba};
        \draw[thick,->](1,0) to  [out=-90,in=0] (.75,-.25) to  [out=180,in=-90](0.5,0);
        \draw[thick,<-](0.75,0) to  [out=-90,in=0] (.5,-.25) to  [out=180,in=-90](.25,0);

        \draw [fill] (.25,0) circle [radius=0.05];
        \draw [fill] (1,0) circle [radius=0.05];
        \draw [fill] (1.25,0) circle [radius=0.05];

        \node [above]at(0,0){$i$};
        \node [above]at(0.25,0){$j$};
        \node [above]at(0.5,0){$k$};
        \node [above]at(0.75,0){$l$};
        \node [above]at(1,0){$r$};
        \node [above]at(1.25,0){$s$};
        \end{scope}
        \node at (2.5,-.15){$\Leftrightarrow$};
         \node[right] at (3.1,-.15) {$[xy]_s^i\cdot [ba]_j^l\cdot[\delta]_r^k\vdash A_i\otimes B^j_{kl}, C^{rs}$};

        \end{tikzpicture}
        \caption{Typing judgement example}
        \label{Typing_judgement_geometrically}
        \end{subfigure}\\
        \begin{subfigure}{1\textwidth}
        \centering
         \begin{tikzpicture}[xscale=1.5,yscale=1.5]

         \begin{scope}[xscale=-1]
        \begin{scope}[shift={(2.6,0)}]

        \begin{scope}[shift={(-2.5,0)}]
          \draw[thick,<-](0.8,0) to  [out=-90,in=180] (1.3,-.25) to  [out=0,in=-90](1.8,0);
          \draw[thick,<-](.5,0) to  [out=-90,in=180] (1.3,-.5) to  [out=0,in=-90](2.1,0);
          \draw[thick,->](0.3,0) to  [out=-90,in=180] (1.3,-.75) to  [out=0,in=-90](2.3,0);
          \draw[thick,->](0,0) to  [out=-90,in=180] (1.3,-1) to  [out=0,in=-90](2.6,0);

                \draw [fill] (0,0) circle [radius=0.05];
                \draw [fill] (0.15,0) circle [radius=0.01];
                \draw [fill] (0.2,0) circle [radius=0.01];
                \draw [fill] (0.1,0) circle [radius=0.01];
                        \draw [fill] (0.3,0) circle [radius=0.05];
                \draw[dashed,-](0,.05)--(0,.15)--(0.3,.15)--(0.3,.05);
                \node at (0.15,.15) [above] {${I}$};

                \begin{scope}[shift={(.5,0)}]
                \draw [fill] (0,0) circle [radius=0.01];
                \draw [fill] (0.15,0) circle [radius=0.01];
                \draw [fill] (0.2,0) circle [radius=0.01];
                \draw [fill] (0.1,0) circle [radius=0.01];
                \draw [fill] (0.3,0) circle [radius=0.01];
                \draw[dashed,-](0,.05)--(0,.15)--(0.3,.15)--(0.3,.05);
                \node at (0.15,.15) [above] {${J}$};
                \end{scope}

            \begin{scope}[shift={(1.8,0)}]
            \draw [fill] (0,0) circle [radius=0.05];
                \draw [fill] (0.15,0) circle [radius=0.01];
                \draw [fill] (0.2,0) circle [radius=0.01];
                \draw [fill] (0.1,0) circle [radius=0.01];
                        \draw [fill] (0.3,0) circle [radius=0.05];
                \draw[dashed,-](0,.05)--(0,.15)--(0.3,.15)--(0.3,.05);
                \node at (0.15,.15) [above] {$J'$};

                \begin{scope}[shift={(.5,0)}]
                \draw [fill] (0,0) circle [radius=0.01];
                \draw [fill] (0.15,0) circle [radius=0.01];
                \draw [fill] (0.2,0) circle [radius=0.01];
                \draw [fill] (0.1,0) circle [radius=0.01];
                \draw [fill] (0.3,0) circle [radius=0.01];
                \draw[dashed,-](0,.05)--(0,.15)--(0.3,.15)--(0.3,.05);
                \node at (0.15,.15) [above] {$I'$};
                \end{scope}
            \end{scope}
          \end{scope}

          \end{scope}
          \end{scope}

               \node at (1.,-.15){$\Leftrightarrow$};
         \node[right] at (2,-.15) {$\delta_{I'\overline{J}}^{\overline{I}J'}\vdash \overline{p}^{I'}_{J'},p_{I}^{J}\quad\quad$};

        \end{tikzpicture}
        \caption{Axiom}
        \label{Axiom_geometrically}
        \end{subfigure}
\begin{subfigure}{1\textwidth}
        \centering
         \begin{tikzpicture}[xscale=1.5,yscale=1.5]
         \begin{scope}[shift={(-3.2,0)}]
         \begin{scope}[shift={(.5,0)}]
             \draw[draw=black](.9,-.2)rectangle(-.2,-.5);\node at(.35,-.35){$t$};

                    \draw [fill] (-.1,0) circle [radius=0.01];
                    \draw [fill] (-0.05,0) circle [radius=0.01];
                    \draw [fill] (0,0) circle [radius=0.01];
                    \draw [fill] (0.05,0) circle [radius=0.01];
                            \draw [fill] (0.1,0) circle [radius=0.01];
                    \draw[thick,-](-.1,0)--(-.1,-.2);
                    \draw[thick,-](.1,0)--(.1,-.2);
                    \node at (0,0) [above] {$\Gamma$};

                    \begin{scope}[shift={(.5,0)}]
                    \draw [fill] (0,0) circle [radius=0.05];
                    %\draw [fill] (0.15,0) circle [radius=0.01];
    %                \draw [fill] (0.2,0) circle [radius=0.01];
    %                \draw [fill] (0.1,0) circle [radius=0.01];
                    \draw [fill] (0.3,0) circle [radius=0.01];
                    \draw[thick,->](0,0)--(0,-.2);
                    \draw[thick,<-](0.3,0)--(0.3,-.2);
                    \node at (0.15,0) [above] {$A$};
    %                \node at (0.3,0) [right] {$I$};
                    \end{scope}

                \begin{scope}[shift={(1.3,0)}]
                \draw[draw=black](.9,-.2)rectangle(-.2,-.5);\node at(.35,-.35){$s$};
                 \draw [fill] (0,0) circle [radius=0.05];
                    %\draw [fill] (0.15,0) circle [radius=0.01];
    %                \draw [fill] (0.2,0) circle [radius=0.01];
    %                \draw [fill] (0.1,0) circle [radius=0.01];
                    \draw [fill] (0.3,0) circle [radius=0.01];
                    \draw[thick,->](0,0)--(0,-.2);
                    \draw[thick,<-](0.3,0)--(0.3,-.2);
                    \node at (0.15,0) [above] {$\overline{A}$};
                    %\node at (0,0) [left] {$I'$};
    %                \node at (0.3,0) [right] {$J'$};
                \end{scope}
      \end{scope}

                %\draw [fill] (0,0) circle [radius=0.05];
    %                \draw [fill] (0.15,0) circle [radius=0.01];
    %                \draw [fill] (0.2,0) circle [radius=0.01];
    %                \draw [fill] (0.1,0) circle [radius=0.01];
    %                        \draw [fill] (0.3,0) circle [radius=0.05];
    %                \draw[dashed,-](0,-.05)--(0,-.15)--(0.3,-.15)--(0.3,-.05);
    %                \node at (0.15,-.15) [below] {$J$};

                    \begin{scope}[shift={(2.6,0)}]
                    \draw [fill] (-.1,0) circle [radius=0.01];
                    \draw [fill] (-0.05,0) circle [radius=0.01];
                    \draw [fill] (0,0) circle [radius=0.01];
                    \draw [fill] (0.05,0) circle [radius=0.01];
                            \draw [fill] (0.1,0) circle [radius=0.01];
                    \draw[thick,-](-.1,0)--(-.1,-.2);
                    \draw[thick,-](.1,0)--(.1,-.2);
                    \node at (0,0) [above] {$\Theta$};

                    \end{scope}
        \end{scope}

        \node at (0,-.25){$\Rightarrow$};

        \begin{scope}[shift={(.65,0)}]
        \draw[draw=black](.9,-.2)rectangle(-.2,-.5);\node at(.35,-.35){$t$};

          \draw[thick,-](0.8,0) to  [out=90,in=180] (1.05,.25) to  [out=0,in=90](1.3,0);
          \draw[thick,-](.5,0) to  [out=90,in=180] (1.05,.5) to  [out=0,in=90](1.6,0);
          %\draw[thick,dotted,-](0.3,0) to  [out=-90,in=180] (1.3,-.75) to  [out=0,in=-90](2.3,0);
          %\draw[thick,dotted,-](0,0) to  [out=-90,in=180] (1.3,-1) to  [out=0,in=-90](2.6,0);

                \draw [fill] (-.1,0) circle [radius=0.01];
                \draw [fill] (-0.05,0) circle [radius=0.01];
                \draw [fill] (0,0) circle [radius=0.01];
                \draw [fill] (0.05,0) circle [radius=0.01];
                        \draw [fill] (0.1,0) circle [radius=0.01];
                \draw[thick,-](-.1,0)--(-.1,-.2);
                \draw[thick,-](.1,0)--(.1,-.2);
                \node at (0,0) [above] {$\Gamma$};

                \begin{scope}[shift={(.5,0)}]
               % \draw [fill] (0,0) circle [radius=0.05];
                %\draw [fill] (0.15,0) circle [radius=0.01];
%                \draw [fill] (0.2,0) circle [radius=0.01];
%                \draw [fill] (0.1,0) circle [radius=0.01];
                %\draw [fill] (0.3,0) circle [radius=0.01];
                \draw[thick,->](0,0)--(0,-.2);
                \draw[thick,<-](0.3,0)--(0.3,-.2);
             %   \node at (0,0) [left] {$J$};
%                \node at (0.3,0) [right] {$I$};
                \end{scope}

            \begin{scope}[shift={(1.3,0)}]
            \draw[draw=black](.9,-.2)rectangle(-.2,-.5);\node at(.35,-.35){$s$};
            % \draw [fill] (0,0) circle [radius=0.05];
                %\draw [fill] (0.15,0) circle [radius=0.01];
%                \draw [fill] (0.2,0) circle [radius=0.01];
%                \draw [fill] (0.1,0) circle [radius=0.01];
               % \draw [fill] (0.3,0) circle [radius=0.01];
                \draw[thick,->](0,0)--(0,-.2);
                \draw[thick,<-](0.3,0)--(0.3,-.2);
                %\node at (0,0) [left] {$I'$};
%                \node at (0.3,0) [right] {$J'$};
            \end{scope}
            %\draw [fill] (0,0) circle [radius=0.05];
%                \draw [fill] (0.15,0) circle [radius=0.01];
%                \draw [fill] (0.2,0) circle [radius=0.01];
%                \draw [fill] (0.1,0) circle [radius=0.01];
%                        \draw [fill] (0.3,0) circle [radius=0.05];
%                \draw[dashed,-](0,-.05)--(0,-.15)--(0.3,-.15)--(0.3,-.05);
%                \node at (0.15,-.15) [below] {$J$};

                \begin{scope}[shift={(2.1,0)}]
                \draw [fill] (-.1,0) circle [radius=0.01];
                \draw [fill] (-0.05,0) circle [radius=0.01];
                \draw [fill] (0,0) circle [radius=0.01];
                \draw [fill] (0.05,0) circle [radius=0.01];
                        \draw [fill] (0.1,0) circle [radius=0.01];
                \draw[thick,-](-.1,0)--(-.1,-.2);
                \draw[thick,-](.1,0)--(.1,-.2);
                \node at (0,0) [above] {$\Theta$};

                \end{scope}
        \end{scope}
        \end{tikzpicture}
        \caption{Cut}
        \label{Cut_geometrically}
        \end{subfigure}
        \caption{Typing judgements geometrically}
   \end{figure}

Typically, the typing judgement
$$[xy]_s^i\cdot [ba]_j^l\cdot[\delta]_r^k\vdash A_i\otimes B_{kl}^j, C^{rs}$$
should be depicted as in Figure \ref{Typing_judgement_geometrically}. We emphasize that there is no reason for pictures to be planar.
Note that in the geometric representation the indices become redundant, because of $\alpha$-conversion of typing judgements.

   With these conventions the axiom $({\rm{Id}})$ becomes nothing but a generalized axiom link of a proof-net as in Figure \ref{Axiom_geometrically}, while the Cut rule corresponds to attaching  cut-links connecting matching indices/vertices (this is the geometric meaning of Proposition \ref{alternative cut}). A simple case, when the cut-type  is  of valency $(1,1)$, so that there is only one upper and one lower index, is shown schematically in Figure \ref{Cut_geometrically}. In general the ``cut-link'' will consist of as many pieces as there are indices in the cut-type. In Figure \ref{Cut_geometrically} we use the format of depicting typing judgements to which we will stick from now on:  we do not write indices in pictures, but indicate occurrences of tensor type symbols.

   As for the $(\otimes)$ rule, it corresponds to putting two pictures side by side, while the $(\wp)$ rule,  on the level of pictures, does strictly nothing.

   Such a pictorial syntax, which we use only for informal illustrations in this work, can be developed systematically and rigorously, allowing dispensing with terms and indices altogether (which is rather convenient in simple cases, but leads to very space-consuming two-dimensional computations in general). This approach was taken in \cite{Slavnov_cowordisms}.

Geometric representation also suggests certain terminology with a proof-net flavour.

 We cannot speak accurately about a tensor {\it type} (which has particular indices) occurring in a typing judgement, at least, when the syntactic representation is not chosen. But occurrences of a tensor type {\it symbol} in a typing judgement are well defined. We say that an occurrence of a tensor type symbol $X$ in a tensor typing judgement $\sigma$ is a {\it link} $X$ in $\sigma$.
We say that a link $Z=X\otimes Y$ in a typing judgement $\sigma$ is {\it splitting}, if $\sigma$ can be written as
$$t\cdot s\vdash \Gamma, X^I_J\otimes Y^K_L,\Theta,$$
 where
\be\label{splitting}
t\vdash\Gamma, X^I_J,\quad s\vdash Y^K_L,\Theta
\ee
are well-defined typing judgements (i.e. indices on the left of the turnstile match those on the right). We say that the link $Z$ {\it splits $\sigma$} into the two typing judgements in (\ref{splitting}).

\subsubsection{Properties of {\bf TTC}}
Geometric representation makes it evident that the tensor type calculus is cut-free. Cut-elimination algorithm repeats cut-elimination for proof-nets. Observe that the case of a cut with an axiom is described by Proposition \ref{alpha conversion} on change of coordinates.
\bl\label{Cut-elimination simple}
Any  typing judgement derivable in {\bf TTC} can be derived  without the Cut rule. $\Box$
\el
\bc Any  typing judgement derivable in {\bf TTC} is regular.
\ec
{\bf Proof} by induction on a cut-free derivation. $\Box$
        \bp\label{par is hypocricy} The following rule is admissible in {\bf TTC}:
          $$
\cfrac{t\vdash\Gamma,A\wp B}{t\vdash\Gamma,A,B}~{(\wp^{-1})}.
$$
\ep
{\bf Proof} Let $A=X^I_J$, $B=Y^K_L$, where $X,Y$ are tensor type symbols.

The rule is obtained by cutting the premise with the  derivable typing judgement
$$\delta^{IJ'}_{I'J}\cdot\delta_{LK'}^{L'K}\vdash \overline{Y}_{\overline{K}}^{\overline{L}}\otimes\overline{X}^{\overline{J}}_{\overline{I}}, X_{J'}^{I'},Y^{K'}_{L'}$$
 and using Proposition \ref{alpha conversion} on change of coordinates. $\Box$

\subsection{Tensor signatures and deduction theorem}
We will be interested in tensor derivations with additional, non-logical axioms.

Let $\Xi$ be a finite multiset of typing judgements (over some sets $P$ and $T$).

We say that a typing judgement $\sigma$ (over $P$ and $T$) is {\it derivable} from $\Xi$ if there exists a  rooted tree,  whose leaves are labeled with {\bf TTC} axioms and elements of $\Xi$ with each element of $\Xi$ used as a label exactly once (counting multiplicity), whose root is labeled with $\sigma$, and each internal vertex is labeled with a tensor typing judgement obtained from the labels of its children by a {\bf TTC} rule. The tree is called a {\it derivation} of $\sigma$ from $\Xi$.

A {\it tensor signature} $\Sigma$ is a tuple $\Sigma=(P,T,\Xi)$), where $P$ is a set of positive atomic type symbols, $T$ is an alphabet of terminal symbols and $\Xi$ is a  set of typing judgements over $P$ and $T$, called {\it axioms of $\Sigma$}.
A typing judgement $\sigma$ is {\it derivable in $\Sigma$} if it is derivable in {\bf TTC} from some finite multiset $\Xi_0$ consisting of elements of $\Xi$.

It would be more traditional to require  that every axiom of a signature has exactly one type to the right of the turnstile. The more restrictive version is essentially equivalent to the more liberal definition that we give, because of the $(\wp)$ and $(\wp)^{-1}$ rules. The liberal version is more convenient in some situations, because derivations are shorter. However when axioms are in the more restrictive format we can use the following.
\bl[``Deduction theorem'' for TTC]\label{deduction theorem}
Let $\Xi$ be  a finite multiset of typing judgements
\be\label{Xi}
\Xi=\{t_{(1)}\vdash(F_{(1)})^{I_1}_{J_1},~\ldots,~t_{(n)}\vdash(F_{(n)})^{I_n}_{J_n}\}.
\ee
where each of $F_{(1)},\ldots,F_{(k)}$ is a single tensor type.
Choose syntactic representatives for elements of $\Xi$ so that all index sequences $I_1,\ldots,I_n$ $J_1,\ldots, J_n$ above are pairwise disjoint.

A  typing judgement $\sigma$
 is derivable
 from  $\Xi$ iff it  has a  syntactic representation of the form $t\vdash\Gamma$, where
\be\label{2'}
t=t_{(0)}\cdot t_{(1)}\cdots t_{(n)},
\ee
$\Gamma$ has no common indices with $I_1,\ldots,I_n$, $J_1,\ldots,J_n$, and
the typing judgement
\be\label{3'}
t_{(0)}\vdash
\overline{(F_{(1)})}_{\overline{I_1}}^{\overline{J_1}}, \ldots, \overline{(F_{(n)})}_{\overline{I_n}}^{\overline{J_n}},\Gamma
\ee
is derivable in {\bf TTC}.
\el
{\bf Proof} If a  representation $t \vdash \Gamma$  as in the formulation exists, then $\sigma$ is derivable from (\ref{Xi}) by iterated cuts with (\ref{3'}).

In the other direction use induction on derivation of $\sigma$.

Most important is the base of induction, when $\sigma$ is an element $t_{(1)}\vdash(F_{(1)})^{I_1}_{J_1}$ of $\Xi$.
For this case choose fresh index sequences $I_1'$, $J_1'$ to replace $I_1$, $J_1$ and put
$t_{(0)}=\delta^{I_1'J_1}_{I_1J_1'}$. By Proposition \ref{alpha conversion} on change of coordinates, the syntactic typing judgement $t_{(0)}t_{(1)}\vdash (F_{(1)})^{I_1'}_{J_1'}$ is just another representation of $\sigma$.

Induction steps are very easy if we use  the Cut rule in the form of Proposition \ref{alternative cut}. $\Box$

\subsection{Grammars}\label{tensor grammars section}
A {\it tensor grammar} $G=(\Sigma, S)$ is a tensor signature $\Sigma=(P,T,\Xi)$, where $T$ and $\Xi$ are finite, together with a positive atomic type symbol $S\in P$ of valency (1,1).
We say that $G$ {\it generates} a word $w\in T^*$  if the typing judgement $$[w]_i^j\vdash S_j^i$$ is derivable in $\Sigma$.
The {\it language} generated by $G$ ({\it language of} $G$) is the set of all words generated by $G$.

Lemma \ref{deduction theorem} (``Deduction theorem'') shows that we can equivalently define the language of $G$ in terms of derivations in the purely logical, cut-free system of {\bf TTC}. This is more convenient for proof-search and theoretical analysis. We stick to our version, because derivations from axioms are much simpler and shorter.

\section{Examples and inadequacy}\label{Examples}
The formalism of tensor grammars can be seen as a surface representation of {\it abstract categorial grammars} (ACG), which will be discussed in the next section.
Derivations of ACG translate to derivations of tensor grammars and this translation is isomorphic on the  level of string languages.

In this section we discuss toy examples, adapted, in fact, from ACG, and then analyse the notorious ``inadequacy'' \cite{Moot_inadequacy} of commutative grammars (such as ACG), which becomes very transparent in tensor representation.

\subsection{``John loves Mary''}\label{simple example}
Consider the  alphabet of atomic type symbols
$P=\{NP,S\}$ with $v(NP)=v(S)=(1,1)$.

\begin{figure}%[p]%[h!]
    \centering
    \begin{subfigure}{1\textwidth}
 \centering
\begin{tikzpicture}
    \begin{scope}[xscale=1]
         \draw[thick,->](0,0) to  [out=-90,in=180] (0.4,-0.8) to  [out=0,in=-90] (.8,0);
         \draw [fill] (0,0) circle [radius=0.05];
          \node[above ] at (0.4,0) {$NP$};
           \node[above] at (0.4,-.6) {$\mbox{John}$};
        \node[below] at (0.4,-0.8){``John''};

        \begin{scope}[shift ={(1.5,0)}]
         \draw[thick,->](0,0) to  [out=-90,in=180] (0.4,-0.8) to  [out=0,in=-90] (0.8,0);
         \draw [fill] (0,0) circle [radius=0.05];
          \node[above ] at (0.4,0) {$NP$};
           \node[above] at (0.4,-.6) {$\mbox{Mary}$};
        \node[below]  at (0.4,-0.8){``Mary''};
        \end{scope}

        \begin{scope}[shift ={(3,0)}]

        \draw[thick,->](0,.7)--(0,0) to  [out=-90,in=180] (1.75,-1.2) to  [out=0,in=-90] (3.5,0)
        --(2,.4)--(2,.7);
         \draw [fill] (0,.7) circle [radius=0.05];
        \draw[thick,<-](.5,.7)--(0.5,0) to  [out=-90,in=180] (0.75,-0.6) to  [out=0,in=-90] (1.5,0)
        --(3,.4)--(3,.7);
        \draw [fill] (3,.7) circle [radius=0.05];
        \draw[thick,<-]
        (3.5,.7)--(3.5,.4)--
        (2,0) to  [out=-90,in=180] (2.5,-0.6) to  [out=0,in=-90] (3,0)
        --(1.5,.4)--(1.5,.7);
        \draw [fill] (1.5,.7) circle [radius=0.05];
        \node[above] at (2.5,-.4) {$\mbox{loves}$};
          \node[above] at (0.3,.7) {$S$};
        \node[above] at (1.8,.7) {$\overline{NP}$};
        \node[above ] at (3.25,.7) {$\overline{NP}$};
        \node  at(1.75,-1.5){``loves''};

        \end{scope}
        \end{scope}
        \end{tikzpicture}
                \caption{Axioms}
\label{John_loves_Mary}
\end{subfigure}
%%%%%%%%%%%%%%%%%%%%%%%%%%%%%%%%%%%%%%%%%%%%%%%%
\begin{subfigure}{1\textwidth}
\centering
\scalebox{.75}[1.]
{
 \begin{tikzpicture}[xscale=1]
\begin{scope}[shift={(4,.7)}]
        \draw[thick,->](0,0) to  [out=-90,in=180] (0.4,-0.8) to  [out=0,in=-90] (0.8,0);
         \draw [fill] (0,0) circle [radius=0.05];
          \node[above ] at (0.4,0) {$NP$};
           \node[above] at (0.4,-.6) {$\mbox{Mary}$};

\end{scope}

\begin{scope}%[shift={(-.25,0)}]

        \draw[thick,->](0,.7)--(0,0) to  [out=-90,in=180] (1.75,-1.2) to  [out=0,in=-90] (3.5,0)
        --(2,.4)--(2,.7);
         \draw [fill] (0,.7) circle [radius=0.05];
        \draw[thick,<-](.5,.7)--(0.5,0) to  [out=-90,in=180] (0.75,-0.6) to  [out=0,in=-90] (1.5,0)
        --(3,.4)--(3,.7);
        \draw [fill] (3,.7) circle [radius=0.05];
        \draw[thick,<-]
        (3.5,.7)--(3.5,.4)--
        (2,0) to  [out=-90,in=180] (2.5,-0.6) to  [out=0,in=-90] (3,0)
        --(1.5,.4)--(1.5,.7);
        \draw [fill] (1.5,.7) circle [radius=0.05];
        \node[above] at (2.5,-.4) {$\mbox{loves}$};
          \node[above] at (0.3,.7) {$S$};
        \node[above] at (1.8,.7) {$\overline{NP}$};
        \node[above ] at (3.25,.7) {$\overline{NP}$};

\end{scope}

{
\node at(5.3,-.5){$\Longrightarrow$};
}

\begin{scope}
\begin{scope}

\begin{scope}[shift={(6.,0)}]

\begin{scope}[shift={(3.8,0.7)}]
        \draw[thick,-](0,0) to  [out=-90,in=180] (0.4,-0.8) to  [out=0,in=-90] (0.8,0);
         %\draw [fill] (0,0) circle [radius=0.05];
          %\node[above ] at (0.4,0) {$NP$};
           \node[above] at (0.4,-.6) {$\mbox{Mary}$};

\end{scope}

\draw[thick](3,.7) to  [out=90,in=180] (3.8,1.3) to  [out=0,in=90] (4.6,.7);
\draw[thick](3.5,.7) to  [out=90,in=180] (3.65,1) to  [out=0,in=90] (3.8,.7);

\begin{scope}%[shift={(-.25,0)}]

        \draw[thick,->](0,.7)--(0,0) to  [out=-90,in=180] (1.75,-1.2) to  [out=0,in=-90] (3.5,0)
        --(2,.4)--(2,.7);
         \draw [fill] (0,.7) circle [radius=0.05];
        \draw[thick,<-](.5,.7)--(0.5,0) to  [out=-90,in=180] (0.75,-0.6) to  [out=0,in=-90] (1.5,0)
        --(3,.4)--(3,.7);
        %\draw [fill] (3,.7) circle [radius=0.05];
        \draw[thick,-]
        (3.5,.7)--(3.5,.4)--
        (2,0) to  [out=-90,in=180] (2.5,-0.6) to  [out=0,in=-90] (3,0)
        --(1.5,.4)--(1.5,.7);
        \draw [fill] (1.5,.7) circle [radius=0.05];
        \node[above] at (2.5,-.4) {$\mbox{loves}$};
          \node[above] at (0.3,.7) {$S$};
        \node[above] at (1.8,.7) {$\overline{NP}$};
        %\node[above ] at (3.25,.7) {$\overline{NP}$};

\end{scope}

\end{scope}
\node [font=\large] at(10.8,-.5) {
$=$
};

\begin{scope}[shift={(10.5,0)}]

    \begin{scope}
     \draw[thick,->](1.1,0) to  [out=-90,in=180] (2.5,-1.25) to  [out=0,in=-90] (3.9,0);
    % to [out=90,in=180] (4,0.5) to [out=0,in=90]

 %    \draw [fill] (2.75,-1.5) circle [radius=0.05];
%    \draw[thick,<-](2.25,-1.5) to  [out=90,in=0] (1.75,-0.75) to  [out=180,in=-90] (1.5,0);
    \draw [fill] (1.1,0) circle [radius=0.05];
     %to [out=90,in=180] (4,2.5) to [out=0,in=90]

    %;

    %\draw[thick,<-](2,0) to  [out=-90,in=180] (2.5,-0.75) to  [out=0,in=-90] (3,0);

    \node[above] at (2.5,-.75) {$``\mbox{loves Mary}''$};

      \node[above] at (1.25,0) {$S$};
      \node[above] at (3.75,0) {$\overline{NP}$};
    %\node[above] at (4.8,0) {$S^\bot$};
    %\node[above] at (6.3,0) {$NP$};

    % \draw[thick,->](4.5,0) to  [out=-90,in=180] (7.25,-2.5) to  [out=0,in=-90] (10,0);
    \draw[thick,<-](1.4,0) to  [out=-90,in=180] (2.5,-0.9) to  [out=0,in=-90] (3.6,0);
    \draw [fill] (3.6,0) circle [radius=0.05];
    \end{scope}
    \end{scope}
    \end{scope}

\end{scope}
\end{tikzpicture}
}
\\
\scalebox{.75}[1.]
{
 \begin{tikzpicture}[xscale=1]
 \begin{scope}[shift={(-.4,0)}]

    \begin{scope}[xscale=-1]
    \begin{scope}[shift={(-3.2,0)}]

    \begin{scope}[shift={(-.5,0)}]
            \draw[thick,<-](0,0) to  [out=-90,in=180] (0.4,-0.8) to  [out=0,in=-90] (0.8,0);
             \draw [fill] (0.8,0) circle [radius=0.05];
              \node[above ] at (0.4,0) {$NP$};
               \node[above] at (0.4,-.6) {$\mbox{John}$};

    \end{scope}

    \begin{scope}[shift={(-.45,0)}]
             \draw[thick,<-](1.1,0) to  [out=-90,in=180] (2.5,-1.25) to  [out=0,in=-90] (3.9,0);
        % to [out=90,in=180] (4,0.5) to [out=0,in=90]

     %    \draw [fill] (2.75,-1.5) circle [radius=0.05];
    %    \draw[thick,<-](2.25,-1.5) to  [out=90,in=0] (1.75,-0.75) to  [out=180,in=-90] (1.5,0);
        \draw [fill] (3.9,0) circle [radius=0.05];
         %to [out=90,in=180] (4,2.5) to [out=0,in=90]

        %;

        %\draw[thick,<-](2,0) to  [out=-90,in=180] (2.5,-0.75) to  [out=0,in=-90] (3,0);

        \node[above] at (2.5,-.75) {$``\mbox{loves Mary}''$};

          \node[above] at (1.25,0) {$\overline{NP}$};
          \node[above] at (3.75,0) {$S$};
        %\node[above] at (4.8,0) {$S^\bot$};
        %\node[above] at (6.3,0) {$NP$};

        % \draw[thick,->](4.5,0) to  [out=-90,in=180] (7.25,-2.5) to  [out=0,in=-90] (10,0);
        \draw[thick,->](1.4,0) to  [out=-90,in=180] (2.5,-0.9) to  [out=0,in=-90] (3.6,0);
        \draw [fill] (1.4,0) circle [radius=0.05];

    \end{scope}

    \end{scope}
    \end{scope}

    {
    \node at(4.4,-.5){$\Longrightarrow$};
    }
 \end{scope}

\begin{scope}[xscale=-1]
\begin{scope}[shift={(-13.7,0)}]
\begin{scope}[shift={(6.2,0)}]

    \begin{scope}[shift={(-1.5,0)}]
     \draw[thick,-](0,0) to  [out=-90,in=180] (0.4,-0.8) to  [out=0,in=-90] (0.8,0);
          \node[above] at (0.4,-.6) {$\mbox{John}$};

   % \draw[thick,->](0,0) -- (0,-1.5);
%    \draw[thick,-](0.5,0) -- (0.5,-1.5);
%      \node[below ] at (0.3,-1.5) {$NP$};
%      \draw [fill] (0.5,-1.5) circle [radius=0.05];
    \end{scope}

    \draw[thick](.0,0) to  [out=90,in=0] (-.35,.5) to  [out=180,in=90] (-.7,0);
    \draw[thick](0.3,0) to  [out=90,in=0] (-.75,1) to  [out=180,in=90] (-1.5,0);

    \begin{scope}[shift={(-1.1,0)}]

     \draw[thick,-](1.1,0) to  [out=-90,in=180] (2.5,-1.25) to  [out=0,in=-90] (3.9,0);
    % to [out=90,in=180] (4,0.5) to [out=0,in=90]

 %    \draw [fill] (2.75,-1.5) circle [radius=0.05];
%    \draw[thick,<-](2.25,-1.5) to  [out=90,in=0] (1.75,-0.75) to  [out=180,in=-90] (1.5,0);

     %to [out=90,in=180] (4,2.5) to [out=0,in=90]

    %;

    %\draw[thick,<-](2,0) to  [out=-90,in=180] (2.5,-0.75) to  [out=0,in=-90] (3,0);

    \node[above] at (2.5,-.75) {$``\mbox{loves Mary}''$};

      \node[above] at (3.7,0) {$S$};
    %\node[above] at (4.8,0) {$S^\bot$};
    %\node[above] at (6.3,0) {$NP$};

    % \draw[thick,->](4.5,0) to  [out=-90,in=180] (7.25,-2.5) to  [out=0,in=-90] (10,0);
    \draw[thick,->](1.4,0) to  [out=-90,in=180] (2.5,-0.9) to  [out=0,in=-90] (3.6,0);
    \draw [fill] (3.9,0) circle [radius=0.05];
    \end{scope}
\end{scope}
\end{scope}
\end{scope}

\begin{scope}[shift={(.2,0)}]
    \node [font=\large] at(9.5,-.5) {
    $=$
    };

    \begin{scope}[xscale=-1]
\begin{scope}[shift={(-24,0)}]
    \begin{scope}[shift={(9.5,.25)}]

        \begin{scope}
         \draw[thick,->](.8,-.5) to  [out=-90,in=180] (2.5,-1.25) to  [out=0,in=-90] (4.2,-.5);
        % to [out=90,in=180] (4,0.5) to [out=0,in=90]

     %    \draw [fill] (2.75,-1.5) circle [radius=0.05];
    %    \draw[thick,<-](2.25,-1.5) to  [out=90,in=0] (1.75,-0.75) to  [out=180,in=-90] (1.5,0);
        \draw [fill] (.8,-.5) circle [radius=0.05];
         %to [out=90,in=180] (4,2.5) to [out=0,in=90]

        %;

        %\draw[thick,<-](2,0) to  [out=-90,in=180] (2.5,-0.75) to  [out=0,in=-90] (3,0);

        \node[above] at (2.5,-1.25) {$``\mbox{John loves Mary}''$};

        %\node[above] at (4.8,0) {$S^\bot$};
        %\node[above] at (6.3,0) {$NP$};

        % \draw[thick,->](4.5,0) to  [out=-90,in=180] (7.25,-2.5) to  [out=0,in=-90] (10,0);

        \end{scope}

    \end{scope}
    \end{scope}
    \end{scope}
\end{scope}
\end{tikzpicture}
}
\caption{Derivation}
\label{cowordism_derivation}
\end{subfigure}
%%%%%%%%%%%%%%%%%%%%%%%%%%%%%%%%%%%%%%%%%%%%%
\caption{Tensor grammar example}
 \end{figure}
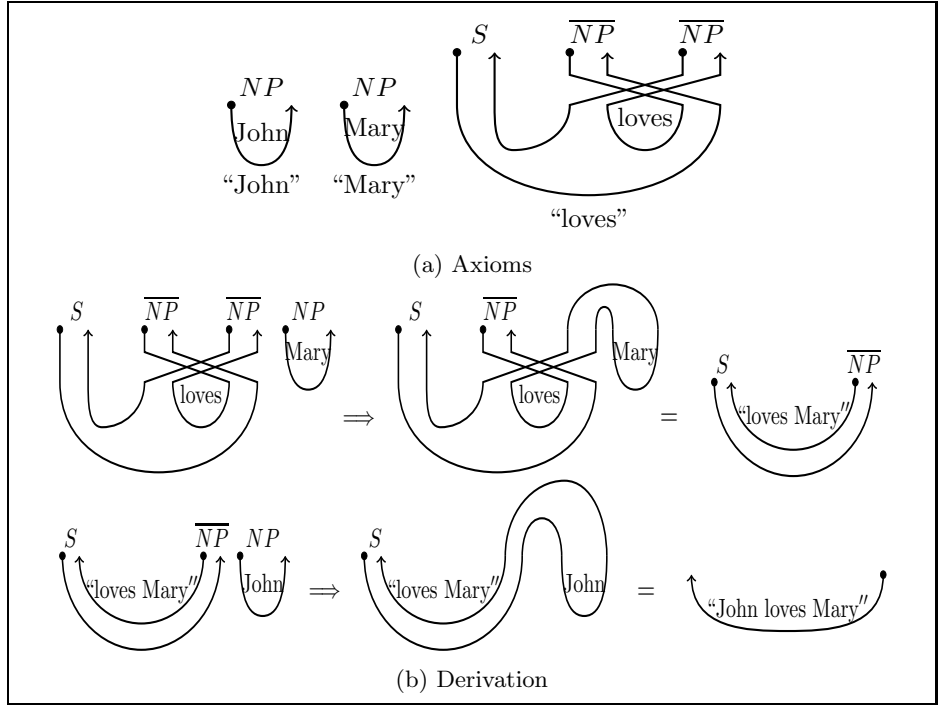

Let the terminal alphabet be $$\{\rm{John},\rm{Mary},\rm{loves}\}$$ and consider the tensor grammar defined by the following axioms
$$
[{\rm{loves}}]_l^r\cdot\delta^k_j\cdot\delta^i_s\vdash
S_i^j,\overline{NP}^l_k,\overline{NP}_r^s,
\quad[{\rm{Mary}}]^i_j\vdash NP^j_i,
\quad
[{\rm{John}}]^i_j\vdash NP^j_i.
$$
If we agree that commas in a sequent are ``hypocrisies'' for $\wp$ and $\wp$ is a symmetrized implication, then the above axioms could be understood as typing declarations
$$\mbox{``John''},\mbox{``Mary''}:NP,\quad \mbox{``loves''}:NP\multimap NP\multimap S.$$ Geometric representation is shown in Figure \ref{John_loves_Mary}.

With this we can derive the typing judgement
$$
\cfrac
            {
            \cfrac
                {
               [{\rm{loves}}]_l^r\cdot\delta^k_j\cdot\delta^i_s\vdash
                S_i^j,\overline{NP}^l_k,\overline{NP}_r^s
                \quad
                [{\rm{Mary}}]_{r}^{s}\vdash NP_{s}^{r}
                }
                {
                [{\rm{loves}}]_l^r\cdot\delta^k_j\cdot\delta^i_s
                \cdot
                [{\rm{Mary}}]_{r}^{s}
                \vdash
                S_i^j,\overline{NP}^l_k
                }({\rm{Cut}})
            }
            {
                \delta^k_j\cdot[{\rm{loves}}\mbox{ }{\rm{Mary}}]_{l}^{i}\vdash
              S_i^j,\overline{NP}^l_k
            }(\equiv)
$$
and then, in a similar way,
$$
\cfrac
    {
     \delta^k_j\cdot[{\rm{loves}}\mbox{ }{\rm{Mary}}]_{l}^{i}\vdash
              S_i^j,\overline{NP}^l_k
              \quad
    [{\rm{John}}]^{l}_{k}\vdash NP^{k}_{l}
    }
    {
     \delta^k_j\cdot[{\rm{loves}}\mbox{ }{\rm{Mary}}]_{l}^{i}
    \cdot
    [{\rm{John}}]^{l}_{k}\vdash
    S_i^j,
    }({\rm{Cut}})
$$
which, after a straightforward computation yields the standard sentence ``John loves Mary''.

The two steps of the derivations are illustrated in Figure \ref{cowordism_derivation}. (It is implied that every individual string in the picture is read from left to right,  but concatenation of individual strings labeling an edge is in the order determined by the edge direction. When there are several terminal symbols in the string, we put it in quotation marks for clarity.)

\subsection{Elaborate example}\label{elaborate example}
We discuss a more elaborate example with adverbs, relativization and medial extraction. Geometric representation would occupy too much space, so we skip it.

Let us add  terminal symbols
$$\{\rm{leaves},\rm{madly},\rm{whom}\}$$ and  axioms
\be\label{leaves}
\delta_j^l\cdot[{\rm{leaves}}]_l^i\vdash S^j_i,\overline{NP}^k_l
\quad
\delta_j^t\cdot\delta_s^k\cdot\delta_l^r\cdot[{\rm{madly}}]_u^i\vdash
S_i^j,\overline{NP}_k^l,{NP}^s_r\otimes \overline{S}^u_t,
\ee
\be\label{whom}
\delta_w^i\cdot\delta_j^k\cdot\delta_u^t\cdot[{\rm{whom}}]_l^v\vdash
NP_i^j,\overline{NP}_k^l,{NP}^u_t\otimes \overline{S}^w_v.
\ee
Just as above, we can read the axioms as typing declarations
$$
\mbox{``leaves''}:NP\multimap S,\quad \mbox{``madly''}:(NP\multimap S)\multimap NP\multimap S,$$
$$
\mbox{``whom''}:(NP\multimap S)\multimap NP\multimap N.
$$

Using the axiom for ``loves'' that we already have and renaming indices to avoid collisions we derive
$$
    \cfrac
    {
        \begin{array}{c}
        \\
        \delta_{j'}^j\delta_k^{k'}\delta_{l'}^l[{\rm{madly}}]_i^{i'}
        \vdash
        S_{i'}^{j'},\overline{NP}_{k'}^{l'},{NP}^k_l\otimes \overline{S}^i_j
        \end{array}
        \cfrac
        {
            \delta_s^i\delta^k_j[{\rm{loves}}]_l^r
            \vdash
            S_i^j,\overline{NP}^l_k,\overline{NP}_r^s
        }
        {
            \delta_s^i\delta^k_j[{\rm{loves}}]_l^r
            \vdash
            S_i^j\wp\overline{NP}^l_k,\overline{NP}_r^s
        }
        (\wp)
    }
    {
    \cfrac
        {
        \delta_{j'}^j\delta_k^{k'}\delta_{l'}^l[{\rm{madly}}]_i^{i'}
        \cdot
        \delta_s^i\delta^k_j[{\rm{loves}}]_l^r
        \vdash
        S_{i'}^{j'},\overline{NP}_{k'}^{l'},\overline{NP}_r^s
        }
        {
        [{\rm{madly}}]_s^{i'}
        \cdot
        [{\rm{loves}}]_{l'}^r\cdot\delta^{k'}_{j'}
        \vdash
        S_{i'}^{j'},\overline{NP}_{k'}^{l'},\overline{NP}_r^s
        }
        (\equiv).
    }
    ({\rm{Cut}})
$$
Then we cut the conclusion with the axiom for ``John''
$$
\cfrac
{
[{\rm{John}}]_{k'}^{l'}\vdash NP^{k'}_{l'}
\quad
 [{\rm{madly}}]_s^{i'}
        \cdot
        [{\rm{loves}}]_{l'}^r\cdot\delta^{k'}_{j'}
        \vdash
        S_{i'}^{j'},\overline{NP}_{k'}^{l'},\overline{NP}_r^s
 }
{
[{\rm{John}}\mbox{ }{\rm{loves}}]^r_{j'}\cdot[{\rm{madly}}]^{i'}_s
\vdash
S_{i'}^{j'},\overline{NP}_r^s
}
({\rm{Cut}})
$$
(we skipped straightforward normalization steps), and rename indices in the conclusion to obtain
$$
[{\rm{John}}\mbox{ }{\rm{loves}}]^u_{v}\cdot[{\rm{madly}}]^{w}_t
\vdash
S_{w}^{v},\overline{NP}_u^t.
$$
Applying the $(\wp)$ rule and cutting with the axiom for ``whom'' we derive from the above the typing judgement
$$
\delta^k_j\cdot [{\rm{whom}\mbox{ }\rm{John}\mbox{ }\rm{loves}\mbox{ } \rm{madly}}]^i_l\vdash
NP_i^j,\overline{NP}_k^l.
$$
Cutting the result with axioms for ``Mary'' and ``leaves'' we can obtain the sentence
$$\rm{Mary}\mbox{ }\rm{whom}\mbox{ }\rm{John}\mbox{ }\rm{loves}\mbox{ }\rm{madly}\mbox{ }\rm{leaves}.$$

\subsection{Inadequacy}\label{Inadequacy}
(The  discussion in the remainder of this section presupposes some familiarity with ACG \cite{deGroote} and Lambek grammars \cite{Lambek}. They will be considered in a greater detail in subsequent sections.)

It can be observed that representation of data in tensor grammars is rather non-economical.

For example, if we encode a transitive verb  in  a  term of type $NP\multimap NP\multimap S$, as is customary in categorial grammars, then we need to keep in memory {\it three} strings (elementary terms), although only one  of them is nonempty. This contrasts with Lambek grammars, where English transitive verb, for example, is customarily  represented as a {\it single} string of type $(NP \backslash S)/ NP$. Of course, the contrast becomes even more striking when we consider more complicated grammatical categories.

The reason is that in tensor grammars (as well as in ACG) all information about positioning words in a sentence is contained in the corresponding term, while in Lambek grammars this information is stored  in the corresponding type, once for all  type elements. Apparently it would be desirable to have some finer structure on tensor types allowing similar economy.

These considerations become even more relevant when we consider complex linguistic phenomena, such as {\it  coordination}. It was very convincingly explained in \cite{Moot_comparing} that, in contrast to Lambek grammars, ACG are in a certain sense {\it inadequate} for modeling (non-constituent) coordination, at least, in the most direct, ``naive'' approach.  This analysis applies to tensor grammars equally well.

In fact, tensor representation makes this ``inadequacy'' very transparent. Consider, as a very simple example, the sentence
\be\label{sentence_with_and}
\mbox{John loves  and Jim hates Mary.}
\ee
The elements that are coordinated are the strings ``John loves'' and ``Jim hates''. In a Lambek grammar both elements will be typed as
\be\label{coordinated_elements}
{\rm{John\mbox{ } loves}},{\rm{Jim \mbox{ }hates}}:S/NP.
\ee
The conjunction ``and'' will be typed as
\be\label{and}
{\rm{and}}:(X\backslash X)/ X, \mbox{ where }X=S/NP,
\ee
and (\ref{sentence_with_and}) is easily derivable.

On the other hand, in an ACG or a tensor grammar the elements in (\ref{coordinated_elements})
 are customarily modeled as terms of type $NP\multimap S$. If we try to mimic (\ref{and}) then we need a coordinating operator of type
 \be\label{AND}
AND:(NP\multimap S)\multimap(NP\multimap S)\multimap(NP\multimap S),
\ee
which will combine two terms of type $NP\multimap S$ into one.

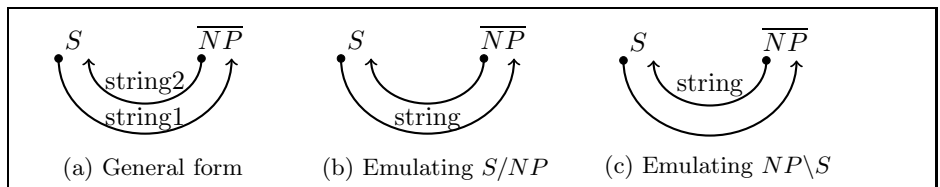
\begin{figure}[t]
\begin{subfigure}{.3\textwidth}
\centering
\begin{tikzpicture}
\begin{scope}[xscale=1]
%\begin{scope}[shift={(-24,0)}]
\draw[thick,<-](0,0) to  [out=-90,in=180] (.75,-.6) to  [out=0,in=-90](1.5,0);
\draw[thick,->](-.4,0) to  [out=-90,in=180] (.75,-1) to  [out=0,in=-90](1.9,0);
        \draw [fill] (-.4,0) circle [radius=0.05];
        \draw [fill] (1.5,0) circle [radius=0.05];
        \node [above] at (-0.2,0) {$S$};
        \node [above] at (1.75,0) {$\overline{NP}$};

        \node  at (0.75,-.8) {string1};
        \node [above] at (0.75,-.6) {string2};
\end{scope}
\end{tikzpicture}
\caption{General form}
\label{implication_general_form}
\end{subfigure}
\begin{subfigure}{.3\textwidth}
\centering
\begin{tikzpicture}
\begin{scope}[xscale=1]
%\begin{scope}[shift={(-24,0)}]
\draw[thick,<-](0,0) to  [out=-90,in=180] (.75,-.6) to  [out=0,in=-90](1.5,0);
\draw[thick,->](-.4,0) to  [out=-90,in=180] (.75,-1) to  [out=0,in=-90](1.9,0);
        \draw [fill] (-.4,0) circle [radius=0.05];
        \draw [fill] (1.5,0) circle [radius=0.05];
        \node [above] at (-0.2,0) {$S$};
        \node [above] at (1.75,0) {$\overline{NP}$};

        \node  at (0.75,-.8) {string};

\end{scope}
\end{tikzpicture}
\caption{Emulating $S/NP$}
\label{implication_/}
\end{subfigure}
\begin{subfigure}{.3\textwidth}
\centering
\begin{tikzpicture}
\begin{scope}[xscale=1]
%\begin{scope}[shift={(-24,0)}]
\draw[thick,<-](0,0) to  [out=-90,in=180] (.75,-.6) to  [out=0,in=-90](1.5,0);
\draw[thick,->](-.4,0) to  [out=-90,in=180] (.75,-1) to  [out=0,in=-90](1.9,0);
        \draw [fill] (-.4,0) circle [radius=0.05];
        \draw [fill] (1.5,0) circle [radius=0.05];
        \node [above] at (-0.2,0) {$S$};
        \node [above] at (1.75,0) {$\overline{NP}$};

        \node [above] at (0.75,-.6) {string
        };
\end{scope}

\end{tikzpicture}
\caption{Emulating $NP\backslash S$}
\label{implication_setmnus}
\end{subfigure}
\caption{Terms of type $NP\multimap S$}
\end{figure}

But the tensor type $NP\multimap S$ is not elementary, i.e. its elements are not strings. Rather, they are ordered pairs of strings (see Figure \ref{implication_general_form}). Coordinating them means gluing two pairs of strings into one pair, and there are too many ways of doing this... But apparently none of these ways corresponds to actual coordination occurring in the language.

Staying with the toy grammar  of the preceding subsection, we can generate at least three different kinds of terms in the type $S^i_j\wp \overline{NP}^k_l$, namely:
$$[\mbox{John loves}]_i^l\cdot\delta_k^j,\quad \delta_i^l\cdot[\mbox{loves Mary}]_k^j,\quad[\mbox{John loves}]_i^l\cdot[\mbox{madly}]_k^j.
$$
Obviously, all three terms have direct linguistic meaning. At the same time they cannot be coordinated with each other.

 On the other hand, we can note that (nonempty) strings in terms of the first kind (Figure \ref{implication_/}) correspond to the Lambek grammar type $S/NP$  and can be coordinated with each other. Similarly,  strings in the second kind (Figure \ref{implication_setmnus}) corresponds to the  type $NP\backslash S$ and can  be coordinated with each other equally well. As for the third kind, its elements cannot be represented as  strings and  cannot be coordinated in a simple way. But a coordinating (i.e. gluing) operator of form (\ref{AND}) will not distinguish between incompatible kinds and thus will inevitably produce ungrammatical sentences (something like ``and Jim hates loves Mary'').

 (Some more complicated analysis will show that treatment of adverbs, such as ``madly'' in the preceding section, tends to be inadequate in a commutative grammar and lead to overgeneration.)

This simple analysis shows once again that the structure of tensor types (or linear implicational types in the case of ACG) is too coarse, at least for simple intuitive modeling of non-constituent coordination. It seems clear that we need type constructors capable of emulating Lambek calculus.

\subsection{Towards Lambek  types}
The tensor representation makes very transparent how the ``non-commutative'' types of Lambek calculus look inside ``commutative''  types of tensor grammars and ACG.

Indeed, let $A$, $B$ be types of Lambek grammar. Their elements are strings, so they can be emulated as tensor types of valency (1,1), say $A_k^l$ and $B^i_j$. Then elements of the complex Lambek grammar types $B/A$  and
$A\backslash B$
 can be represented as elements of the tensor type $B^i_j\wp \overline{A}^k_l$ of the form, respectively,
\be\label{towards lambek}
[u]_l^i\cdot\delta_k^j \mbox{ and }[u]_k^j\cdot\delta_l^i,
\ee
as Figures \ref{implication_/}, \ref{implication_setmnus} suggest.
 It is easily computed that elements of the first form act (by means of the Cut rule) on elements of $A_k^l$ by multiplication (concatenation) on the left, and elements of the second form, by multiplication on the right.

 The two formats in (\ref{towards lambek}) identify two {\it subtypes} of the implicational tensor type that correspond to two implicational types of Lambek calculus. Note that, if we restrict to either of  these subtypes, then only two of the four available type indices become relevant, while two other indices are automatically connected with an empty string (Kronecker delta). This suggests that the new type constructors for emulating Lambek calculus should {\it bind indices in types}.

 In the next section we accurately discuss tensor representation of ACG. After that we develop {\it extended tensor type calculus} for embedding Lambek grammars basing on the above considerations.

\section{Representing abstract categorial grammars}\label{ACG section}
In this section we assume that the reader is familiar with basic notions of  $\lambda$-calculus, see \cite{Barendregt} for a reference.

\subsection{Linear $\lambda$-calculus}
Given a  set $X$ of {\it variables} and a   set $C$ of  {\it constants}, with $C\cap X=\emptyset$, the set $\Lambda(X,C)$ of {\it linear $\lambda$-terms} is defined by the following.
\begin{itemize}
\item Any $a\in X\cup C$ is in $\Lambda(X,C)$;
\item if $t,s\in\Lambda(X,C)$ are linear $\lambda$-terms whose sets of free variables are disjoint then $(t\cdot s)\in\Lambda(X,C)$;
\item if $t\in\Lambda(X,C)$, and $x\in X$ occurs freely in $t$ exactly once then $(\lambda x.t)\in \Lambda(X,C)$.
\end{itemize}

We use  common notational conventions such as omitting dots and outermost brackets and writing  iterated applications as
\be\label{iterated application}
(tsk)=(ts)k.
\ee

Given a  set $P$ of {\it propositional symbols} or {\it atomic types}, the set $Tp_\multimap(P)$ of  {\it linear implicational types} over $P$  is defined by induction.
\begin{itemize}
\item Any $A\in P$ is in $Tp_\multimap(P)$;
\item if $A,B\in Tp_\multimap(P)$, then $(A\multimap B)\in Tp_\multimap(N)$.
\end{itemize}
Since we are not going to discuss any non-linear fragment, the title ``linear'' in the context of $\lambda$-calculus will usually be omitted.

A {\it $\lambda$-typing assumption} is an expression of the form $x:A$, where $x\in X$  and $A\in Tp_\multimap(P)$. A (linear) {\it $\lambda$-typing context} is a finite set
$$
x_{1}:A_{1},\ldots,x_{n}:A_{n}
$$
 of typing assumptions, where $x_1,\ldots,x_n\in X$ are pairwise distinct.

A (linear) {\it $\lambda$-typing judgement} is an expression of the form
$\Gamma\vdash t:A$, where $t\in \Lambda(X,C)$, $A\in Tp_\multimap(P)$, and $\Gamma$ is a typing context.
We will say that two $\lambda$-typing judgements
$$\Gamma\vdash t_i:A,~i=1,2,$$
are $\beta\eta$-equivalent if the terms $t_1,t_2$ are $\beta\eta$-equivalent.

 $\lambda$-Typing judgements are derived from the following type inference rules, which happen to be the natural deduction rules for {\it implicational linear logic} ({\bf ILL}) decorated with $\lambda$-terms:
$${x:A\vdash x:A}~(\mbox{Id}),$$
$$\frac{\Gamma\vdash t:A\multimap B\quad\Gamma'\vdash s:A}{\Gamma,\Gamma'\vdash (ts):B}~(\multimap\rm{E}), \quad\frac{\Gamma,x:A\vdash t:B}{\Gamma\vdash (\lambda x.t):A\multimap B}~(\multimap\rm{I}).$$

A {\it $\lambda$-signature}  $\Sigma$ is a triple $\Sigma=(P,C,\mathfrak{T})$, where $P$ is a finite set of atomic types, $C$ is a finite set of constants and $\mathfrak{T}$ is a function assigning to each constant $c\in C$ a linear implicational type $\mathfrak{T}(c)\in Tp(P)$.
Typing judgements of the form
$\vdash c:\mathfrak{T}(c)
$,
 where $c\in C$, are called {\it signature axioms} of $\Sigma$.

Given a $\lambda$-signature  $\Sigma$, we say that a typing judgement $\Gamma\vdash t:A$ is {\it derivable in} $\Sigma$ if it is derivable from axioms of $\Sigma$ by rules of linear $\lambda$-calculus. We write in this case $\Gamma\vdash_\Sigma t:A$.

The following theorems are standard and proven by induction on derivation.
\bp\label{substitution}
For any $\lambda$-signature $\Sigma$, the set of derivable typing judgements is closed under the {\it substitution}
rule
$$\frac{\Gamma\vdash t:A\quad x:A,\Gamma'\vdash s:B}{\Gamma,\Gamma'\vdash s[x:=t]:B}~(\rm{Subst}).\quad \Box$$
\ep
\bl[``Deduction theorem'' for $\lambda$-calculus]\label{deduction theorem lambda}
 Let $\Sigma$ be a $\lambda$-signature.

 Let $$\Xi =\{\vdash t_{(1)}:A_{(1)},\ldots,\vdash t_{(n)}:A_{(n)}\} $$
be a finite multiset of $\lambda$-typing judgements.

A typing judgement  ${\Gamma}\vdash t:A $ is derivable in $\Sigma$ from elements of $\Xi$ using each element exactly once  iff there exists a term $t'$ and a
typing judgement $\sigma'$ of the form
 $$x_1:A_{(1)},\ldots,x_n:A_{(n)},\Gamma\vdash t': A$$
 derivable
in $\Sigma$   such that
$$t'[x_1:=t_{(1)},\ldots, x_n:=t_{(n)}]\sim_\beta t.\quad \Box$$
\el

\subsection{Translating $\lambda$-signatures}
Linear $\lambda$-calculus and {\bf ILL} translate to {\bf MLL} as conservative fragments. This translation can be lifted to the level of {\bf TTC}, provided that we decide how to assign valencies to atomic types. This lifting is very formal and does not use indices in any essential way.

So, let $\Sigma =(P,C,\mathfrak{T})$ be a $\lambda$-signature, and assume that every element $p\in P$ is assigned a valency $v(p)\in{\bf N}^2$.

Then we get  an embedding $$tr:Tp_\multimap(P)\to Tp_\otimes(P)$$ of implicational types to tensor type symbols  by the following induction:
\be\label{lambda types translation}
tr(p)=p\mbox{ for } p\in P,\quad tr(A\multimap B)=tr(B)\wp\overline{tr(A)}.
\ee

By abuse of notation, in the following we will denote an implicational type $A$ and the corresponding tensor type symbol $tr(A)$ the same.
Also, for a tensor type $A$ we say that $A$ is a {\it positive implicational type} if it is a translation of an implicational type, and we say that it is a {\it negative implicational type} if its dual is a positive implicational type.

Now let a finite alphabet $T$ of terminal symbols be given, and
assume that each axiom $\alpha$ of $\Sigma$ is assigned a tensor typing judgement over $P$ and $T$, its {\it tensor translation}, in such a way that if $\alpha$ has the form $\vdash c:A$, then its tensor translation has the form $\widetilde \alpha\vdash A_I^J$ (where $\widetilde \alpha$ is some tensor term).

Let $\Xi$ be the set of all tensor translations of axioms of $\Sigma$.
This defines a tensor signature $\widetilde \Sigma=(P,T,\Xi)$, the {\it tensor translation} of $\Sigma$.
We are going to extend tensor translation from axioms to all  $\lambda$-typing judgements derivable in $\Sigma$.

A $\lambda$-typing judgement $\sigma$ of the form
$$x_1:A_{(1)},\ldots,x_n:A_{(n)}\vdash t:A$$ derivable in $\Sigma$
will be translated to tensor typing judgement $Tr(\sigma)$ with links
$$\overline{A_{(1)}},\ldots,\overline{A_{(n)}},{A},$$
derivable in $\widetilde \Sigma$.
Definition is by induction on derivation of $\sigma$  and simply repeats translation of {\bf ILL} to {\bf MLL} without any essential use of indices. So we describe it very briefly.

 The axiom $x:A\vdash x:A$ translates as
  $\delta_{\overline{I}J'}^{I'\overline{J}}\vdash A_{I'}^{J'},\overline{A}^{I}_{J}$, which is  easily derivable in {\bf TTC} (if $A$ is atomic, this is just an axiom).

 If $\sigma$ is obtained from a derivable typing judgement $\sigma'$   by the $(\multimap{\rm{I}})$ rule introducing the type $A\multimap B$, then we translate it to
  %a $\lambda$-typing judgement derivable in $\Sigma$
  the typing judgement obtained from $Tr(\sigma')$ by the $(\wp)$ rule introducing the link $B\wp\overline{A}$.

If $\sigma$ is obtained from some derivable typing judgments $\tau$ and $\rho$ with types $A$ and $A\multimap B$ respectively to the right of the turnstile by the $(\multimap{\rm{E}})$ rule, then $Tr(\sigma)$ is obtained by applying the $(\wp)^{-1}$ rule to $Tr(\rho)$ (in the only possible way) and then cutting the result with $Tr(\tau)$ on the links $A,\overline{A}$.

It will be convenient to have a special title for tensor typing judgements that look like translated $\lambda$-typing judgements.
Given a translation of implicational types to tensor types, we say that a tensor typing judgement is {\it polarized} if it has one occurrence of a positive implicational type symbol and all other links are negative  implicational type symbols.

\subsubsection{Properties of tensor translation}
We list some useful properties of $\lambda$-signature tensor translations. The statements below are, essentially, well-known theorems about plain translation of  {\bf ILL} to  {\bf MLL}.
\bl\label{substituition=cut}
In the setting as above, if we have $\lambda$-typing judgements $\rho$, $\tau$, respectively
$$\Gamma\vdash t:A,\quad
x:A,\Gamma'\vdash s:B,
$$
 derivable in the signature $\Sigma$, and a $\lambda$-typing judgement $\sigma$ is obtained from $\rho$, $\tau$ by the substitution rule, then $Tr(\sigma)$ is obtained from $Tr(\tau)$, $Tr(\rho)$ by cut on the links $A,\overline A$. $\Box$
 \el
  \bl\label{beta-eta-equivalence preservation}
   Tensor translations of $\beta\eta$-equivalent $\lambda$-typing judgements coincide. $\Box$
  \el
Lemma \ref{substituition=cut} is proven by induction on derivation, and Lemma \ref{beta-eta-equivalence preservation} by induction on the definition of $\beta\eta$-equivalence using Lemma \ref{substituition=cut}.

The following is a more technical result (proved by induction on derivation using Lemma \ref{substituition=cut}) that will be used in the sequel.
\bl\label{main for translation}
Assume that we have $m$ constants $c_1,\ldots, c_m$ ($m$ may be zero) and $n$ $\lambda$-terms $t_1,\ldots, t_n$, $n>0$,
and
let $\sigma$ be a $\lambda$-typing judgement of the form
$
\Gamma,x:G\multimap H\vdash  t:F,
$
where $$t=c_1(\ldots(c_m(x\cdot t_1\cdots t_n)\ldots),$$
derivable in $\Sigma$.

Then there exist typing judgements $\sigma_1,\sigma_2$, respectively
$$
\Gamma_1\vdash t_1:G,
\quad
y:H,\Gamma_2\vdash t' :F
$$
where $$t'=c_1(\ldots(c_m(y\cdot t_2\cdots t_n)\ldots),$$
derivable in $\Sigma$, such that $\sigma$ can be obtained from $\sigma_1,\sigma_2$ as
  $$
\cfrac
    {
       \cfrac
           {
           x:G\multimap H\vdash x:G\multimap H
           \quad
           \sigma_1
           }
           {
           x:G\multimap H,\Gamma_1\vdash x\cdot t_1:H
           }
           ~(\multimap{\rm{E}}
       )
    \quad
    \begin{array}{c}
       \\
       \sigma_2
    \end{array}
    }
    {
     \sigma
    }({\rm{Subst}}),
   $$
      and the tensor translation $Tr(\sigma)$ can be obtained from $Tr(\sigma_1)$, $Tr(\sigma_2)$
   by the $(\otimes)$ rule introducing the link $G\otimes\overline{H}$. $\Box$
\el

\subsubsection{Variations on conservativity of tensor translations}\label{variations on conservativity section}
In the case of no non-logical axioms, i.e. when we translate linear $\lambda$-calculus to {\bf TTC}, tensor translation of typing judgements is, up to $\beta\eta$-equivalence, a bijection with its image. This is, again, a theorem about plain {\bf MLL} where indices play no role, but we reproduce basic arguments, in order to reuse them for the special case of {\it string signatures}.

In the proposition below we collect some simple statements, essentially,  about plain {\bf MLL} proof-nets, which are proven by induction on cut-free derivations in {\bf TTC}.
\bp\label{proof-net properties}
Let $\sigma$ be  a typing judgement derivable in {\bf TTC}.
\begin{enumerate}[(i)]
\item If $\sigma$ has no $\wp$-links and is not an axiom, then $\sigma$ has a splitting link.
\item  If there is  a  link $Z=X\otimes Y$ splitting $\sigma$ into
typing judgements $\sigma_1$, $\sigma_2$ then both $\sigma_1$, $\sigma_2$ are derivable in {\bf TTC}.
\item If $\sigma$ is polarized and has a splitting link, then it is split into  polarized typing judgements. $\Box$
\end{enumerate}
\ep
\bl\label{inverse translation}
There exists an {\it inverse translation} map $Tr^{-1}(.)$ from polarized typing judgements derivable in {\bf TTC} to $\lambda$-typing judgements derivable in linear $\lambda$-calculus satisfying $Tr(Tr^{-1})(\sigma)=\sigma$.
\el
{\bf Proof}
The map is defined  by induction on the number of connectives occurring in $\sigma$, using Proposition \ref{proof-net properties}.

 If  $\sigma$ has a $\wp$-link, then it can be obtained from some $\sigma_0$ by the $(\wp)$-rule and the inverse translation $Tr^{-1}(\sigma)$ is obtained from $Tr^{-1}(\sigma_0)$ by the $(\multimap{\rm{I}})$ rule.

Otherwise $\sigma$ has a splitting link $G\otimes\overline{H}$ ($=\overline{G\multimap H}$) which splits it
into polarized typing judgements $\sigma_1$, $\sigma_2$, which already have inverse translations  of the forms
$$\Gamma_1\vdash t:G,\quad y:H,\Gamma_2\vdash s:C$$
respectively (where $C$ is the unique positive implicational link of $\sigma$).

  We define $Tr^{-1}(\sigma)$ as the conclusion of the following $\lambda$-calculus derivation
   $$
\cfrac
    {
       \cfrac
       {
       x:G\multimap H\vdash x:G\multimap H
       \quad
       Tr^{-1}(\sigma_1)
       }
       {
       x:G\multimap H,\Gamma_1\vdash x\cdot t:H}~(\multimap{\rm{E}}
       )
    \quad
    \begin{array}{c}
       \\
       Tr^{-1}(\sigma_2)
    \end{array}
    }
    {
     x:G\multimap H,\Gamma_1,\Gamma_2\vdash s[y:=x\cdot t]:{{F}},
    }({\rm{Subst}}).
   $$
   One shows by induction on  derivation of $\sigma$ that this definition does not depend on the choice of the splitting
 link.
  $\Box$

 It can be shown that the above defined map is indeed, up to $\beta\eta$-equivalence, an  inverse of tensor translation, which means that,  in the axiom-free case,  tensor translation is a conservative {\it embedding}. For a general $\lambda$-signature we can say only that tensor translations are conservative.
 \bp\label{surjectivity of tensor translation}
 Let $\Sigma$ be a  $\lambda$-signature translated to a tensor signature $\widetilde{\Sigma}$. Then for any tensor typing judgement $\pi$ derivable in $\widetilde\Sigma$ there exists a $\lambda$-typing judgement $\sigma$ derivable in $\Sigma$ whose translation is $\pi$.
 \ep
 {\bf Proof} This  follows immediately from Lemmas \ref{deduction theorem}, \ref{deduction theorem lambda} (``Deduction theorems''), Lemma \ref{inverse translation} above and Lemma \ref{substituition=cut} (that substitution translates to cut). $\Box$

 Below we discuss {\it string signatures}, for which Lemma \ref{inverse translation} can be generalized and tensor translation is indeed an embedding, as in the axiom-free case.

\subsection{String signatures}
Let $T$ be a finite alphabet.

The {\it $\lambda$-string signature} $Str_{\lambda,T}$ {\it over} $T$ is the linear $\lambda$-signature with a single atomic type $O$, the alphabet $T$ as the set of constants and the  typing assignment
$$\mathfrak{T}(c)=O\multimap O\mbox{ }\forall c\in T.$$

We denote the implicational type $O\multimap O$ as $str$. This is the type of strings.
Any word $w=a_1\ldots a_n$ in the alphabet $T$ can be represented as the  term
$$
\rho(w)=(\lambda x.a_1(\ldots(a_n(x))\ldots)),
$$
for which we have
$$\vdash_{Str_{\lambda,T}}\rho(w):str.$$

 We define the {\it tensor string signature} $Str_{\otimes,T}$ over $T$ as  the tensor signature  with the single positive atomic type
$O$ of valency $v(O)=(1,0)$, the alphabet $T$ of terminal symbols and axioms
$$[c]_i^j\vdash  O^i\wp\overline{O}_j\mbox{ }\forall c\in T.$$
We denote the set of tensor type symbols
generated by $\{O,\overline{O}\}$ as $Tp_\otimes(O)$ and
we identify $Tp_\multimap(O)$ with a subset of $Tp_\otimes(O)$  using standard translation  (\ref{lambda types translation}) of implicational types to tensor types. In particular, we denote $str^i_j=O^i\wp\overline{O}_j$.

The signature $Str_{\otimes,T}$ is, obviously, a tensor translation of $Str_{\lambda,T}$ where each axiom
$\vdash c:str$ of  $Str_{\lambda,T}$ is translated as $[c]_i^j\vdash str^i_j$.
We are going to show now that this translation is an embedding.

\subsubsection{Inverse translation  of  tensor string signature}
Given a terminal alphabet $T$, the {\it big string tensor signature} $STR_{\otimes,T }$ over $T$
is defined by the single positive atomic type symbol $O$ of valency $(1,0)$ and the infinite set of axioms
$$
\{[w]^j_i\vdash O^i,\overline{O}_j~|w\in T^*\}.
$$
Obviously $STR_{\otimes,T}$ and $Str_{\otimes,T}$ have identical sets of derivable typing judgements. It is also easy to prove that $STR_{\otimes,T}$ is {\it cut-free}.

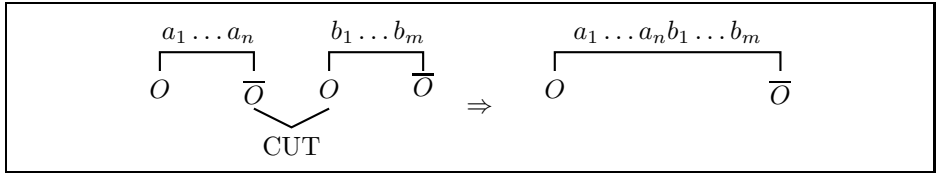
\begin{figure}
\centering
\begin{tikzpicture}[yscale=1]
  \node at (-.25,0)[below]{$O$};
  \node at (1,0)[below]{$\overline{O}$};
  \draw[thick](-.25,.0)--(-.25,.25)--(1,.25)--(1,0);
  \node at (.4,.25)[above]{$a_1\ldots a_n$};
  \draw[thick](1,-.5)--(1.5,-.75)--(2,-.5);

  \begin{scope}[shift={(2.25,0)}]
  \node at (-.25,0)[below]{$O$};
  \node at (1,0.1)[below]{$\overline{O}$};
  \draw[thick](-.25,.0)--(-.25,.25)--(1,.25)--(1,0);

  \node at (.4,.25)[above]{$b_1\ldots b_m$};
  \end{scope}
  \node at (1.5,-.75)[below]{CUT};

  \node at (4.,-.5) {$\Rightarrow$};

  \begin{scope}[shift={(5,0)}]
  \node at (0,0)[below]{$O$};
  \node at (3,0)[below]{$\overline{O}$};
  \draw[thick](0,.0)--(0,.25)--(3,.25)--(3,0);
  \node at (1.5,.25)[above]{$a_1\ldots a_nb_1\ldots b_m$};

  \end{scope}
\end{tikzpicture}
\caption{$STR_{\otimes,T}$ axiom links}
\label{STR links}
\end{figure}

(In fact, axioms of $STR_{\otimes,T}$ can be seen as most usual {\bf MLL} proof-net axiom links, only decorated with words. Then all derivable typing judgements become usual proof-nets with a slightly fancier cut-reduction rule for axiom links shown in Figure \ref{STR links}.
Thus, reasoning about typing judgements derivable in $STR_{\otimes,T}$ repeats reasoning about {\bf MLL} proof-nets, possibly with tiny additional details at the axiom level.)

In particular, we  have by induction on cut-free derivations in $STR_{\otimes, T}$ that Proposition \ref{proof-net properties} of Section \ref{variations on conservativity section} holds when {\bf TTC} is replaced by $STR_{\otimes, T}$ (hence, by $Str_{\otimes, T}$). This allows us defining inverse translation from polarized typing judgements of $Str_{\otimes,T}$ by induction the number of connectives occurring in a typing judgement, just as for pure {\bf TTC}.
\bl\label{inverse translation for strings}
There exists an {\it inverse translation} map $Tr^{-1}(.)$ from polarized typing judgements derivable in $STR_{\otimes,T}$ to $\lambda$-typing judgements derivable in $Str_{\lambda,T}$ satisfying $Tr(Tr^{-1})(\sigma)=\sigma$.
\el
{\bf Proof} Same as  Lemma \ref{inverse translation} except that for the basis of induction we have to consider all axioms of $STR_{\otimes,T}$.

If $\sigma$ is an axiom $[c_1\ldots c_m]^j_i\vdash O^i,\overline{O}_j$ (where $m$ is possibly zero) we take as $Tr^{-1}(\sigma)$
the typing judgement
 $x:O\vdash c_1\cdot(\ldots (c_mx)\ldots)$. $\Box$
 \bl
 If $\sigma$ is a $\lambda$-typing judgement of the form
 $\Gamma\vdash t:C$ derivable in $Str_{T,\otimes}$, where $t$ is a $\beta$-normal $\eta$-long $\lambda$-term, then
 $Tr^{-1}(Tr(\sigma))=\sigma$.
 \el
 {\bf Proof} Induction on $t$.

 The base case is when $t=c_1\cdot(\ldots (c_mx)\ldots)$ for some $c_1,\ldots c_m\in T$  (where $m$ may be zero). Then $Tr(\sigma)$ is an axiom of $STR_{\otimes,T}$ and the statement holds.

 If $t=\lambda x.t'$, then $t'$ is also $\beta$-normal $\eta$-long and the statement immediately follows from the induction hypothesis.

 The only other possibility is that $t=c_1\cdot(\ldots (c_m(xt_1\cdots t_n))\ldots)$,
 where $x$ is a variable, $n>0$ and $t_1,\ldots,t_n$ are $\beta$-normal $\eta$-long. Then the statements follows from the induction hypothesis and Lemma \ref{main for translation}. $\Box$ %%%%%%%%%%%%%%%%%%%%%%%%%%%%%%%%%%%%%%%%%%%%%%%%%

It follows that, for string signatures, tensor translations are {\it embeddings}.
\bc\label{conservativity for string signature}
For any polarized tensor typing judgement $\sigma$  derivable in
$Str_{\otimes,T}$ there exists a $\lambda$-typing judgement $\pi$ derivable  in $Str_{\lambda,T}$ whose tensor translation is $\sigma$. This $\pi$ is unique up to $\beta\eta$-equivalence of $\lambda$-terms. $\Box$
\ec

\subsection{Translating string ACG}
\subsubsection{String ACG}
Given two linear $\lambda$-signatures $\Sigma_i=(N_i,C_i,\mathfrak{T}_i)$, $i=1,2$, a {\it map of $\lambda$-signatures} $$\phi:\Sigma_1\to\Sigma_2$$ is a
pair $\phi=(F,G)$, where
\begin{itemize}
  \item $F:Tp_\multimap(\Sigma_1)\to  Tp_\multimap(\Sigma_2)$ is a function satisfying the homomorphism property $$F(A\multimap
      B)=F(A)\multimap F(B),$$
  \item $G:C_1\to\Lambda(X,C_2)$ is a function such that
       for any $c\in C_1$ it holds that $\vdash_{\Sigma_2}G(c):F(\mathfrak{T}(c))$.
\end{itemize}

The map $G$ above extends inductively to a map
$$G:\Lambda(X,C_1)\to\Lambda(X,C_2)$$
  by
  $$G(x)=x,\mbox{ }x\in X,$$
   $$G(ts)=(G(t)G(s)),\quad G(\lambda x.t)=(\lambda x.G(t)).$$

For economy of notation, we write $\phi(A)$ for $F(A)$ when $A\in Tp_\multimap(C_1)$, and we write $\phi(t)$ for $G(t)$ when $t\in\Lambda(X,C_1)$.

A {\it string abstract categorial grammar (string ACG)} $G$ is a tuple $G=(\Sigma,T, \phi,S)$, where
\begin{itemize}
  \item $\Sigma$, the {\it abstract signature}, is a linear $\lambda$-signature;
    \item $\phi:\Sigma\to Str_{\lambda,T}$, the {\it lexicon}, is a map of $\lambda$-signatures;
  \item $S$, the {\it  sentence type}, is an atomic type of $\Sigma$ with $\phi(S)=O\multimap O$.
\end{itemize}

The {\it string language} $L(G)$ generated by a string ACG $G$ is the set of words
$$L(G)=\{w\in T^*|\exists t~\rho(w)\sim_{\beta\eta}\phi(t),~\vdash_\Sigma t:S\}.$$

\subsubsection{Translation}
Now let a string ACG $G=(\Sigma,T,\phi,S)$  be given.

Let $P$ be the set of atomic types of $\Sigma$.

Let  $p\in P$. We have an implicational type $\phi(p)$ of $Str_{\lambda,T}$ and its tensor translation (which we denote by the same expression) in ${Str_{\otimes,T}}$.  We assign to $p$  the valency $v(p)=v(\phi(p))$ and consider it as a positive atomic type symbol for a new tensor signature. In this way we get the new set $Tp_\otimes(P)$ of tensor type symbols. We identify implicational types $Tp_\multimap(P)$ with a subset of $Tp_\otimes(P)$ using (\ref{lambda types translation}), as always.

The map $\phi:Tp_\multimap(P)\to Tp_\multimap(O)$ extends to a map
\be\label{phi for tensor}
\phi:Tp_\otimes(P)\to Tp_{\otimes}(O)
\ee
in the obvious way, setting
$$\phi(\overline{A})=\overline{\phi(A)},\quad \phi(A\otimes B)=\phi(A)\otimes\phi(B).$$
Note that  map  (\ref{phi for tensor}) preserves valencies of tensor types.

In order to translate $\Sigma$ to a tensor signature we need to translate
 its axioms.
 Roughly speaking,
  an axiom $\alpha$ of the form
 $\vdash c:A$
  should become a tensor typing judgement of the form $\widetilde{\alpha}\vdash A$, but this is not quite accurate. Indeed, we identified $A$ only with a tensor type {\it symbol}. In order to turn it into a legitimate tensor type, which can be written in a tensor typing judgement, we need to write it with indices, as $A^{i_1\ldots i_k}_{j_1\ldots j_n}$, where $v(A)=(k,n)$.

 Now, we have the $\lambda$-typing judgement $\phi(\alpha)$ of the form
$\vdash\phi(c):\phi(A)$,
which is derivable in $Str_{\lambda,T}$,
and we have its  translation $Tr(\alpha)$ in ${Str_{\otimes,T}}$ of the form
$$
\widetilde{\phi(\alpha)}\vdash \phi(A)^{i_1\ldots i_k}_{j_1\ldots j_n},
$$
where $FSub(\widetilde{\phi(\alpha)})=\{i_1,\ldots,i_k\}$, $FSup(\widetilde{\phi(\alpha)})=\{j_1,\ldots,j_n\}$
(remember that  $v(A)=v(\phi(A))$).

We put $\widetilde{\alpha}=\widetilde{\phi(\alpha)}$ and define the  tensor translation $Tr(\alpha)$ of the axiom $\alpha$
as the tensor typing judgement
$$
 \widetilde{\phi(\alpha)}\vdash A^{i_1\ldots i_k}_{j_1\ldots j_n},
 $$
 which is well-defined.

 We define the  tensor  signature $\widetilde \Sigma$ over $P$ and $T$ by taking as axioms all translations of axioms of $\Sigma$.
 \bp\label{tensor translation commutes with phi}
    Let a $\lambda$-typing judgement $\sigma$ of the form
    $$x_1:A_{(1)},\ldots,x_n:A_{(n)}\vdash t:A$$ be    derivable in $\Sigma$, and
    let its tensor translation $Tr(\sigma)$  in $\widetilde \Sigma$ be
    $$\widetilde{\sigma}\vdash (\overline{A_{(1)}})^{I_1}_{J_1},\ldots,(\overline{A_{(n)}})^{I_n}_{J_n},A^{I}_{J}$$
     Then the $\lambda$-typing judgement $\phi(\sigma)$ defined as
     $$x_1:\phi(A_{(1)}),\ldots,x_n:\phi(A_{(n)})\vdash\phi( t):\phi(A),$$
     which is derivable in $Str_{\lambda,T}$, translates  in ${Str_{\otimes,T}}$ as
     $$\widetilde{\sigma}\vdash (\overline{\phi(A_{(1)})})^{I_1}_{J_1},\ldots,(\overline{\phi(A_{(n)})})^{I_n}_{J_n},(\phi(A))^{I}_{J}.$$
       \ep
{\bf Proof} Straightforward induction on derivation of $\sigma$. $\Box$

 We define a tensor grammar $\widetilde G$ as $\widetilde G=(\widetilde \Sigma, S)$. We say that $\widetilde G$ is the {\it tensor representation} of $G$.

 Proposition \ref{surjectivity of tensor translation} (on conservativity of tensor translations), Proposition \ref{tensor translation commutes with phi} above and Corollary \ref{conservativity for string signature} (that tensor translations of string signatures are embeddings) immediately yield the following.
 \bt\label{conservativity for ACG}
 In the setting as above, tensor translation of $\Sigma$ to $\widetilde\Sigma$ is an embedding on the string level.

  A polarized tensor  typing judgement $\pi$ of the form $t\vdash \Gamma$ is derivable in $\widetilde\Sigma$ if it is the translation of some $\lambda$-typing judgement $\sigma$ derivable in $\Sigma$ such that the string level image,  $\lambda$-typing judgement $\phi(\sigma)$, in $Str_{\lambda,T}$ translates to the tensor typing judgment $t\vdash\phi(\Gamma)$ in $Str_{\otimes, T}$.

  Moreover, if there are two $\lambda$-typing judgements $\sigma_1$, $\sigma_2$ translating to $\pi$, then their string level images $\phi(\sigma_1)$, $\phi(\sigma_2)$ coincide up to $\beta\eta$-equivalence of $\lambda$-terms.  $\Box$
 \et
 \bc\label{ACG generates the same language}
 For any ACG $G$, the  string language $L(G)$ of $G$ coincides with the language $L(\widetilde G)$ of its tensor representation $\widetilde G$. $\Box$
 \ec

\subsubsection{Example}
For an illustration, consider the standard toy example, when the abstract signature has atomic types $NP$, $S$, and the lexicon $\phi$ maps them to the string type, i.e. $\phi(NP)=\phi(S)=O\multimap O$.

Assume that we have an intransitive English verb, for example, ``leaves''. The customary ACG type for that is $NP\multimap S$, and the corresponding  axiom
$$\vdash{\rm{leaves}}:NP\multimap S$$
in the abstract signature is mapped by the lexicon to
$$\phi({\rm{leaves}})=\lambda xy.x({\rm{leaves}}\cdot y).$$
We will translate this to an axiom for a tensor signature. The axiom will be as in (\ref{leaves}) in Section \ref{elaborate example}.

We identify $NP$ and $S$ as atomic tensor type symbols of valency $(1,1)$. The axiom for ``leaves'' should be of the form $t\vdash S\wp\overline{NP}$ decorated with four indices. For simplicity, we will omit the outermost $\wp$ connective and write the above as $t\vdash S,\overline{NP}$.

The derivation of $\phi({\rm{leaves}})\vdash\phi(NP\multimap S)$ in the string signature is as follows:
$$
\cfrac
    {
    \begin{array}{c}
                \\
        x:O\multimap O\vdash x:O\multimap O
    \end{array}
    ~
    \cfrac
        {\vdash{\rm{leaves}}:O\multimap O
        \;
        y:O\vdash y:O
        }
        {
        y:O\vdash{\rm{leaves}}\cdot y:O
        }(\multimap{\rm{E}})
    }
    {
    \cfrac
        {
        x:O\multimap O,y:O\vdash x({\rm{leaves}}\cdot y)
        }
        {
        \cfrac
            {
            x:O\multimap O\vdash \lambda y.x({\rm{leaves}}\cdot y):O\multimap O
            }
            {
            \vdash \lambda xy.x({\rm{leaves}}\cdot y):(O\multimap O)\multimap O\multimap O
            }(\multimap{\rm{I}}).
        }(\multimap{\rm{I}})
    }(\multimap{\rm{E}})
$$
We translate the assumptions
$$x:O\multimap O\vdash x:O\multimap O,\quad \vdash{\rm{leaves}}:O\multimap O,\quad y:O\vdash y:O$$
of the above derivation as, respectively
        $$\delta^{u j}_{iv}\vdash  O^v\wp\overline{O}_u,O^i\otimes\overline{O}_j,
        \quad
                        [{\rm{leaves}}]^t_u\vdash O^u\wp\overline{O}_t,
                        \quad
                        \delta^s_t\vdash O^t,\overline{O}_s.$$

The derivation itself translates as
$$
        \cfrac
            {
            \begin{array}{c}
                \\
                \\
            \cfrac
                {
                \delta^{u j}_{iv}\vdash  O^v\wp\overline{O}_u,O^i\otimes\overline{O}_j
                }
                {
                \delta^{u j}_{iv}\vdash  O^v,\overline{O}_u,O^i\otimes\overline{O}_j
                }(\wp^{-1})
            \end{array}
            \;
            \cfrac
                {
                \cfrac
                    {
                        [{\rm{leaves}}]^t_u\vdash O^u\wp\overline{O}_t
                    }
                    {
                        [{\rm{leaves}}]^t_u\vdash O^u,\overline{O}_t
                    }(\wp^{-1})
                \,
                \begin{array}{c}
                        \\
                        \delta^s_t\vdash O^t,\overline{O}_s
                \end{array}
                }
                {
                    \cfrac
                        {
                        [{\rm{leaves}}]^t_u\cdot\delta^s_t
                        \vdash
                        O^u,\overline{O}_s
                        }
                        {
                        [{\rm{leaves}}]^s_u\vdash O^u,\overline{O}_s
                        }(\equiv)
                }({\rm{Cut}})
            }
            {
            \cfrac
                {
                \cfrac
                    {
                 \delta^{u j}_{iv}
                    \cdot
                    [{\rm{leaves}}]^s_u
                    \vdash
                   O^v,\overline{O}_s,O^i\otimes\overline{O}_j
                    }
                    {
                    \delta_v^j\cdot[{\rm{leaves}}]_i^s
                    \vdash
                       O^v,\overline{O}_s,O^i\otimes\overline{O}_j
                    }(\equiv)
                }
                {
                \delta_v^j\cdot[{\rm{leaves}}]_i^s
                    \vdash
                       O^v\wp\overline{O}_s,O^i\otimes\overline{O}_j
                    }(\wp).
            }({\rm{Cut}})
    $$
We read the  conclusion as
$$\delta_v^j\cdot[{\rm{leaves}}]_i^s
                    \vdash
                       \phi( S)^v_s,\phi(\overline{NP})^i_j.$$
According to the general prescription, the tensor axiom for ``leaves'' takes the form
$$\delta_v^j\cdot[{\rm{leaves}}]_i^s
                    \vdash
                       S^v_s,\overline{NP}^i_j.$$

It might be entertaining to work out the tensor translation of an ACG term corresponding to the relative pronoun ``whom''. Typically, the constant in the abstract signature axiom
$$\vdash {\rm{whom}}:(NP\multimap S)\multimap NP\multimap NP$$
will be translated by the lexicon to the term
$$\phi({\rm{whom}})=\lambda fxt.(x(\mbox{who}( (f\cdot(\lambda y.y))t))).$$
The tensor translation will be   (\ref{whom}) in Section \ref{elaborate example}, and, arguably, this is more readable.

 \section{Extended tensor grammars}\label{extended things section}
We proceed  now to extending tensor grammars with {\it binding operators} in order to accommodate noncommutative constructions. The two new operators will be denoted as $\nabla$ and $\triangle$.

Loosely speaking, {\it extended tensor types} are defined by the same rules as tensor types plus a rule for binding operators:
\begin{itemize}
  \item If $A$ is a type  with $i\in FSub(A)$, $j\in FSup(A)$, then $\nabla^i_j A $, $\triangle^i_j A $ are types.
\end{itemize}
However,  as is usual with binding operators, we need to be accurate and treat all ``bureaucratic'' details carefully.

\subsection{Extended tensor types}
We first define {\it extended tensor type expressions}.

Extended tensor type expressions have upper and lower indices that may have multiple occurrences. An index occurrence may be {\it free} or {\it bound}, and any index may have at most one free occurrence. A bound index occurrence is bound by some binding operator. Also, a  binding operator occurring in an extended type expression $A$ has a {\it scope}, which is a set of free index occurrences in $A$.

For an extended tensor type expression we will denote the set of indices having a free occurrence as an upper, respectively lower, index in $A$ as $FSup(A)$, respectively $FSub(A)$. As usual, we write $F\mathit{Ind}(A)=FSub(A)\cup FSup(A)$.

    Given a set $P$ of  positive atomic type symbols together with the valency function $v:P\to{\bf N}^2$, we simultaneously define {\it extended tensor type expressions}, {\it free} and {\it bound index occurrences} and {\it binding operator scopes} by the following rules.
\begin{itemize}
  \item If $p\in Lit=P\cup\overline{P}$ with $v(p)=(m,n)$ and $i_1,\ldots, i_m,j_1,\ldots,j_n$ are pairwise distinct elements of $\mathit{Ind}$ then $A=p_{j_1\ldots j_n}^{i_1\ldots i_m}$ is an extended tensor  type expression, and every index occurrence is free.
          \item If $A$, $B$ are extended tensor  type expressions  with $F\mathit{Ind}(A)\cap F\mathit{Ind}(B)=\emptyset$
  then $C=A\otimes B$, respectively $C=A\wp B$, is an extended tensor  type expression.

  A free index occurrence in $A$ or $B$ remains  free in $C$.
  An occurrence bound by a binding operator $Q$ in $A$ or $B$ remains bound by $Q$ in $C$.
  A free occurrence in the scope of a binding operator $Q$ in $A$ or $B$ remains in the scope of $Q$ in $C$.
  \item If $A$ is an extended tensor type expression with $i\in FSub(A)$, $j\in FSup(A)$, then $A'=\nabla^i_j A $, respectively $A'=\triangle^i_j A $, is an extended tensor type expression.

       A free occurrence of an  index  $k\not=i,j$  in $A$ remains  free in  $A'$.
      An occurrence  bound by a binding operator $Q$ in $A$ remains  bound by $Q$ in $A'$.  A free occurrence of  an index $k\not=i,j$  in the scope of a binding operator $Q$ in $A$ remains in the scope of $Q$ in $A'$.

       The unique  free occurrences of $i$ and $j$ in $A$ are bound in $A'$ by the binding operator $\nabla_j^i$, respectively $\triangle^i_j$, in the beginning of $A'$. Any free  index occurrence in $A'$ is in the scope of the binding operator $\nabla_j^i$, respectively $\triangle^i_j$, in the beginning of $A'$.
\end{itemize}

We extend the notion of {\it $\alpha$-equivalence} from term expressions and syntactic typing judgements to extended tensor type expressions, denoting it as $\sim_\alpha$. It is the smallest equivalence relation  satisfying the following.
 \begin{itemize}
   \item Let $A$ be an extended tensor type expression,  $Q=\nabla^i_j$, respectively $Q=\triangle^i_j$, be an occurrence of a binding operator in $A$,  $(i',j')$ be a pair of indices not occurring freely
       in $A$  in the scope of $Q$.
        Let $A'$ be obtained from $A$ by simultaneously replacing $Q$ with $\nabla_{j'}^{i'}$, respectively $\triangle_{j'}^{i'}$ and replacing the unique occurrences of $i$ and $j$  bound by $Q$ in $A$ with $i'$ and $j'$ respectively. Then $A\sim_\alpha A'$.
 \end{itemize}

An {\it extended tensor type} is an equivalence class for the $\alpha$-equivalence of extended tensor type expressions.

Note that the sets of free indices are invariant under $\alpha$-equivalence of extended tensor type expressions, so they are well-defined for extended tensor types as well.

Finally, in an analogy with tensor type symbols, we want to define {\it extended tensor type symbols} as extended tensor types with all  free indices erased... but locations of bound indices kept intact.

In order to do this accurately we introduce the last equivalence relation in this section, defined on extended tensor types by the following rules.
\begin{itemize}
  \item If $A\sim_\alpha A'$ then $A\sim A'$.
  \item Let $i\in F\mathit{Ind}(A)$ and $i'\not\in F\mathit{Ind}(A)$ be such that if the unique free occurrence of $i$ in $A$ is in the scope of a binding operator $Q$, then $i'$ is not bound by $Q$ in $A$. Let $A'$ be obtained from $A$ by replacing the unique free occurrence of $i$ in $A$ with $i'$. Then $A\sim A'$.
\end{itemize}
We define {\it extended tensor type symbols} as
equivalence classes for the above relation.

 We denote the set of extended tensor type symbols over $P$ as $Tp_{\otimes,\nabla}(P)$.

  Just as in the case of tensor types, an extended tensor type can be reconstructed from its symbol by specifying the ordered sets of free upper and lower indices. Accordingly, we continue to use the notation $A_J^I$ to denote the extended tensor type defined by the extended tensor type  symbol $A$ and ordered index sets $I$ and $J$.

 {\it Type valency} for extended tensor type symbols is defined by formula (\ref{valency}) supplemented with
 $$
 v(\nabla^i_jA)= v(\triangle^i_jA)=v(A)-(1,1).
$$

Finally, the {\it dual type } $\overline A$ of an extended tensor type  $A$ is defined by (\ref{dual type}) supplemented with
\be\label{extended type dual}
\overline{\nabla_i^jA}= \triangle^i_j\overline{A},\quad
\overline{\triangle_i^jA}= \nabla^i_j\overline{A}.
\ee

\subsection{Extended typing rules}
{\it Extended tensor sequents} and
{\it extended tensor typing judgements} are defined exactly as ordinary tensor sequents and tensor typing judgements, by replacing every occurrence of the adjective ``tensor'' with ``extended tensor''.

{\it Extended tensor type calculus} ({\bf ETTC}) is defined by the rules of {\bf TTC} supplemented with the following
$$\cfrac{
    \delta_\beta^\alpha\cdot t\vdash\Gamma,A
    }
    {
    t\vdash\Gamma,\nabla^\alpha_\beta A
    }~(\nabla),
\quad
\cfrac{
    t\vdash\Gamma,A
    }
    {
    \delta^\beta_\alpha\cdot t\vdash\Gamma,\triangle_\beta^\alpha A
    }~(\triangle),
$$
where  it is assumed that $\alpha\in FSub(A),\beta\in FSup(A)$.

Note that the term $t$ in the premise of the $(\triangle)$ rule above must have { free} occurrences of $\alpha$ and $\beta$ (as an upper and a lower index respectively).
 In the conclusion,  the term $\delta^\beta_\alpha\cdot t$ has these occurrences {\it bound}. Thus free indices to the left and to the right of the turnstile do match, and the typing judgement  is well defined.

 {\it Lambek restriction} on the $(\nabla)$ rule consists in requiring that the context $\Gamma$ is not empty.
   {\it Extended tensor type calculus with Lambek restriction} (${\bf ETTC}_-$) is defined by the rules of {\bf ETTC} with Lambek restriction on  the $(\nabla)$ rule.

\subsubsection{Meaning of extended types}
Extended tensor typing judgements can be represented geometrically in the same way as ordinary tensor typing judgements: terms to the left of the turnstile corresponds to graphs, and {\it free} indices correspond to vertices. The sequent to the right of the turnstile still induces an ordering of free indices. As for bound indices, they are not in the picture.

\begin{figure}
\begin{subfigure}{1\textwidth}
        \centering
         \begin{tikzpicture}[xscale=1.5,yscale=1.5]
         \draw[draw=black](1.2,-.2)rectangle(-.5,-.5);\node at(.35,-.35){$t$};

         \begin{scope}[shift={(.1,0)}]
            \begin{scope}[shift={(-.3,.2)}]
                \draw [fill] (-.1,0) circle [radius=0.01];
                \draw [fill] (-0.05,0) circle [radius=0.01];
                \draw [fill] (0,0) circle [radius=0.01];
                \draw [fill] (0.05,0) circle [radius=0.01];
                        \draw [fill] (0.1,0) circle [radius=0.01];
                \draw[thick,dashed, -](-.2,0)--(-.2,.25) -- (.2,.25)-- (.2,0);
                        \node at (0,.25) [above] {$\Gamma$};
            \end{scope}

                \draw[thick, -](-.5,.2)--(-.5,-.2);
                \draw[thick, -](-.1,.2)--(-.1,-.2);

         \end{scope}

                \begin{scope}[shift={(.5,.2)}]

                \draw[thick,-](-.25,0)--(-.25,-.4);
                \draw[thick,-](.55,0)--(.55,-.4);

                \draw [fill] (-.1,0) circle [radius=0.01];
                \draw [fill] (-.15,0) circle [radius=0.01];
                \draw [fill] (0.15,0) circle [radius=0.01];
                \draw [fill] (0.2,0) circle [radius=0.01];
                \draw [fill] (0.1,0) circle [radius=0.01];
                \draw [fill] (.4,0) circle [radius=0.01];
                \draw [fill] (.45,0) circle [radius=0.01];

                \draw[thick,dashed, -](-.25,0)--(-.25,.25) -- (.55,.25)-- (.55,0);
                        \node at (.15,.25) [above] {$A$};

                \begin{scope}[shift={(0,-.05)}]
                    \draw [fill] (0.3,0) circle [radius=0.01];
                    \draw [fill] (0,0) circle [radius=0.05];
                    \node at (0,0) [above] {$\beta$};
                    \node at (0.3,.0) [above] {$\alpha$};
                    \draw[thick,->](0,0) to  [out=-90,in=180] (.15,-.25) to  [out=0,in=-90](0.3,0);
                \end{scope}

                \end{scope}

\node at(1.6,-.35){$\Rightarrow$};

            \begin{scope}[shift={(2.5,0)}]
         \draw[draw=black](1.2,-.2)rectangle(-.5,-.5);\node at(.35,-.35){$t$};

         \begin{scope}[shift={(.1,0)}]
            \begin{scope}[shift={(-.3,.2)}]
                \draw [fill] (-.1,0) circle [radius=0.01];
                \draw [fill] (-0.05,0) circle [radius=0.01];
                \draw [fill] (0,0) circle [radius=0.01];
                \draw [fill] (0.05,0) circle [radius=0.01];
                        \draw [fill] (0.1,0) circle [radius=0.01];
                \draw[thick,dashed, -](-.2,0)--(-.2,.1) -- (.2,.1)-- (.2,0);
                        \node at (0,.1) [above] {$\Gamma$};
            \end{scope}

                \draw[thick, -](-.5,.2)--(-.5,-.2);
                \draw[thick, -](-.1,.2)--(-.1,-.2);

         \end{scope}

                \begin{scope}[shift={(.5,.2)}]

                \draw[thick,-](-.25,0)--(-.25,-.4);
                \draw[thick,-](.55,0)--(.55,-.4);

                \draw [fill] (-.1,0) circle [radius=0.01];
                \draw [fill] (-.15,0) circle [radius=0.01];
                \draw [fill] (0.15,0) circle [radius=0.01];
                \draw [fill] (0.2,0) circle [radius=0.01];
                \draw [fill] (0.1,0) circle [radius=0.01];
                \draw [fill] (.4,0) circle [radius=0.01];
                \draw [fill] (.45,0) circle [radius=0.01];

                \draw[thick,dashed, -](-.25,0)--(-.25,.1) -- (.55,.1)-- (.55,0);
                        \node at (.15,.1) [above] {$\nabla_\beta^\alpha A$};

                \begin{scope}[shift={(0,-.05)}]
                    %\draw [fill] (0.3,0) circle [radius=0.01];
%                    \draw [fill] (0,0) circle [radius=0.05];
                    %\node at (0,0) [above] {$\beta$};
%                    \node at (0.3,.01) [above] {$\alpha$};
%                    \draw[thick,->](0,0) to  [out=-90,in=180] (.15,-.25) to  [out=0,in=-90](0.3,0);
                \end{scope}

                \end{scope}

                \end{scope}

        \end{tikzpicture}
        \caption{$(\nabla)$ rule}
        \label{nabla_geometrically}
        \end{subfigure}
\begin{subfigure}{1\textwidth}
        \centering
         \begin{tikzpicture}[xscale=1.5,yscale=1.5]
         \draw[draw=black](1.2,-.2)rectangle(-.5,-.5);\node at(.35,-.35){$t$};

         \begin{scope}[shift={(.1,0)}]
            \begin{scope}[shift={(-.3,.2)}]
                \draw [fill] (-.1,0) circle [radius=0.01];
                \draw [fill] (-0.05,0) circle [radius=0.01];
                \draw [fill] (0,0) circle [radius=0.01];
                \draw [fill] (0.05,0) circle [radius=0.01];
                        \draw [fill] (0.1,0) circle [radius=0.01];
                \draw[thick,dashed, -](-.2,0)--(-.2,.25) -- (.2,.25)-- (.2,0);
                        \node at (0,.25) [above] {$\Gamma$};
            \end{scope}

                \draw[thick, -](-.5,.2)--(-.5,-.2);
                \draw[thick, -](-.1,.2)--(-.1,-.2);

         \end{scope}

                \begin{scope}[shift={(.5,.2)}]

                \draw[thick,-](-.25,0)--(-.25,-.4);
                \draw[thick,-](.55,0)--(.55,-.4);

                \draw[thick,->](0,-.05)--(0,-.4);
                \draw[thick,<-](.3,-.05)--(.3,-.4);

                \draw [fill] (-.1,0) circle [radius=0.01];
                \draw [fill] (-.15,0) circle [radius=0.01];
                \draw [fill] (0.15,0) circle [radius=0.01];
                \draw [fill] (0.2,0) circle [radius=0.01];
                \draw [fill] (0.1,0) circle [radius=0.01];
                \draw [fill] (.4,0) circle [radius=0.01];
                \draw [fill] (.45,0) circle [radius=0.01];

                \draw[thick,dashed, -](-.25,0)--(-.25,.25) -- (.55,.25)-- (.55,0);
                        \node at (.15,.25) [above] {$A$};

                \begin{scope}[shift={(0,-.05)}]
                    \draw [fill] (0.3,0) circle [radius=0.01];
                    \draw [fill] (0,0) circle [radius=0.05];
                    \node at (0,0) [above] {$\beta$};
                    \node at (0.3,.0) [above] {$\alpha$};
                    %\draw[thick,->](0,0) to  [out=-90,in=180] (.15,-.25) to  [out=0,in=-90](0.3,0);
                \end{scope}

         \end{scope}

\node at(1.6,-.35){$\Rightarrow$};

            \begin{scope}[shift={(2.5,0)}]
         \draw[draw=black](1.2,-.2)rectangle(-.5,-.5);\node at(.35,-.35){$t$};

         \begin{scope}[shift={(.1,0)}]
            \begin{scope}[shift={(-.3,.2)}]
                \draw [fill] (-.1,0) circle [radius=0.01];
                \draw [fill] (-0.05,0) circle [radius=0.01];
                \draw [fill] (0,0) circle [radius=0.01];
                \draw [fill] (0.05,0) circle [radius=0.01];
                        \draw [fill] (0.1,0) circle [radius=0.01];
                \draw[thick,dashed, -](-.2,0)--(-.2,.25) -- (.2,.25)-- (.2,0);
                        \node at (0,.25) [above] {$\Gamma$};
            \end{scope}

                \draw[thick, -](-.5,.2)--(-.5,-.2);
                \draw[thick, -](-.1,.2)--(-.1,-.2);

         \end{scope}

                \begin{scope}[shift={(.5,.2)}]

                \draw[thick,-](-.25,0)--(-.25,-.4);
                \draw[thick,-](.55,0)--(.55,-.4);

                \draw[thick,-](0,-.05)--(0,-.4);
                \draw[thick,-](.3,-.05)--(.3,-.4);

                \draw [fill] (-.1,0) circle [radius=0.01];
                \draw [fill] (-.15,0) circle [radius=0.01];
                \draw [fill] (0.15,0) circle [radius=0.01];
                \draw [fill] (0.2,0) circle [radius=0.01];
                \draw [fill] (0.1,0) circle [radius=0.01];
                \draw [fill] (.4,0) circle [radius=0.01];
                \draw [fill] (.45,0) circle [radius=0.01];

                \draw[thick,dashed, -](-.25,0)--(-.25,.25) -- (.55,.25)-- (.55,0);
                        \node at (.15,.25) [above] {$\triangle_\beta^\alpha A$};

                \begin{scope}[shift={(0,-.05)}]
                    %\draw [fill] (0.3,0) circle [radius=0.01];
                    %\draw [fill] (0,0) circle [radius=0.05];
                    %\node at (0,0) [above] {$\beta$};
                    %\node at (0.3,.01) [above] {$\alpha$};
                    \draw[thick,-](0,0) to  [out=90,in=180] (.15,.2) to  [out=0,in=90](0.3,0);
                \end{scope}

         \end{scope}
                \end{scope}

        \end{tikzpicture}
        \caption{$(\triangle)$ rule}
        \label{delta_geometrically}
        \end{subfigure}
        \caption{Extended typing rules geometrically}
        \label{Extended typing rules geometrically}
   \end{figure}
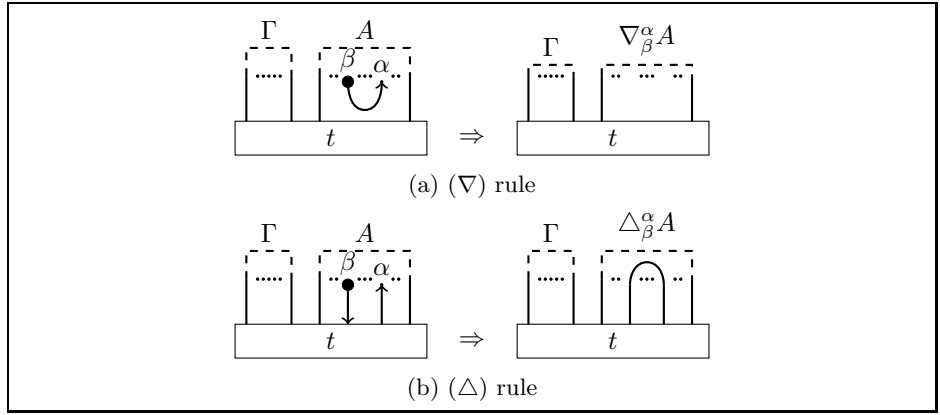

With these conventions, the extended typing rules are shown (rather schematically) in Figure \ref{Extended typing rules geometrically}.
      The $(\triangle)$ rule simply glues together the two bound indices/vertices. The $(\nabla)$ rule is applicable only in the case when the corresponding indices/vertices are connected with an edge carrying no label. Then this edge (together with its endpoints) is erased from the picture completely. The information about the erased edge is stored in the introduced type (or link).

      Thus, terms of type $\nabla^\alpha_\beta A$ encode the {\it subtype} of $A$ consisting of terms/graphs with that specific form: vertices corresponding to $\alpha$ and $\beta$ are connected with an edge, and the connecting edge carries no label.
      It can be observed that
      we have the admissible rule
      $$
      \cfrac{t\vdash \Gamma,\nabla^\alpha_\beta A}{\delta^\alpha_\beta \cdot t\vdash \Gamma,A}~(\nabla^{-1}),
      $$ so that ``decoding'' from $\nabla^\alpha_\beta A$ to $A$ is always possible.
      \nb
      For any extended tensor type symbol $A$ the following rule is admissible in {\bf ETTC}.
       $$
      \cfrac{t\vdash \Gamma,\nabla^\alpha_\beta (A^{I\beta K}_{J\alpha L})}{\delta^\alpha_\beta \cdot t\vdash \Gamma,(A^{I\beta K}_{J\alpha L})}~(\nabla^{-1}).
      $$
      \nbe
      {\bf Proof}
            We have the derivable typing judgement
            $$
      \delta^{\alpha'  \beta }_{\alpha  \beta' }\cdot
      \delta^{IKJ'L'}_{I'K'JL}
      \vdash \overline{A}_{\overline{K}{\beta }\overline{I}}^{\overline{L}{\alpha  }\overline{J}}, A^{I'\beta' K'}_{J'\alpha'  L'},
      $$
      and after an application of the $(\triangle)$ rule this yields
      $$
      \delta^\alpha _\beta \cdot\delta^{\alpha' \beta }_{\alpha \beta' }\cdot
      \delta^{IKJ'L'}_{I'K'JL}
      \vdash \triangle^\beta _\alpha (\overline{A}_{\overline{K}\beta \overline{I}}^{\overline{L}{\alpha }\overline{J}}), A^{I'\beta' K'}_{J'\alpha' L'}.
      $$
                  Cutting this with
                  $t\vdash \Gamma,\nabla^\alpha _\beta (A^{I\beta K}_{J\alpha L})$
                  (and slightly permuting factors in the term expression to the left of the turnstile)
                  we obtain the typing  judgement
                  $$\delta^{\alpha' \beta }_{\alpha \beta' }\cdot
     \delta^{IKJ'L'}_{I'K'JL}\cdot\delta^{\alpha}_{\beta}\cdot t\vdash A^{I'\beta' K'}_{J'\alpha' L'},$$
      which is precisely
      $\delta^\alpha _\beta \cdot t\vdash \Gamma,A^{I\beta K}_{J\alpha L}$ by Proposition \ref{alpha conversion} on change of coordinates. $\Box$

\subsubsection{Cut elimination}
\bl\label{Cut-elimination extended}
Any  typing judgement derivable in {\bf ETTC} ($ {\bf ETTC}_-$) is   derivable without the Cut rule.
\el
{\bf Proof} Similar to  Lemma \ref{Cut-elimination simple}.

We skip details with commutation or rules and work out the new cut-elimination step corresponding to the new operators.

Let a derivation end with the Cut rule of the form
$$
\cfrac{
    \cfrac{
        \delta^\alpha_\beta\cdot t\vdash\Gamma, A
        }
        {
        t\vdash\Gamma,\nabla^\alpha_\beta A
        }~{ }(\nabla)
    \quad
    \cfrac{
        s\vdash\Theta,\overline{A}
        }
        {
        \delta_{\beta}^{\alpha}\cdot s\vdash\Theta,\triangle^{\beta}_{\alpha}(\overline{A})
        }~(\triangle)
    }
    {
    \delta_{\beta}^{\alpha}\cdot t s\vdash\Gamma,\Theta
    } ~(\rm{Cut}).
$$
This  reduces to the derivation with a cut on $A$:
$$
\cfrac{
 \delta_{\beta}^{\alpha}\cdot t\vdash\Gamma, A
\quad
s\vdash\Theta,\overline{A}
   }
    {
    \delta_{\beta}^{\alpha}\cdot t\cdot s\vdash\Gamma,\Theta
    } ~(\rm{Cut}),
$$
which results just in the same typing judgement. $\Box$

\subsubsection{Variations on deduction theorem}
Lemma \ref{deduction theorem}, ``Deduction theorem'', strictly speaking, is not valid for {\bf ETTC}. (In the same way as the usual Deduction theorem is not valid for first order logic.)
Indeed, if we have a tensor signature with the axiom
$$\delta_{ik}^{jl}\vdash p^{ik}_{jl},$$ then we  can immediately derive the typing judgement
$$\delta_{i}^{j}\vdash \nabla_k^l p^{ik}_{jl},$$ but for no tensor term $t$ whatsoever the judgement
$$t\vdash\overline{p}^{l'j'}_{k'i'},\nabla_k^l p^{ik}_{jl}$$ is derivable in {\bf ETTC}.

However, there is a  limited form of Lemma \ref{deduction theorem}, which will be sufficient for our purposes.

Let us say that a tensor term $t$ is {\it lexical} if it cannot be written in the form $t=t'\cdot\delta^i_j$, where both indices $i,j$ are free. (In the geometric language, a term is lexical if in the corresponding graph every edge is labeled with a nonempty word.) An extended tensor typing judgement $t\vdash A$ is {\it lexicalized} if the term $t$ is lexical.
\bl\label{extended deduction theorem}
Let $\Xi$ be  a finite multiset of lexicalized extended tensor typing judgements
$$\Xi=\{t_{(1)}\vdash(\Gamma_{(1)})^{I_1}_{J_1},~\ldots,~t_{(n)}\vdash(\Gamma_{(n)})^{I_n}_{J_n}\}.$$
Choose syntactic representatives for elements of $\Xi$ so that all index sequences $I_1,\ldots,I_n$ $J_1,\ldots, J_n$ above are pairwise disjoint.

An extended tensor  typing judgement
 is derivable
 from  $\Xi$ iff it has a  syntactic representation of the form
$$
t_{(0)}\cdot t_{(1)}\cdots t_{(n)}\vdash\Gamma
$$
where $t_{(0)}$ is a closed term, and
the typing judgement
\be\label{4'}
t_{(0)}\vdash
\overline{(F_{(1)})}_{\overline{I_1}}^{\overline{J_1}}, \ldots, \overline{(F_{(n)})}_{\overline{I_n}}^{\overline{J_n}},\Gamma
\ee
is derivable in {\bf ETTC}.
\el
{\bf Proof} Similar to Lemma \ref{deduction theorem}.

An important point is the new step in the induction on derivation from $\Xi$ when the last rule is the $(\nabla)$ rule. Then we have a typing judgement $t\vdash\Gamma,\nabla^\alpha_\beta A$ obtained from $t'\vdash\Gamma,A$, where the latter  is derivable from $\Xi$ and $t'=t\cdot\delta^\alpha_\beta $.

Using  the induction hypothesis we may assume without loss of generality that
$t'=t_{(0)}'\cdot t_{(1)}\cdots t_{(n)}$
where $t_{(0)}'$ is  closed and the typing judgement
\be\label{5'}
t_{(0)}'\vdash
\overline{(F_{(1)})}_{\overline{I_1}}^{\overline{J_1}}, \ldots, \overline{(F_{(n)})}_{\overline{I_n}}^{\overline{J_n}},\Gamma, A
\ee
is derivable.

Now, the indices $\alpha,\beta$ are free in $t'$, and the terms $t_{(1)},\ldots, t_{(n)}$ are lexical, so $\delta^\alpha_\beta $ cannot be a factor in any of them. It follows that we have a representation $t_{(0)}'=t_{(0)}\cdot \delta^\alpha_\beta $, for some closed term $t_{(0)}$ and $t=t_{(0)}\cdot t_{(1)}\cdots t_{(n)}$ with (\ref{4'}) derivable from (\ref{5'}) by the $(\nabla)$ rule.
 $\Box$

\subsubsection{Lambek style types and examples of admissible rules}
Let us say that an extended tensor type symbol $A$ is {\it Lambek-style} if its valency is $(1,1)$.
An extended tensor type is Lambek style if its symbol is Lambek style.

For Lambek style type symbols $A$, $B$ we define new Lambek style types
  \be\label{translating lambek types}
 (A\bullet B)^i_j=\triangle_\beta^\alpha(A^i_\alpha\otimes
 B^\beta_j),\quad
(B/A)^i_j=\nabla^\alpha_\beta(B^i_\alpha\wp \overline{A}_j^\beta ),\quad
(A\backslash B)^i_j=\nabla^\alpha_\beta(\overline{A}_\alpha^i\wp B_j^\beta ).
\ee
This corresponds to operations of Lambek calculus, which will be discussed shortly. (It is interesting to compare (\ref{translating lambek types}) with formulas for embedding {\bf LC} into {\bf MILL1} in \cite{MootPiazza}, \cite{Moot_comparing}, which look tantalizingly similar.)
As a preparation, let us derive some admissible rules for Lambek style types that correspond to slash elimination rules of {\bf LC} in the natural deduction formulation.
\nb
For Lambek style types $A,B$ the following rules are admissible in {\bf ETTC}:
$$
\cfrac
    {
    t\vdash \Theta, (B/A)^i_j \quad s\vdash A_\beta^j,\Gamma,
    }
    {
     ts\cdot\delta^\alpha_\beta\vdash\Theta,  B_\alpha^i,\Gamma
    }(/{\rm{E}}),
\quad
\cfrac
    {
    s\vdash\Gamma, A^\alpha_i\quad t\vdash (A\backslash B)^i_j,\Theta
    }
    {
    \delta^\alpha_\beta\cdot st\vdash\Gamma,  B^\beta_j,\Theta
    }(\backslash{\rm{E}}).
$$
\nbe
{\bf Proof}
The rule $(\backslash{\rm{E}})$ is obtained as
$$
\cfrac
    {
    \begin{array}{c}
       \\
       \\
       s\vdash\Gamma, A^\alpha_i
    \end{array}
       \quad
    \cfrac
        {
        \cfrac
            {
            t\vdash \nabla^\alpha_\beta(\overline{A}^i_\alpha\wp B^\beta_j),\Theta
            }
            {
            \delta^\alpha_\beta t\vdash \overline{A}^i_\alpha\wp B^\beta_j,\Theta
            }(\nabla^{-1})
        }
        {
        \delta^\alpha_\beta t\vdash \overline{A}^i_\alpha, B^\beta_j,\Theta
        }(\wp^{-1})
    }
    {
    \delta^\alpha_\beta\cdot st\vdash\Gamma,  B^\beta_j,\Theta
    }({\rm{Cut}}).
$$
The other one is similar. $\Box$
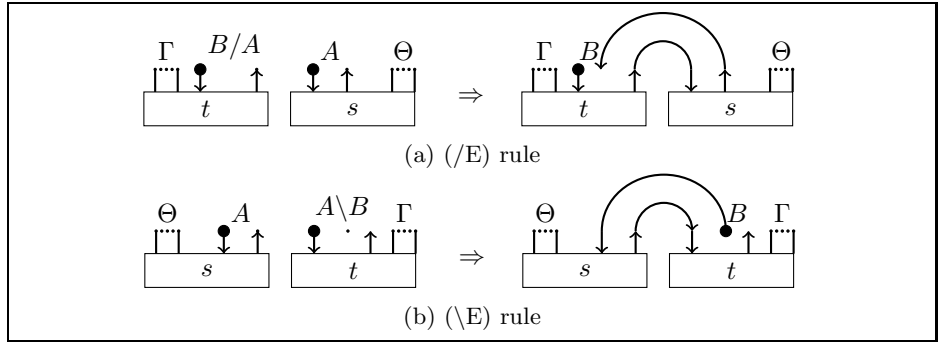
\begin{figure}
\begin{subfigure}{1\textwidth}
        \centering
         \begin{tikzpicture}[xscale=1.5,yscale=1.5]
         \begin{scope}[shift={(-3.2,0)}]
         \begin{scope}[shift={(.5,0)}]
             \draw[draw=black](.9,-.2)rectangle(-.2,-.5);\node at(.35,-.35){$t$};

              %\draw[thick,dotted,-](0.3,0) to  [out=-90,in=180] (1.3,-.75) to  [out=0,in=-90](2.3,0);
              %\draw[thick,dotted,-](0,0) to  [out=-90,in=180] (1.3,-1) to  [out=0,in=-90](2.6,0);

                    \draw [fill] (-.1,0) circle [radius=0.01];
                    \draw [fill] (-0.05,0) circle [radius=0.01];
                    \draw [fill] (0,0) circle [radius=0.01];
                    \draw [fill] (0.05,0) circle [radius=0.01];
                            \draw [fill] (0.1,0) circle [radius=0.01];
                    \draw[thick,-](-.1,0)--(-.1,-.2);
                    \draw[thick,-](.1,0)--(.1,-.2);
                    \node at (0,0) [above] {$\Gamma$};

                    \begin{scope}[shift={(.5,0)}]
                    \draw [fill] (-.2,0) circle [radius=0.05];
                    %\draw [fill] (0.15,0) circle [radius=0.01];
    %                \draw [fill] (0.2,0) circle [radius=0.01];
    %                \draw [fill] (0.1,0) circle [radius=0.01];
                    \draw [fill] (0.3,0) circle [radius=0.01];
                   % \node at(-.3,0)[above]{$j$};
%                    \node at(.3,0)[above]{$\alpha$};
                    \draw[thick,->](-.2,0)--(-.2,-.2);
                    \draw[thick,<-](0.3,0)--(0.3,-.2);
                    \node at (0.1,0) [above] {$B/A$};
    %                \node at (0.3,0) [right] {$I$};
                    \end{scope}

                \begin{scope}[shift={(1.3,0)}]
                \draw[draw=black](.9,-.2)rectangle(-.2,-.5);\node at(.35,-.35){$s$};
                 \draw [fill] (0,0) circle [radius=0.05];
                    %\draw [fill] (0.15,0) circle [radius=0.01];
    %                \draw [fill] (0.2,0) circle [radius=0.01];
    %                \draw [fill] (0.1,0) circle [radius=0.01];
                    \draw [fill] (0.3,0) circle [radius=0.01];
                    \draw[thick,->](0,0)--(0,-.2);
                    \draw[thick,<-](0.3,0)--(0.3,-.2);
                    \node at (0.15,0) [above] {${A}$};
                    %\node at (0,0) [left] {$I'$};
    %                \node at (0.3,0) [right] {$J'$};
                \end{scope}
      \end{scope}

                %\draw [fill] (0,0) circle [radius=0.05];
    %                \draw [fill] (0.15,0) circle [radius=0.01];
    %                \draw [fill] (0.2,0) circle [radius=0.01];
    %                \draw [fill] (0.1,0) circle [radius=0.01];
    %                        \draw [fill] (0.3,0) circle [radius=0.05];
    %                \draw[dashed,-](0,-.05)--(0,-.15)--(0.3,-.15)--(0.3,-.05);
    %                \node at (0.15,-.15) [below] {$J$};

                    \begin{scope}[shift={(2.6,0)}]
                    \draw [fill] (-.1,0) circle [radius=0.01];
                    \draw [fill] (-0.05,0) circle [radius=0.01];
                    \draw [fill] (0,0) circle [radius=0.01];
                    \draw [fill] (0.05,0) circle [radius=0.01];
                            \draw [fill] (0.1,0) circle [radius=0.01];
                    \draw[thick,-](-.1,0)--(-.1,-.2);
                    \draw[thick,-](.1,0)--(.1,-.2);
                    \node at (0,0) [above] {$\Theta$};

                    \end{scope}
        \end{scope}

        \node at (0,-.25){$\Rightarrow$};

        \begin{scope}[shift={(.65,0)}]
        \draw[draw=black](.9,-.2)rectangle(-.2,-.5);\node at(.35,-.35){$t$};

          \draw[thick,-](0.8,0) to  [out=90,in=180] (1.05,.25) to  [out=0,in=90](1.3,0);
          \draw[thick,<-](.5,0) to  [out=90,in=180] (1.05,.5) to  [out=0,in=90](1.6,0);
          %\draw[thick,dotted,-](0.3,0) to  [out=-90,in=180] (1.3,-.75) to  [out=0,in=-90](2.3,0);
          %\draw[thick,dotted,-](0,0) to  [out=-90,in=180] (1.3,-1) to  [out=0,in=-90](2.6,0);

                \draw [fill] (-.1,0) circle [radius=0.01];
                \draw [fill] (-0.05,0) circle [radius=0.01];
                \draw [fill] (0,0) circle [radius=0.01];
                \draw [fill] (0.05,0) circle [radius=0.01];
                        \draw [fill] (0.1,0) circle [radius=0.01];
                \draw[thick,-](-.1,0)--(-.1,-.2);
                \draw[thick,-](.1,0)--(.1,-.2);
                \node at (0,0) [above] {$\Gamma$};

                \begin{scope}[shift={(.5,0)}]
                \draw [fill] (-.2,0) circle [radius=0.05];
                %\draw [fill] (0.15,0) circle [radius=0.01];
%                \draw [fill] (0.2,0) circle [radius=0.01];
%                \draw [fill] (0.1,0) circle [radius=0.01];
                %\draw [fill] (0.3,0) circle [radius=0.01];
                \node at(-.1,0)[above]{$B$};
                \draw[thick,->](-.2,0)--(-.2,-.2);
                \draw[thick,<-](0.3,0)--(0.3,-.2);
             %   \node at (0,0) [left] {$J$};
%                \node at (0.3,0) [right] {$I$};
                \end{scope}

            \begin{scope}[shift={(1.3,0)}]
            \draw[draw=black](.9,-.2)rectangle(-.2,-.5);\node at(.35,-.35){$s$};
            % \draw [fill] (0,0) circle [radius=0.05];
                %\draw [fill] (0.15,0) circle [radius=0.01];
%                \draw [fill] (0.2,0) circle [radius=0.01];
%                \draw [fill] (0.1,0) circle [radius=0.01];
               % \draw [fill] (0.3,0) circle [radius=0.01];
                \draw[thick,->](0,0)--(0,-.2);
                \draw[thick,<-](0.3,0)--(0.3,-.2);
                %\node at (0,0) [left] {$I'$};
%                \node at (0.3,0) [right] {$J'$};
            \end{scope}
            %\draw [fill] (0,0) circle [radius=0.05];
%                \draw [fill] (0.15,0) circle [radius=0.01];
%                \draw [fill] (0.2,0) circle [radius=0.01];
%                \draw [fill] (0.1,0) circle [radius=0.01];
%                        \draw [fill] (0.3,0) circle [radius=0.05];
%                \draw[dashed,-](0,-.05)--(0,-.15)--(0.3,-.15)--(0.3,-.05);
%                \node at (0.15,-.15) [below] {$J$};

                \begin{scope}[shift={(2.1,0)}]
                \draw [fill] (-.1,0) circle [radius=0.01];
                \draw [fill] (-0.05,0) circle [radius=0.01];
                \draw [fill] (0,0) circle [radius=0.01];
                \draw [fill] (0.05,0) circle [radius=0.01];
                        \draw [fill] (0.1,0) circle [radius=0.01];
                \draw[thick,-](-.1,0)--(-.1,-.2);
                \draw[thick,-](.1,0)--(.1,-.2);
                \node at (0,0) [above] {$\Theta$};

                \end{scope}
        \end{scope}
        \end{tikzpicture}
        \caption{$(/{\rm{E}})$ rule}
        %\label{Cut_geometrically}
        \end{subfigure}
    \begin{subfigure}{1\textwidth}
        \centering
         \begin{tikzpicture}[xscale=1.5,yscale=1.5]
         \begin{scope}[shift={(-3.2,0)}]
         \begin{scope}[shift={(.5,0)}]
             \draw[draw=black](.9,-.2)rectangle(-.2,-.5);\node at(.35,-.35){$s$};

              %\draw[thick,dotted,-](0.3,0) to  [out=-90,in=180] (1.3,-.75) to  [out=0,in=-90](2.3,0);
              %\draw[thick,dotted,-](0,0) to  [out=-90,in=180] (1.3,-1) to  [out=0,in=-90](2.6,0);

                    \draw [fill] (-.1,0) circle [radius=0.01];
                    \draw [fill] (-0.05,0) circle [radius=0.01];
                    \draw [fill] (0,0) circle [radius=0.01];
                    \draw [fill] (0.05,0) circle [radius=0.01];
                            \draw [fill] (0.1,0) circle [radius=0.01];
                    \draw[thick,-](-.1,0)--(-.1,-.2);
                    \draw[thick,-](.1,0)--(.1,-.2);
                    \node at (0,0) [above] {$\Theta$};

                    \begin{scope}[shift={(.5,0)}]
                    \draw [fill] (-0,0) circle [radius=0.05];
                    %\draw [fill] (0.15,0) circle [radius=0.01];
    %                \draw [fill] (0.2,0) circle [radius=0.01];
    %                \draw [fill] (0.1,0) circle [radius=0.01];
                    \draw [fill] (0.3,0) circle [radius=0.01];
                   % \node at(-.3,0)[above]{$j$};
%                    \node at(.3,0)[above]{$\alpha$};
                    \draw[thick,->](0,0)--(0,-.2);
                    \draw[thick,<-](0.3,0)--(0.3,-.2);
                    \node at (0.15,0) [above] {$A$};
    %                \node at (0.3,0) [right] {$I$};
                    \end{scope}

                \begin{scope}[shift={(1.3,0)}]
                \draw[draw=black](.9,-.2)rectangle(-.2,-.5);\node at(.35,-.35){$t$};
                 \draw [fill] (0,0) circle [radius=0.05];
                    %\draw [fill] (0.15,0) circle [radius=0.01];
    %                \draw [fill] (0.2,0) circle [radius=0.01];
    %                \draw [fill] (0.1,0) circle [radius=0.01];
                    \draw [fill] (0.3,0) circle [radius=0.01];
                    \draw[thick,->](0,0)--(0,-.2);
                    \draw[thick,<-](0.5,0)--(0.5,-.2);
                    \node at (0.25,0) [above] {${A}\backslash B$};
                    %\node at (0,0) [left] {$I'$};
    %                \node at (0.3,0) [right] {$J'$};
                \end{scope}
      \end{scope}

                %\draw [fill] (0,0) circle [radius=0.05];
    %                \draw [fill] (0.15,0) circle [radius=0.01];
    %                \draw [fill] (0.2,0) circle [radius=0.01];
    %                \draw [fill] (0.1,0) circle [radius=0.01];
    %                        \draw [fill] (0.3,0) circle [radius=0.05];
    %                \draw[dashed,-](0,-.05)--(0,-.15)--(0.3,-.15)--(0.3,-.05);
    %                \node at (0.15,-.15) [below] {$J$};

                    \begin{scope}[shift={(2.6,0)}]
                    \draw [fill] (-.1,0) circle [radius=0.01];
                    \draw [fill] (-0.05,0) circle [radius=0.01];
                    \draw [fill] (0,0) circle [radius=0.01];
                    \draw [fill] (0.05,0) circle [radius=0.01];
                            \draw [fill] (0.1,0) circle [radius=0.01];
                    \draw[thick,-](-.1,0)--(-.1,-.2);
                    \draw[thick,-](.1,0)--(.1,-.2);
                    \node at (0,0) [above] {$\Gamma$};

                    \end{scope}
        \end{scope}

        \node at (0,-.25){$\Rightarrow$};

        \begin{scope}[shift={(.65,0)}]
        \draw[draw=black](.9,-.2)rectangle(-.2,-.5);\node at(.35,-.35){$s$};

          \draw[thick,->](0.8,0) to  [out=90,in=180] (1.05,.25) to  [out=0,in=90](1.3,0);
          \draw[thick,-](.5,0) to  [out=90,in=180] (1.05,.5) to  [out=0,in=90](1.6,0);
          %\draw[thick,dotted,-](0.3,0) to  [out=-90,in=180] (1.3,-.75) to  [out=0,in=-90](2.3,0);
          %\draw[thick,dotted,-](0,0) to  [out=-90,in=180] (1.3,-1) to  [out=0,in=-90](2.6,0);

                \draw [fill] (-.1,0) circle [radius=0.01];
                \draw [fill] (-0.05,0) circle [radius=0.01];
                \draw [fill] (0,0) circle [radius=0.01];
                \draw [fill] (0.05,0) circle [radius=0.01];
                        \draw [fill] (0.1,0) circle [radius=0.01];
                \draw[thick,-](-.1,0)--(-.1,-.2);
                \draw[thick,-](.1,0)--(.1,-.2);
                \node at (0,0) [above] {$\Theta$};

                \begin{scope}[shift={(.5,0)}]
                %draw [fill] (0,0) circle [radius=0.05];
                %\draw [fill] (0.15,0) circle [radius=0.01];
%                \draw [fill] (0.2,0) circle [radius=0.01];
%                \draw [fill] (0.1,0) circle [radius=0.01];
                %\draw [fill] (0.3,0) circle [radius=0.01];
                %\node at(.15,0)[above]{$B$};
                \draw[thick,->](0,0)--(0,-.2);
                \draw[thick,<-](0.3,0)--(0.3,-.2);
             %   \node at (0,0) [left] {$J$};
%                \node at (0.3,0) [right] {$I$};
                \end{scope}

            \begin{scope}[shift={(1.3,0)}]
            \draw[draw=black](.9,-.2)rectangle(-.2,-.5);\node at(.35,-.35){$t$};
             \draw [fill] (.3,0) circle [radius=0.05];
                %\draw [fill] (0.15,0) circle [radius=0.01];
%                \draw [fill] (0.2,0) circle [radius=0.01];
%                \draw [fill] (0.1,0) circle [radius=0.01];
               % \draw [fill] (0.3,0) circle [radius=0.01];
                \draw[thick,->](0,0)--(0,-.2);
                \draw[thick,<-](0.5,0)--(0.5,-.2);
                \node at (.4,0) [above] {$B$};
%                \node at (0.3,0) [right] {$J'$};
            \end{scope}
            %\draw [fill] (0,0) circle [radius=0.05];
%                \draw [fill] (0.15,0) circle [radius=0.01];
%                \draw [fill] (0.2,0) circle [radius=0.01];
%                \draw [fill] (0.1,0) circle [radius=0.01];
%                        \draw [fill] (0.3,0) circle [radius=0.05];
%                \draw[dashed,-](0,-.05)--(0,-.15)--(0.3,-.15)--(0.3,-.05);
%                \node at (0.15,-.15) [below] {$J$};

                \begin{scope}[shift={(2.1,0)}]
                \draw [fill] (-.1,0) circle [radius=0.01];
                \draw [fill] (-0.05,0) circle [radius=0.01];
                \draw [fill] (0,0) circle [radius=0.01];
                \draw [fill] (0.05,0) circle [radius=0.01];
                        \draw [fill] (0.1,0) circle [radius=0.01];
                \draw[thick,-](-.1,0)--(-.1,-.2);
                \draw[thick,-](.1,0)--(.1,-.2);
                \node at (0,0) [above] {$\Gamma$};

                \end{scope}
        \end{scope}
        \end{tikzpicture}
        \caption{$(\backslash{\rm{E}})$ rule}
        %\label{Cut_geometrically}
        \end{subfigure}
        \caption{Slash elimination rules}
        \label{Slash elimination rules}
   \end{figure}

Figure \ref{Slash elimination rules} shows the geometric meaning of the above admissible rules. The most instructive case is when the contexts $\Gamma,\Theta$ are empty and the terms/graphs $t,s$ are single edges. Then the two rules amount to concatenating edges in different order.

\subsection{Grammars}\label{extended grammars section}
There are two versions of extended tensor type calculus, with and without Lambek restriction. Accordingly, there are two possible versions of corresponding grammars. We will say loosely {\it extended tensor grammar} without specifying the version.

Obviously, {\it extended tensor signatures} and {\it extended tensor grammars} are defined exactly the same as tensor signatures and tensor grammars with the adjective ``tensor'' replaced with ``extended tensor'' and the title {\bf TTC} replaced with {\bf ETTC} (respectively, ${\bf ETTC}_-$).

Important remarks are the following.

Firstly, in the extended tensor setting we may restrict to {\it lexicalized grammars}, i.e. those with lexicalized axioms. Indeed, a non-lexicalized axiom must have representation
$$
t\cdot\delta^\alpha_\beta\vdash \Gamma, A^{I\beta I'}_{J\alpha J'} \mbox{ or }t\cdot\delta^\alpha_\beta\vdash \Gamma, A^{I\beta I'}_L,
B_{J\alpha J'}^K,
$$
and it can be replaced with, respectively,
$$
t\vdash \Gamma, \nabla^\alpha_\beta A^{I\beta I'}_{J\alpha J'} \mbox{ or }t\vdash \Gamma, \nabla^\alpha_\beta (A^{I\beta I'}_L\wp B_{J\alpha J'}^K),
$$
 which has fewer non-lexical (Kronecker delta) components.
The old and the new grammars will be strictly equivalent thanks to the  $(\nabla)$, $(\nabla)^{-1}$,
$(\wp)$ and $(\wp)^{-1}$ rules. (Theoretically, there is a possibility that there is no $t$ in the axiom above, i.e. the  whole term to the left of the turnstile  is a Kronecker delta. Then the above argument does not work; but this seems a very artificial situation.)

Secondly, thanks to Lemma \ref{extended deduction theorem}, when dealing with lexicalized extended tensor grammars, we can replace derivations from non-logical axioms with derivations in pure cut-free {\bf ETTC}, similarly to the case of tensor grammars. This is very convenient for proof search and for theoretical analysis.

We now proceed to proving formally that Lambek calculus and Lambek grammars translate to {\bf ETTC} and extended tensor grammars as {\it conservative} fragments.

\section{Embedding Lambek calculus}\label{Lambek section}
\subsection{Lambek calculus}
Given a  set $P$ of atomic types, the set $Tp_{\backslash,/,\bullet}(P)$ of  {\it Lambek types} over $P$,   is defined by induction.
\begin{itemize}
\item Any $A\in P$ is in $Tp_{\backslash,/,\bullet}(P)$;
\item if $A,B\in Tp_{\backslash,/,\bullet}(P)$, then $(A\backslash B),(A/ B), (A\bullet B)\in Tp_{\backslash,/,\bullet}(P)$.
\end{itemize}

A {\it Lambek sequent} is an expression of the form $\Gamma\vdash A$, where the {\it context} $\Gamma$ is a finite {\it sequence} of Lambek  types, and $A$ is a Lambek type.

The Lambek sequent above satisfies the {\it Lambek restriction} if $\Gamma\not=\emptyset$.

 Lambek sequents are derived from the following  inference rules of {\it Lambek calculus} ({\bf LC }) \cite{Lambek}:
$${A\vdash A}~(\mbox{Id}),\quad\frac{\Gamma\vdash A\quad \Theta_1, A,\Theta_2\vdash B}{\Theta_1,\Gamma,\Theta_2\vdash B}~(\rm{Cut}),$$
$$\frac{\Gamma, A\vdash B}{\Gamma \vdash B/A}~(/\rm{R}), \quad\frac{\Gamma\vdash A\quad \Theta_1, B,\Theta_2\vdash C}{\Theta_1,B/A,\Gamma,\Theta_2\vdash C}~(/\rm{L}),$$
$$\frac{A,\Gamma\vdash B}{\Gamma \vdash A\backslash B}~(\backslash \rm{R}), \quad\frac{\Gamma\vdash A\quad \Theta_1, B,\Theta_2\vdash C}{\Theta_1,\Gamma,A\backslash B ,\Theta_2\vdash C}~(\backslash\rm{L}),$$
$$ \quad\frac{\Theta_1\vdash A\quad \Theta_2\vdash B}{\Theta_1,\Theta_2\vdash A\bullet B}~(\bullet\rm{R}),\frac{\Theta_1,A,B,\Theta_2\vdash C}{\Theta_1,A\bullet B,\Theta_2\vdash C}~(\bullet \rm{L}).$$

{\it Lambek calculus with  Lambek restriction} is obtained  by adding to the
($\backslash\rm{L}$) and ($/\rm{L}$) rules above the requirement  $\Gamma\not=\emptyset$. We will denote this system as ${\bf LC}_-$.

It is easy to see that all sequents derivable in ${\bf LC}_-$ satisfy  Lambek restriction.

(Our notation and terminology is somewhat unconventional. Usually, it is the version  {\it with Lambek restriction} that is taken as default.)

It is well known that Lambek calculus is cut-free: any sequent derivable in {\bf LC} (${\bf LC}_-$) is derivable without the Cut rule \cite{Lambek}.

\subsection{Translating Lambek calculus}\label{translating LC}
Given a set of atomic types $P$,
we assign to every element $p\in P$  valency $v(p)=(1,1)$ and treat it as an atomic tensor type symbol.
Then formulas (\ref{translating lambek types}) define by induction an embedding
$$tr:Tp_{\backslash,/,\bullet}(P)\to Tp_{\otimes,\nabla}(P)$$
of Lambek types to extended tensor type symbols,
where any Lambek type is translated to an extended tensor type symbol of valency (1,1), i.e., to a Lambek style type symbol.
We will abuse notation and  denote a Lambek type and its extended tensor translation by the same expression.

 If $A$ is a Lambek type, then we will call its image $A$ in $Tp_{\otimes,\nabla}$ a {\it positive Lambek type symbol}, and the dual $\overline A\in Tp_{\otimes,\nabla}$ of its image, a {\it negative Lambek type symbol}.

  We are going to translate {\bf LC} derivations to derivable extended tensor typing judgements. For that purpose we introduce some more terminology.

Let
\be\label{type sequence in a lambek cycle}
A_{(1)},\ldots,A_{(n)}
\ee
 be a cyclically ordered sequence of  Lambek style type symbols.

\begin{figure}
\begin{subfigure}{1\textwidth}
\centering
 \begin{tikzpicture}[xscale=1.5,yscale=1.5]
    \begin{scope}

    \begin{scope}%[xscale=-1]
    \begin{scope}%[shift={(-3.2,0)}]
        %\draw  (0,0) circle [radius=0.05];
        \draw[thick,->](0,0) to  [out=-90,in=180] (1.6,-.5) to  [out=0,in=-90](3.2,0);
        \draw [fill] (0,0) circle [radius=0.05];
       \node[right] at(.05,0){$\overline{A_{n}}$};
       \draw[thick,<-](.7,0) to  [out=-90,in=180] (.8,-.25) to  [out=0,in=-90](.9,0);
       \draw [fill] (.9,0) circle [radius=0.05];
       \node[right] at(1,0){$\ldots$};
        \node[right] at(1.85,0){$\overline{A_{1}}$};
       \draw[thick,<-](1.6,0) to  [out=-90,in=180] (1.7,-.25) to  [out=0,in=-90](1.8,0);
       \draw [fill] (1.8,0) circle [radius=0.05];
       \node[right] at(2.75,0){${A}$};
       \draw[thick,<-](2.5,0) to  [out=-90,in=180] (2.6,-.25) to  [out=0,in=-90](2.7,0);
       \draw [fill] (2.7,0) circle [radius=0.05];
       \end{scope}
       \end{scope}
%        %\node at(.75,.35){a};
        \node[right] at(3.4,0){$\Leftrightarrow$};
        \node[right] at(3.9,0){$A_{1},\ldots A_n\vdash A$};

        \end{scope}
       % \node[right] at (1.5,.25){Boundary: $X$, $|X|=6$, $X_l=\{2,5,6\}$.};
%         \node[right] at (1.5,.0) {Edges:  $[6,xy,1]$, $[2,ba,4]$, $[5,\epsilon,3]$.};

        \end{tikzpicture}
        \caption{Lambek cycle of a Lambek sequent}
        \label{cycle_example}
\end{subfigure}
\centering
\begin{subfigure}{1\textwidth}
\centering
 \begin{tikzpicture}[xscale=1.5,yscale=1.5]
    \begin{scope}
        %\draw  (0,0) circle [radius=0.05];
        \draw[thick,->](0,0) to  [out=-90,in=180] (1.1,-.5) to  [out=0,in=-90](2.2,0);
        \draw [fill] (0,0) circle [radius=0.05];
        \node[above] at(0,0.){$\beta$};
        \node[above] at(2.2,0.){$\alpha$};
       \node[right] at(.1,0){$\overline{A}$};
       \draw[thick,<-](.6,0) to  [out=-90,in=180] (.75,-.25) to  [out=0,in=-90](.9,0);
       \draw [fill] (.9,0) circle [radius=0.05];
       \node[ right] at(.95,0){$\overline{\Gamma}$};
        \node[right] at(1.75,0){${B}$};
        \node[above] at(1.7,0.){$i$};
       \node[above] at(.6,0.){$j$};

       \draw[thick,<-](1.4,0) to  [out=-90,in=180] (1.55,-.25) to  [out=0,in=-90](1.7,0);
       \draw [fill] (1.7,0) circle [radius=0.05];

%        %\node at(.75,.35){a};
        \node[right] at(2.4,0){$=$};

        \begin{scope}[shift={(2.9,0)}]
        \draw[thick,->](0,0) to  [out=-90,in=180] (1.1,-.5) to  [out=0,in=-90](2.2,0);
        \draw [fill] (0,0) circle [radius=0.05];
        \node[above] at(1.4,0.){$\alpha$};
        \node[above] at(1.7,0.){$\beta$};
       \node[right] at(.1,0){$\overline{\Gamma}$};
       \draw[thick,<-](.6,0) to  [out=-90,in=180] (.75,-.25) to  [out=0,in=-90](.9,0);
       \draw [fill] (.9,0) circle [radius=0.05];
       \node[ right] at(.95,0){$B$};
       \node[above] at(.9,0.){$i$};
       \node[above] at(2.2,0.){$j$};
        \node[right] at(1.75,0){$\overline{A}$};
       \draw[thick,<-](1.4,0) to  [out=-90,in=180] (1.55,-.25) to  [out=0,in=-90](1.7,0);
       \draw [fill] (1.7,0) circle [radius=0.05];
       \node[right] at(2.4,0){$\Rightarrow$};
             \end{scope}

                \node at(2.4,-1){$\cdots$ $(\nabla)$ rule $\cdots$};

        \begin{scope}[shift={(1.1,-1.5)}]
        \draw[thick,->](0,0) to  [out=-90,in=180] (1.3,-.7) to  [out=0,in=-90](2.6,0);
        \draw [fill] (0,0) circle [radius=0.05];
        \node[right] at(.1,0){$\overline{\Gamma}$};
       \draw[thick,<-](.6,0) to  [out=-90,in=180] (.75,-.25) to  [out=0,in=-90](.9,0);
       \node[above] at(.9,0.){$i$};
        \node[above] at(2.6,0.){$j$};
       \draw [fill] (.9,0) circle [radius=0.05];
       \node[ right] at(1,-.2){$\underbrace{\nabla^\alpha_\beta(B^i_\alpha\wp\overline{A}^\beta_j)}_{=B/A}$};

        \end{scope}

        \end{scope}

        \end{tikzpicture}
        \caption{Translating the ($/R$) rule}
        \label{slash on the right intro picture}
\end{subfigure}
%%%%%%%%%%%%%%%%%%%%%%%%%%%%%%%%%%%%%%%%%%%%%%%%%%%%%%%%%%%%%%%%%%%%%%%%%%%%%%%%%%%%%%%
%%%%%%%%%%%%%%%%%%%%%%%%%%%%%%%%%%%%%%%%%%%%%%%%%%%%%%%%%%%%%%%%%%%%%%%%%%%%%%%%%%%%%%%%%%
\begin{subfigure}{1\textwidth}
\centering
 \begin{tikzpicture}[xscale=1.5,yscale=1.5]
    \begin{scope}
        %\draw  (0,0) circle [radius=0.05];
        \draw[thick,->](0,0) to  [out=-90,in=180] (.65,-.5) to  [out=0,in=-90](1.3,0);
        \draw [fill] (0,0) circle [radius=0.05];

        \node[above] at(1.3,0.){$\beta$};
       \node[right] at(.1,0){$\overline{\Gamma}$};
       \draw[thick,<-](.5,0) to  [out=-90,in=180] (.65,-.25) to  [out=0,in=-90](.8,0);
       \draw [fill] (.8,0) circle [radius=0.05];
       \node[ right] at(.85,0){${A}$};

       \node[above] at(.8,0.){$j$};

    \begin{scope}[shift={(1.7,0)}]
       \draw[thick,->](0,0) to  [out=-90,in=180] (1.5,-.5) to  [out=0,in=-90](3,0);
        \draw [fill] (0,0) circle [radius=0.05];
        \node[above] at(.9,0.){$\alpha$};

       \node[right] at(.1,0){$\overline{\Theta_2}$};
       \draw[thick,<-](.6,0) to  [out=-90,in=180] (.75,-.25) to  [out=0,in=-90](.9,0);
       \draw [fill] (.9,0) circle [radius=0.05];
       \node[ right] at(.95,0){$\overline{B}$};
        \node[right] at(1.75,0){$\overline{\Theta_1}$};

        \draw[thick,<-](2.2,0) to  [out=-90,in=180] (2.35,-.25) to  [out=0,in=-90](2.5,0);
       \draw [fill] (2.5,0) circle [radius=0.05];
       \node[ right] at(2.55,0){${C}$};

        \node[above] at(1.4,0.){$i$};

       \draw[thick,<-](1.4,0) to  [out=-90,in=180] (1.55,-.25) to  [out=0,in=-90](1.7,0);
       \draw [fill] (1.7,0) circle [radius=0.05];
       \node[right] at(3.2,0){$=$};
    \end{scope}

\begin{scope}[shift={(0,-1)}]
\draw[thick,->](0,0) to  [out=-90,in=180] (.65,-.5) to  [out=0,in=-90](1.3,0);
        \draw [fill] (0,0) circle [radius=0.05];

        \node[above] at(1.3,0.){$\beta$};
       \node[right] at(.1,0){$\overline{\Gamma}$};
       \draw[thick,<-](.5,0) to  [out=-90,in=180] (.65,-.25) to  [out=0,in=-90](.8,0);
       \draw [fill] (.8,0) circle [radius=0.05];
       \node[ right] at(.85,0){${A}$};

       \node[above] at(.8,0.){$j$};

    \begin{scope}[shift={(1.7,0)}]
       \draw[thick,->](0,0) to  [out=-90,in=180] (1.5,-.5) to  [out=0,in=-90](3,0);
        \draw [fill] (0,0) circle [radius=0.05];
        \node[above] at(0,0.){$\alpha$};

       \node[right] at(.1,0){$\overline{B}$};
       \draw[thick,<-](.6,0) to  [out=-90,in=180] (.75,-.25) to  [out=0,in=-90](.9,0);
       \draw [fill] (.9,0) circle [radius=0.05];
       \node[ right] at(.95,0){$\overline{\Theta_1}$};
        \node[right] at(1.75,0){${C}$};

        \draw[thick,<-](2.2,0) to  [out=-90,in=180] (2.35,-.25) to  [out=0,in=-90](2.5,0);
       \draw [fill] (2.5,0) circle [radius=0.05];
       \node[ right] at(2.55,0){$\overline{\Theta_2}$};

        \node[above] at(.6,0.){$i$};

       \draw[thick,<-](1.4,0) to  [out=-90,in=180] (1.55,-.25) to  [out=0,in=-90](1.7,0);
       \draw [fill] (1.7,0) circle [radius=0.05];
       \node[right] at(3.2,0){$\Rightarrow$};
    \end{scope}

    \node at(2.35,-1){$\cdots$ $(\triangle)$ rule $\cdots$};

    \begin{scope}[shift={(0,-1.7)}]
\draw[thick,->](0,0) to  [out=-90,in=180] (2.5,-.75) to  [out=0,in=-90](5,0);
        \draw [fill] (0,0) circle [radius=0.05];

       \node[right] at(.1,0){$\overline{\Gamma}$};
       \draw[thick,<-](.5,0) to  [out=-90,in=180] (.65,-.25) to  [out=0,in=-90](.8,0);
       \draw [fill] (.8,0) circle [radius=0.05];
       \node[ right] at(.9,-.2){$\underbrace{\triangle_\alpha^\beta(A^j_\beta\otimes \overline{B}^\alpha_i)}_{=\overline{B/A}}$};

       \node[above] at(.8,0.){$j$};

    \begin{scope}[shift={(2,0)}]

       \draw[thick,<-](.6,0) to  [out=-90,in=180] (.75,-.25) to  [out=0,in=-90](.9,0);
       \draw [fill] (.9,0) circle [radius=0.05];
       \node[ right] at(.95,0){$\overline{\Theta_1}$};
        \node[right] at(1.75,0){${C}$};

        \draw[thick,<-](2.2,0) to  [out=-90,in=180] (2.35,-.25) to  [out=0,in=-90](2.5,0);
       \draw [fill] (2.5,0) circle [radius=0.05];
       \node[ right] at(2.55,0){$\overline{\Theta_2}$};

        \node[above] at(.6,0.){$i$};

       \draw[thick,<-](1.4,0) to  [out=-90,in=180] (1.55,-.25) to  [out=0,in=-90](1.7,0);
       \draw [fill] (1.7,0) circle [radius=0.05];
       %\node[right] at(3.2,0){$\Rightarrow$};
    \end{scope}

\end{scope}
\end{scope}
\end{scope}

\end{tikzpicture}
        \caption{Translating the ($/L$) rule}
        \label{slash on the left intro picture}
\end{subfigure}
\caption{Embedding Lambek calculus}
\end{figure}
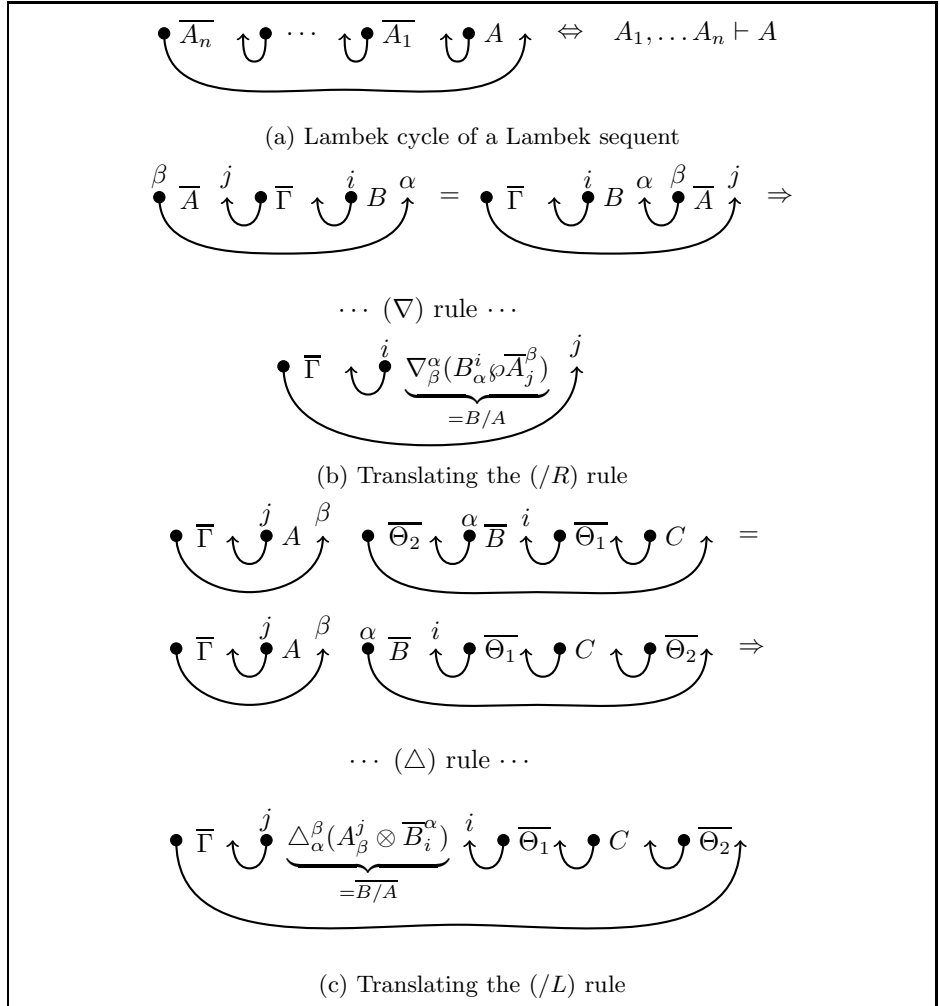

The {\it Lambek cycle} of (\ref{type sequence in a lambek cycle}) is the typing judgement
\be\label{Lambek cycle}
t\vdash (A_{(1)})^{i_1}_{j_1},\ldots,(A_{(n)})^{i_n}_{j_n},\mbox{ where } t=\delta^{j_n}_{i_1}\cdot\delta^{j_1}_{i_2}\cdots\delta^{j_{n-1}}_{i_n}.
\ee
The {\it Lambek cycle of a Lambek sequent} $A_{(1)},\ldots, A_{(n)}\vdash A$ is the  Lambek cycle of the cyclically ordered sequent
$\overline{A_{(n)}},\ldots,\overline{A_{(1)}},  A$,
  see  Figure \ref{cycle_example}.
  Lambek cycles, essentially, are proof-nets for cyclic logic  \cite{LamarcheRetore}.

 The meaning of Lambek cycles is summarized in the following proposition, which must be obvious from Figure \ref{cycle_example}.
  \nb\label{meaning of Lambek cycles}
 If we have an extended tensor signature $\Xi$    over the terminal alphabet $T$, words $w_1,\ldots,w_n\in T^*$
  and  Lambek style types $A_{(1)},\ldots,A_{(n)},A$  of $\Xi$, such that the typing judgements
  $$[w_1]^j_i\vdash (A_{(1)})^i_j,\ldots,[w_n]_i^j\vdash (A_{(n)})_j^i$$ and the Lambek cycle of $\overline{A_{(n)}},\ldots,\overline{A_{(1)}},A$ are derivable in $\Xi$, then
  the typing judgement
  $$[w_1\ldots w_n]^j_i\vdash A^i_j$$ is derivable in $\Xi$. $\Box$
 \nbe

\bl\label{LC2Tensor}
For  a Lambek sequent $\Gamma\vdash A$  derivable in {\bf LC}, respectively ${\bf LC}_-$, its  Lambek cycle
 is derivable in {\bf ETTC}, respectively ${\bf ETTC}_-$.
\el
{\bf Proof}
Induction on a cut-free derivation $\pi$ of $\Gamma\vdash A$.
Single-premise rules of {\bf LC} translate as the $(\wp)$ rule followed by the $(\nabla)$ rule of {\bf ETTC}, and
two-premise rules of {\bf LC} translate as the $(\otimes)$ rule followed by the $(\triangle)$ rule of {\bf ETTC}.

The cases of ($/\rm{R}$) and ($/\rm{L}$)  are shown in Figures \ref{slash on the right intro picture}, \ref{slash on the left intro picture} respectively. Remaining four rules are treated very similarly. $\Box$

\subsubsection{Conservativity of translation}
Let us say that an extended tensor  type (formula) $X$ is {\it locally connected}  if
all its subformulas of the form $A\wp B$ occur in  subformulas of the form $\nabla^\alpha_\beta(A\wp B)$, where $\alpha\in \mathit{FInd}(A)$,
$\beta\in \mathit{FInd}(B)$ or vice versa.
An extended tensor typing judgement $\sigma$ is locally connected, if
all type symbols occurring in $\sigma$ are locally connected.

Note that  all positive and negative Lambek types are locally connected.
\bp\label{locally connected}
Let $\sigma$ be a locally connected typing judgement derivable in {\bf ETTC} (or  ${\bf ETTC}_-$). Then any $\otimes$-link in $\sigma$  splits it into derivable typing judgements.
\ep
{\bf Proof} Induction on the number of connectives and binding operators occurring in $\sigma$.

We consider the last rule in a cut-free derivation of $\sigma$ and if this is not the $(\nabla)$-rule introducing a link $F$ of the form $F=\nabla^\alpha_\beta(A\wp B)$, we use the induction hypothesis applied to the premises of the rule.

Now consider the problematic case. Let $\sigma$ be obtained from a typing judgement $\sigma'$ by the $(\nabla)$-rule introducing the link $F$ as above.
Then we apply the $(\wp)^{-1}$ rule to $\sigma'$ and get a locally connected typing judgement $\sigma''$ with links $A$, $B$, with fewer connectives and binding operators than $\sigma$.

If $X$ is a $\otimes$-link in $\sigma$ then $X$ is a link in $\sigma''$ as well, and by the induction hypothesis it splits $\sigma''$ into derivable typing judgements $\tau$ and $\rho$. But the links $A,B$  cannot get apart after splitting, because there is a subterm (i.e. edge) $\delta_\beta^\alpha$ connecting them in $\sigma'$, hence in $\sigma''$. So assume that both $A,B$ occur, say,  in $\tau$. We apply to $\tau$ the $(\wp)$ rule and then the $(\nabla)$ rule introducing $F$ and obtain a derivable typing judgement $\tau'$. Then $X$ splits $\sigma$ into $\rho$ and $\tau'$. $\Box$
\bl\label{conservativity for lambek}
Let $\sigma$ be a typing judgement derivable in {\bf ETTC}, respectively ${\bf ETTC}_-$, such that it has exactly one occurrence of a positive Lambek type symbol, and all other links in $\sigma$ are negative Lambek type symbols. Then $\sigma$ is a translation of a sequent derivable in {\bf LC}, respectively ${\bf LC}_-$.
\el
{\bf Proof} Induction on the number of connectives and binding operators occurring in $\sigma$.

If $\sigma$ has a $\nabla$ link $F$ then it necessarily has the form $F=\nabla^\alpha_\beta(A_\alpha\wp B^\beta)$ and corresponds to one of three possible forms: $\overline{X\bullet Y}$ (if both $A$, $B$ are negative),
$X/Y$ (if $A$ is positive and $B$ negative), or $X\backslash Y$ (if $A$ is negative and $B$ positive).

 We apply the $\nabla^{-1}$ rule followed by the $(\wp^{-1})$ rule to obtain a derivable typing judgement $\sigma'$ with links $A$, $B$ and fewer connectives and binding operators. By the induction hypothesis, $\sigma'$ is the translation of some {\bf LC} derivable  sequent $\Gamma$. Then it is immediate that $\sigma$ is the translation of the sequent $\Gamma'$ obtained from $\Gamma$ by one of the three single premise rules in accordance with one of the three possible forms of $F$.

 If $\sigma$ has no $\nabla$ links and is not an axiom, consider the last rule in a cut-free derivation of $\sigma$. This must be the $(\triangle)$ rule introducing a link $F$ of the form $F=\triangle^\alpha_\beta(A_\alpha\otimes B^\beta)$, which also has three possible forms (we skip listing them). Then the premise $\sigma'$ is locally connected and has a unique $\otimes$-link. By Proposition \ref{locally connected} above, it is split by the $\otimes$-link into two typing judgements $\rho$ and $\tau$, with $\rho$ getting the link $A$ and $\tau$ getting $B$.

  The induction hypothesis applies to $\rho$ and to $\tau$, and they are translations of {\bf LC} derivable sequents $\Gamma$ and $\Theta$.
  It is easy to check that in all three possible cases the typing judgement $\sigma$ is  the translation of the sequent obtained from $\Gamma$ and $\Theta$ by one of the three two-premise rules of {\bf LC} in accordance with one of the three possible forms of $F$. $\Box$

\subsection{Translating Lambek grammars}
\subsubsection{Lambek grammars}
A {\it Lambek grammar} $G$ is a tuple $G=(P,T,Lex, S)$,
where
\begin{itemize}
  \item $P$ is a set of atomic types;
  \item $T$ is a finite alphabet of {\it terminal symbols};
  \item $Lex$, the {\it Lexicon}, is a finite set of expressions of the form $a:A$, called {\it axioms}, where $a\in T$, $A\in Tp_{\backslash,/,\bullet}$;
      \item $S$, the {\it sentence type}, is an element of $P$.
\end{itemize}

A word $w\in T^*$ is in the {\it type} $A$ {\it generated} by $G$, if there exist axioms
$$a_{1}:A_{1},\ldots, a_{n}:A_{n}\in Lex$$
such that  the sequent $A_{1},\ldots,A_{n}\vdash A$ is derivable in ${\bf LC}_-$, and $$w=a_{1}\ldots a_{n}.$$

The {\it language of} $G$ is the set of all words in the sentence type generated by $G$.

A {Lambek grammar without Lambek restriction} is defined as above, replacing ${\bf LC}_-$  with ${\bf LC}$.

\subsubsection{Translation}
Let a Lambek Grammar $G=(P,T,Lex, S)$ be given.

We treat elements of $P$ as atomic tensor  type symbols and translate all Lambek types to extended tensor type symbols using formulas (\ref{translating lambek types}).

The axioms of $G$ are immediately translated to extended tensor typing judgements.
An axiom of the form $w:A\in Lex$ translates as $[w]^j_i\vdash A_j^i$, where we denote a Lambek type and its translation the same.

Then the set  $\{[w]^i_j\vdash A_i^j|~w:A\in Lex\}$ of these translation defines an extended tensor signature $\widetilde{Lex}$ over $P$ and $T$.

We define the extended tensor grammar translating $G$ as $\widetilde{G}=(\widetilde{Lex},S)$.
It is easy to see that the translation is an {\it embedding}.
\bt\label{Lambek grammars embedding theorem}
For any Lambek type $A\in Tp_{\backslash,/,\bullet}$, a word $w\in T^*$
is in the Lambek type $A$ generated by $G$ iff the extended tensor typing judgement $[w^i_j]\vdash A_i^j$ is derivable in the signature $\widetilde{Lex}$.
\et
{\bf Proof} Lemma \ref{extended deduction theorem} (``Deduction theorem'' for lexicalized signatures),   Corollary \ref{conservativity  for lambek} (on conservativity of translation from {\bf LC} to {\bf ETTC}) and Proposition \ref{meaning of Lambek cycles}  (on the meaning of Lambek cycles). $\Box$
\bc
The languages generated by $G$ and $\tilde{G}$ coincide. $\Box$
\ec

\section{Conclusion}
In this work we defined the cut-free system of  tensor type calculus ({\bf TTC}), which is a system of typed proof terms, generalizing proof-nets, based on classical linear logic (Section \ref{tensor types section}).

We defined tensor grammars based on this system and proved that string abstract categorial grammars embed into tensor grammars  up to $\beta\eta$-equivalence of the surface level terms (Theorem \ref{conservativity for ACG}).

 Next we  defined the cut-free system {\bf ETTC} of extended tensor type calculus, containing {\bf TTC} as a subsystem, which allows binding operators on types (Section \ref{extended things section}). We defined extended tensor grammars, which allow emulating operations of Lambek calculus as well. We proved that Lambek calculus grammars embed into extended tensor grammars as conservative fragments (Theorem \ref{Lambek grammars embedding theorem}).

 We leave proving similar embedding results for displacement grammars and hybrid type logical grammars for future work. We also hope to be able to answer in the future  the important open question of how do {\bf ETTC} and extended type grammars relate to grammatical formalisms based on first order multiplicative linear logic.

 \end{document}